\documentclass{amsart}
\usepackage{amsmath,amsthm}
\usepackage{amsfonts,amssymb}

\usepackage{enumerate}

\newtheorem{thm}{Theorem}[section]
\newtheorem{cor}[thm]{Corollary}
\newtheorem{lem}[thm]{Lemma}
\newtheorem{prop}[thm]{Proposition}
\newtheorem{exam}[thm]{Example}

\newtheorem{defn}[thm]{Definition}

\theoremstyle{remark}
\newtheorem{rem}{Remark}[section]

 \def\tr{{\triangle}}
 \def\ll{\lesssim}

\def\sph{\mathbb{S}^{d-1}}
\def\f{\frac}
\def\vi{\varphi}

 \def\a{{\alpha}}
 \def\b{{\beta}}
 \def\g{{\gamma}}
 \def\k{{\kappa}}
 \def\t{{\theta}}
 \def\l{{\lambda}}
 \def\d{{\delta}}
 \def\o{{\omega}}
 \def\s{{\sigma}}
 \def\va{\varepsilon}
 \def\la{{\langle}}
 \def\ra{{\rangle}}
 \def\ve{{\varepsilon}}

 \def\CD{{\mathcal D}}
 
 \def\CH{{\mathcal H}}

 \def\CV{{\mathcal V}}
 
 \def\BB{{\mathbb B}}

 \def\II{{\mathbb I}}
 \def\JJ{{\mathbb J}}
 \def\NN{{\mathbb N}}

 \def\RR{{\mathbb R}}
  \def\SS{{\mathbb S}}
 \def\ZZ{{\mathbb Z}}
        \def\grad{\operatorname{grad}}

\newcommand{\wt}{\widetilde}
\newcommand{\wh}{\widehat}

\begin{document}

\title[modulus of smoothness and approximation on sphere and ball]
{Moduli of Smoothness and Approximation on the Unit Sphere and the Unit Ball}
\author{Feng Dai}
\address{Department of Mathematical and Statistical Sciences\\
University of Alberta\\, Edmonton, Alberta T6G 2G1, Canada.}
\email{dfeng@math.ualberta.ca}
\author{Yuan Xu}
\address{Department of Mathematics\\ University of Oregon\\
    Eugene, Oregon 97403-1222.}\email{yuan@math.uoregon.edu}

\date{\today}
\keywords{modulus of smoothness, $K$-functional, best approximation, sphere,
Euler angles, ball}
\subjclass[2000]{42B15, 41A17, 41A63}
\thanks{Work of the  first  author  was  supported  in part
 by the NSERC grant of
 Canada. Work of the second author was supported in part by the National
Science Foundation under Grant DMS-0604056}

\begin{abstract}
A new modulus of smoothness based on the Euler angles is
introduced on the unit sphere and is shown to satisfy all the
usual characteristic properties of moduli of smoothness, including
direct and inverse theorem for the best approximation by
polynomials and its equivalence to a $K$-functional, defined via
partial derivatives in Euler angles. The set of results on the
moduli on the sphere serves as a basis for defining new moduli of
smoothness and their corresponding $K$-functionals on the unit
ball, which are used to characterize the best approximation by
polynomials on the ball.
\end{abstract}

\maketitle

 \tableofcontents

\section{Introduction}
\setcounter{equation}{0}

A central problem for approximation theory is to characterize the best
approximation by polynomials via moduli of smoothness, or via $K$-functionals.
In this paper we consider the setting of best
approximation by polynomials on the unit sphere and the unit ball
$$
\SS^{d-1} = \{x: \|x\| =1\} \qquad \hbox{and} \qquad \BB^d=\{x: \|x\| \le 1\}
$$
of $\RR^d$, where $\|x\|$ denotes the usual Euclidean norm.

\subsection{Approximation on the unit sphere}
On the unit sphere, we consider the best approximation by
polynomials in the space $L^p(\SS^d)$, $1 \le p < \infty$, or
$C(\SS^{d-1})$ for $p = \infty$, with norm denoted by $\| \cdot
\|_{p} := \|\cdot\|_{L^p(\SS^{d-1})}$, $1 \le p \le \infty$, in
the Part 1. Let $\Pi_n^d$ denote the space of polynomials of total
degree $n$ in $d$ variables and $\Pi_n(\SS^{d-1}) := \Pi_n^d
\vert_{\SS^{d-1}}$,  the space of spherical polynomials, or
equivalently polynomials restricted on the sphere. In the
following we shall write $\Pi_n^d$ for $\Pi_n^d(\SS^{d-1})$
whenever it causes no confusing. The quantity of best
approximation is defined by
\begin{equation}\label{bestapp}
    E_n (f)_{p} : = \inf_{ g  \in \Pi_{n-1}^d} \| f - g \|_{p}, \qquad 1 \le p \le \infty.
\end{equation}

The first modulus of smoothness that characterizes $E_n(f)_p$ on the
sphere is defined via the spherical means
\begin{equation} \label{s-means}
 S_\theta f(x) := \frac{1}{\s_{d-1} (\sin \theta)^{d-2}}
     \int_{\langle x, y\rangle = \cos\theta} f(y) d\s_{x,\theta}(y),
\end{equation}
where $d \s_{x,\theta}$ is the Lebesgue measure on  $\{ y\in
\SS^{d-1}: \langle x, y\rangle = \cos\theta\}$ and $\s_{d}  =
2 \pi^{d/2}/\Gamma(d/2)$ (\cite[p. 216]{BBP}, \cite[p. 475]{P}, \cite[p. 288]{Rud}).
For $r > 0$ and $t> 0$, this modulus of smoothness is defined by
\begin{equation}\label{omega*}
  \o_r^*(f, t)_p := \sup_{|\t| \le t} \| (I - S_\t)^{r/2} f\|_p,
\end{equation} where $(I-S_\t)^{ r/2}$ is defined in terms of
infinite series when $ r/2$ is not an integer
\mbox{(\cite[p.183]{WL})}. After earlier studies by several
authors (see, for example, \cite{BBP, LN, P}), Rustamov \cite[p.
315]{Rus} finally established, for $1<p<\infty$,  both direct and
inverse theorems for the polynomial best approximation, as well as
the equivalence of $\o_r^*(f,t)_p$ to the $K$-functional
\begin{equation}\label{K-func*}
  K_r^*(f,t)_p := \inf_g \left\{ \|f-g\|_p +
  t^r \|(- \Delta_0)^{r/2}g\|_p \right\},
\end{equation}
 where $\Delta_0$, given in \eqref{Laplace-Betrami} below, is the
Laplace-Beltrami operator on the sphere and the infimum is taken
over all $g$ for which $(-\Delta_0)^{r/2} g \in L^p$. The proofs
of these results for the full range of $1\leq p\leq \infty$ can be
found in \cite[p.195-216]{WL}.
 The study
of $\o_r^*(f, t)_p$ and  $K_r^*(f, t)_p$ relies heavily on the
fact that both $(I-S_\t)^{r/2}$ and $(-\Delta_0)^{r/2}$ are
multiplier operators of Fourier series in spherical harmonics.
This approach has been extended in \cite[p. 15]{X05a} to the
setting of weighted space $L^p(\SS^{d-1}, h_\k^2)$, where $h_\k$
is a weight function invariant under a finite reflection group.

The second modulus of smoothness on the sphere is defined via rotation,
\begin{equation}\label{T_Q}
   T_Q f(x) := f(Qx), \qquad Q \in SO(d),
\end{equation}
where $SO(d)$ denotes the group of rotations on $\RR^d$, or the group of
orthogonal matrices of determinant 1 in $\RR^d$. For $t > 0$, define
$$
  O_t := \left \{Q \in SO(d):  \max_{x \in \SS^{d-1}} d(x,Qx) \le t \right \},
$$
where $d(x,y): = \arccos \la x,y \ra$ denotes the geodesic distance on
$\SS^{d-1}$. For $r > 0$ and $t > 0$ define
\begin{equation}\label{omegaD}
\wt \o_r(f,t)_p : = \sup_{Q \in O_t} \| \tr_Q^r f\|_p,
    \qquad\hbox{where}\quad  \tr_Q^r: = (I - T_Q)^r.
\end{equation}
For $r =1$ and $p=1$ this modulus of smoothness was introduced and
used by Calder\`{o}n, Weiss and Zygmund (\cite{CWZ}) and further
studied in \cite{KW}. For other spaces, including
$L^p(\SS^{d-1})$, $p>0$, these moduli were introduced and
investigated in \cite{Di1}. The direct and weak converse theorems
for $L^p(\SS^{d-1})$, $1 \le p \le \infty$ were given in \cite[p.
23]{Di2} and \cite[p. 197]{Di1}, respectively. An easier proof of
the direct result applicable to a more general class of spaces was
given in \cite{DD1}. In \cite[(9.1)]{DDH} it is shown that $\wt
\o_r(f,t)_p$ is equivalent to $\o^\ast_r(f,t)_p$ when $1 < p <
\infty$, whereas the equivalence fails for $p =1$ and $p =\infty$
\cite{Di3}.

The modulus of smoothness defined via multipliers, such as $\o^*_r (f, t)_p$,
allows an easy access to a neat theory, but it is hard to compute and more
difficult to follow because of its dissimilarity to the traditional modulus of
smoothness defined  via differences of function evaluations. The modulus
$\wt \o_r(f,t)_p$, on the other hand, is closer to the traditional form, as
distance on the sphere is measured by geodesic distance; in fact, the
authors in \cite{CWZ} considered it the most natural definition on the sphere.
However, the supremum over $O_t$ makes it difficult, if at all possible, to
compute even for simple functions.

In the present paper we shall introduce another modulus of smoothness,
denoted by $\o_r(f,t)_p$, that can be expressed as forward differences
in Euler angles. More precisely,
\begin{equation}
  \o_r(f,t)_p :=\sup_{|\t| \le t} \max_{1 \le i <j \le d} \|\tr_{i,j,t}^r f\|_p,
\end{equation}
where $\tr_{i,j,t}^r$ denotes the $r$-th forward difference in the
Euler angle $\t_{i,j}$. These angles can be described (see next
section and \cite[Chapt. 9]{Vi}) by rotations on two dimensional
planes, so that the new modulus of smoothness can be defined
through a collection of two dimensional forwarded differences,
which are well understood and can be easily computed. Examples of
functions will be given in Section 4, where $\o_r(f,t)_p$ is
computed whereas we do not see how to compute $\wt \o_r(f,t)_p$ or
$\o^*_r(f,t)_p$. Both direct and inverse theorems will be
established in terms of this new modulus of smoothness. We will
also define a new $K$-functional, using the derivatives with
respect to the Euler angles, and show that it is equivalent to
$\o_r(f,t)_p$. Comparing to the other moduli of smoothness, we
shall prove that $\o_r (f,t)_p$ is bounded by both $\wt
\o_r(f,t)_p$ and $\o_r^*(f,t)_p$, and is equivalent to them for $1
< p < \infty$ and $r = 1, 2$. The strength of the new modulus of
smoothness lies in its computability. We will give examples to
show how the asymptotic order of $\o_r(f,t)_p$ can be determined.

By taking the norm in a weighted $L^p$ space, we can also define our
modulus of smoothness for a doubling weight. Best approximation in a
weighted space with respect to a doubling weight was first investigated,
on an interval, in \cite{MT1,MT2}. It was studied in \cite{Dai2} on the sphere
in terms of weighted version of the modulus of smoothness $\o_r(f,t)_p$.
We shall show that the results in \cite{Dai2} can be established using our
new modulus of smoothness.

\subsection{Approximation on the unit ball}
For the unit ball $\BB^d$, a modulus of smoothness, denoted by
$\o_r^*(f,t)_{p,\mu}$, as it is in the spirit of $\o_r^*(f,t)_p$
in \eqref{omega*},  is introduced in \cite[p. 503]{X05b} for
$L^p(\BB^d,W_\mu)$, where
\begin{equation}\label{weight}
W_\mu(x) := (1-\|x\|^2)^{\mu-1/2}, \qquad
   \mu \ge 0,
\end{equation}
in terms of the generalized translation operator of the orthogonal series,
given explicitly in \cite[p. 500]{X05b}, and used to characterize the best
approximation by polynomials. There were also earlier results in \cite[p. 164]{Ra}
of direct theorem given in terms of $\sup_{\|h\| \le t} |f(x+ h) - f(x)|$, which,
however, does not take into account the boundary of $\BB^d$ and, hence,
does not have a matching inverse theorem. At the moment, the modulus
$\o_r^*(f,t)_{p,\mu}$ is the only one that gives both direct and inverse
theorem for $d > 1$.

There is a close connection between analysis on the sphere and on
the unit ball $\BB^d = \{x: \|x\| \le 1\}$ (\cite[Sect. 4]{X05a}, \cite[Sect. 4]{X06}).
Our results on the sphere can be used to define a new modulus of smoothness
and a new $K$-functional on the unit ball with weight function $W_\mu$
for $\mu$ being a half integer, and all results on the sphere can be carried
over to the weighted approximation on the ball. It is worth to mention that
our results appear to be new even in the case of $d=1$, which corresponds
to the extensively studied case of best approximation in $L^p([-1,1],W_\mu)$,
and our new modulus of smoothness takes the form
$$
  \o_r(f,t)_{p,\mu} = \sup_{|\t| \le t} \left(\int_{\BB^2}
    \left |\overrightarrow{\tr}_\t^r f(x \cos (\cdot) + y \sin (\cdot)) \right|^p
      (1-x^2-y^2)^{\mu-1} dxdy\right)^{p}.
$$
There are several well-studied modulus of smoothness in this setting of one
variable. Among others the most established one is due to Ditzian and Totik
\cite[p. 11]{Di-To}, defined in the unweighted case ($\mu =1/2$ in \eqref{weight}) by
\begin{equation}\label{1-2omegaDT}
  \o_\varphi^r (f,t)_p\equiv \wh \o_r (f,t)_p : = \sup_{0 < \t \le t}
       \|\wh \tr^r_{\t \varphi} f\|_{L^p[-1,1]}, \quad
  1 \le p \le \infty,
\end{equation}
where $\varphi(x) = \sqrt{1-x^2}$, $\wh\tr^r_h f(x)$ is the $r$-th
central difference which equals $0$ when $x\pm \f h2\notin [-1,1]$
(see \cite[p. 11]{Di-To} or \eqref{modu-DT} below for details) and
we have dropped $\mu =1/2$ in the notation of the norm. It turns
out that
$$
   \o_r(f,t)_{p,1/2} \le c\, \o_\varphi^r(f,t)_p, \qquad 1 \le p \le \infty.
$$
The comparison between the last two moduli of smoothness and the
$K$-functional equivalent to $\o_\varphi^r(f,t)_{p,\mu}$ suggests yet another
pair of modulus of smoothness and $K$-functional on the unit ball, which
can be regarded as a natural extension to those defined by Ditzian-Totik.
In the unweighted case, it is defined by
\begin{equation}\label{1-2omegaDTball}
    \o_\varphi^r (f,t)_{L^p(\BB^d)}\equiv \wh\o_r (f,t)_p: =
     \sup_{|\t| \le t} \left \{ \max_{1\le i< j \le d}
     \|\tr^r_{i,j,\t} f\|_p, \max_{1 \le i \le d}
      \|\wh \tr^r_{\t \varphi e_i} f\|_p  \right\},
\end{equation}
where $\varphi(x) = \sqrt{1-\|x\|^2}$,  $e_i$ denotes the $i$-th
coordinator vector, and $\|\cdot\|_p$ is the  $L^p$ norm computed
with respect to the Lebesgue measure on $\BB^d$ (see the
definition in Section 7.3 for details). A corresponding
$K$-functional can also be defined. We are able to prove the main
results normally associated with moduli of smoothness and
$K$-functionals for this pair.

These new moduli of smoothness and $K$-functionals provide, we believe,
a satisfactory solution for the problem of characterizing the best
approximation on the unit ball, and new tools for gauging the smoothness
of functions on the unit ball. Our computation examples give the asymptotic
order of the moduli of smoothness for several functions, which are not
intuitively evident, and will be hard to distinguish without the new moduli
of smoothness.

\subsection{Organization of the paper}
The paper is naturally divided into three parts. Part 1 deals with approximation
on the sphere, whereas Part 2 deals with approximation on the ball. Part 3 contains
examples of functions for which the asymptotic order of new moduli of smoothness
and best approximation by polynomials are determined.

Throughout this paper we denote by $c, c_1, c_2, ...$ generic constants that may
depend on fixed parameters, whose value may vary from line to line. We write
$A \ll B$ if $A \le c B$ and $A \sim B$ if $A \ll B$ and $B \ll A$.

\part{Approximation on the Unit Sphere}

This part is organized as follows. The new modulus of smoothness and
$K$-functional on the sphere are defined and studied in Section 2, their
equivalence and the characterization of best approximation in terms of them
are proved in Section 3. Finally, in Section 4, we discuss the weighted
approximation with respect to a doubling weight.

\section{A New Modulus of Smoothness and $K$-functional}
\setcounter{equation}{0}

\subsection{Euler angles and Laplace-Betrami operators}
In the case of $d=3$, the Euler angles are often used to describe
motions in the Euclidean space and are well-known to physicists
and people working in computer graphics. We shall need the
definition for $d \ge 3$, for which we follow \cite[p.438]{Vi}.

Let $e_1, \cdots, e_d$ denote the standard orthogonal basis in $\RR^d$.
For $1\leq i \neq j \leq d$ and $t \in \RR$, we denote by $Q_{i,j,t}$ a rotation
by the angle $t$ in the $(x_i, x_j)$-plane, oriented such that the rotation
from the vector $e_i$ to the vector $e_j$ is assumed to be positive. For example,
the action of the rotation $Q_{1,2, t}\in SO(d)$ is given  by
\begin{align}\label{EulerAngle1}
Q_{1,2,t}(x_1,\ldots, x_d) = & (x_1 \cos t -x_2 \sin t, x_1 \sin t + x_2 \cos t, x_3,
   \cdots, x_d) \\
  = & (s \cos (\phi + t), s \sin (\phi +t), x_3, \cdots, x_d), \notag
\end{align}
where $(x_1,x_2) = s (\cos \phi, \sin \phi)$, and other $Q_{i,j,t}$ are defined likewise.
It is known \cite[p. 438]{Vi} that every rotation $Q\in SO(d)$ can be presented in the form
\begin{equation}\label{EulerAngle2}
   Q = Q_{d-1} Q_{d-2}\cdots Q_1, \quad \hbox{where} \quad
   Q_k = Q_{1, 2, \t_1^k}Q_{2, 3, \t_2^k}\cdots Q_{k, k+1, \t_k^k}
\end{equation}
for some $\t_1^k \in [0, 2\pi)$ and $\t_2^k, \cdots, \t_k^k  \in [0,\pi)$, and the
representation \eqref{EulerAngle2} is unique  for almost all elements $Q$ of
$SO(d)$. The numbers
$$
    \t_j^k, \qquad 1 \le j \le k, \quad 1 \le k \le d-1
$$
are called the Euler angles of the rotation $Q$. There are a total $d(d-1)/2$ Euler
angles, which agrees with the dimension of $SO(d)$.

We note that an Euler angle comes from a two dimensional rotation. The
following simple fact is useful in our development below.

\begin{lem} \label{lem:EulerAngle}
Suppose that $1\leq i\neq j\leq d$ and $x, y \in \sph$ differ only
at their $i$-th and $j$-th components. Then $y = Q_{i,j,t}\, x$
with the angle $t$ satisfying
$$
   \cos t = (x_i y_i + x_j y_j) / s^2  \quad \hbox{and} \quad
      t \sim \|x-y\| /s \quad \hbox{with}\quad s: = \sqrt{x_i^2+ x_j^2}.
$$
\end{lem}

\begin{proof}
Since $x$ and $y$ differ at exactly two components, they differ by a two
dimensional rotation. Moreover, as $x_i^2+x_j^2 = y_i^2+y_j^2$, the
formula for $\cos t$ is the classical formula for the angle between two
vectors in $\RR^2$. We also have
$$
  t^2 \sim 4 \sin^2 \tfrac{t}2 = 2 (1-\cos t) = \|(x_i,x_j)- (y_i,y_j)\|^2 /s^2
       = \|x-y\|^2 /s^2,
$$
where the first $\|\cdot \|$ is the Euclidean norm of $\RR^2$ and the
second one is of $\RR^d$.
\end{proof}

To each $Q \in SO(d)$ corresponds an operator $L(Q)$ in the space
$L^2(\SS^{d-1})$, defined by $L(Q) f(x):= f(Q^{-1} x)$ for $x \in
\SS^{d-1}$. Since $L(Q_1Q_2) = L(Q_1)L(Q_2)$, it is a group
representation of $SO(d)$. In terms of Euler angles, the
infinitesimal operator of $L(Q_{i,j,t})$ has the form
\begin{equation}\label{partial_ij}
 D_{i,j}: = \frac{\partial}{\partial t} \left[ L(Q_{i,j,t})\right] \Big \vert_{t =0}
       = x_j \frac{\partial}{\partial x_i}   - x_i \frac{\partial}{\partial x_j}, \qquad
        1 \le i < j \le d,
\end{equation}
where the second equation follows from \eqref{EulerAngle1}. For
more details, see \cite[Chapt. IX]{Vi}. In particular, it is easy
to verify that, taking $(i,j) = (1,2)$ as an example,
\begin{equation}\label{partial_ij2}
 D_{1,2}^r f(x) = \left(-\frac{\partial}{\partial \phi}\right)^r
  f(s \cos \phi, s \sin \phi, x_3,\ldots,x_d),
 \end{equation}
where  $(x_1,x_2) = (s \cos \phi, s \sin \phi)$. It turns out that these operators
are closely related to the Laplace-Betrami operator on the sphere.
Let $\Delta = \frac{\partial^2}{\partial x_1^2}+ \ldots + \frac{\partial^2}{\partial x_d^2}$
be the usual Laplace operator. The Laplace-Betrami operator is defined
by the relation
\begin{equation} \label{Laplace-Betrami}
    \Delta_0 f(x) = \Delta \left[ f \left(\frac{y}{\|y\|}\right) \right](x),\   \   \    \ x \in\SS^{d-1},
\end{equation}
where the Laplace operator $\Delta$ acts on the  variables $y$. The explicit
formula of $\Delta_0f(x)$, $x \in \SS^{d-1}$, is often given in terms of differential operators in
spherical coordinates as in \cite[p. 494]{Vi}. It turns out, however, that
it also satisfies a decomposition,
\begin{equation} \label{Laplace-Betrami2}
  \Delta_0 = \sum_{1 \le i < j \le d} D_{i,j}^2.
\end{equation}
When applied to a function, the right hand side of this decomposition is defined for all
$x \in \RR^d$, but the above identity holds for those $x$ restricted to $\SS^{d-1}$.
The point  is that each operator $D_{i,j}$ in this decomposition  commutes with the
Laplace-Beltrami operator $\Delta_0$.  This decomposition must be
classical but we are not aware of a convenient reference. It can
be easily verified, however. In fact, a straightforward computation of
$ \sum_{i=1}^d \partial_i^2 f( y/\|y\|)$ shows, by \eqref{Laplace-Betrami}, that
\begin{equation} \label{Laplace-Betrami3}
  \Delta_0 = \Delta - \sum_{i=1}^d \sum_{j=1}^d x_i x_j \partial_i \partial_j
      - (d-1) \sum_{i=1}^d x_i \partial_i,
\end{equation}
where $\partial_i$ denotes the $i$-th partial derivative, and the right hand side
of \eqref{Laplace-Betrami2} gives the same formula as an other straightforward
computation shows.


\subsection{New modulus of smoothness and $K$-functional}
For each $Q \in SO(d)$, we have defined $\tr^r_Q f = (I - T_Q)^r f$ in
\eqref{omegaD}. For the rotations $Q_{i,j,\t}$ in the Euler angles, we
shall denote
$$
        \tr_{i,j,\t}^r : =\tr_{Q_{i,j,\t}}^r, \qquad 1 \le i \ne j \le d
$$
for convenience. Since $Q_{i,j,\t} = Q_{j,i,-\t}$, we have $ \tr_{i,j,\t}^r
=  \tr_{j,i,-\t}^r$. Let $\overrightarrow{\tr}_\t^r$ denote the forward
difference operator acting on $f: \RR \mapsto \RR$, defined by
$\overrightarrow{\tr}_\t f(t) := f(t+\t) - f(t)$ and $\overrightarrow{\tr}_\t^r
:= \overrightarrow{\tr}_\t^{r-1} \overrightarrow{\tr}_\t$; then
$$
 \overrightarrow{\tr}_\t^r f(t) = \sum_{j=0}^r (-1)^j \binom{r}{j} f(t + \t j).
$$
Because of \eqref{EulerAngle1}, it follows that $\tr_{i,j,\t}^r$
can be expressed in the forward difference. For instance,  take
$(i,j) =(1,2)$ as  example,
\begin{equation} \label{Delta_ij}
\tr_{1,2, \t}^r f(x) =  \overrightarrow{\tr}_\t^r
   f \left (x_1 \cos(\cdot)- x_2 \sin (\cdot), x_1 \sin(\cdot) + x_2 \cos (\cdot),
        x_3,\ldots,x_d \right),
\end{equation}
where $\overrightarrow{\tr}_\t^r$ is acted on the variable
$(\cdot)$, and is evaluated at $t=0$.

\begin{defn} \label{def:modulus}
For $r \in \NN$, $t > 0$, and $f \in L^p(\SS^{d-1})$, $1 \le p < \infty$, or
$f \in C(\SS^{d-1})$ for $p = \infty$, define
\begin{equation} \label{eq:modulus}
 \o_r (f,t)_p : = \sup_{ |\t| \le t} \max_{1 \le i < j \le d}
          \left \|\Delta_{{i,j,\t}}^r f \right \|_p.
\end{equation}
For $r =1$ we write $\o (f,t)_p := \o_1(f,t)_p$.
\end{defn}

Let us remark that this modulus of smoothness is not rotationally
invariant, that is, if we define $f_Q(x) = f(Qx)$, then $\o_r(f_Q,
t)_p$ is in general different from $\o_r (f,t)_p$, whereas both
$\o_r^*(f,t)_p$ and $\wt \o_r(f,t)_p$ are rotationally invariant.
Moreover, $\wt \o_r(f,t)_p$ does not depend on the choice of the
the orthogonal basis of $\RR^d$, whereas, on the face of it,  the
new modulus $\o_r(f,t)_p$ relies on the standard basis
$e_1,\ldots, e_d$ of $\RR^d$, but do not, we note, on the order of
$e_1,\ldots, e_d$. As will be shown later in Subsection 3.4, the
three moduli are nevertheless closely related (see Corollary 3.11
for details) and the result below shows that $\o_r(f,t)_p$ is
smaller than $\wt \o_r(f,t)_p$.

Recall $O_t = \{Q \in SO(d):  \max_{x \in \SS^{d-1}} d(x,Qx) \le t\}$. For $d=2$,
there is only one Euler angle $\t$; a rotation $Q \in SO(2)$ belongs to $O_t$ if
and only if its Euler angle $\t$ satisfies $|\sin \t| \le \sin t$. Hence, for $d =2$,
$\o_r(f,t)_p$ agrees with $\wt \o_r(f,t)_p$ in \eqref{omegaD}. This, however,
does not extend to $d \ge 3$. A rotation $Q \in SO(d)$ belonging to $O_t$ may
not be easily characterized by its Euler angles. For example, for $d =3$, the
rotation $Q_{1,2, 2\pi -\t} Q_{2,3, t} Q_{1,2,\t}$ is in $O_t$ for all $\t \in (0,2\pi)$,
as can be easily seen from \eqref{EulerAngle1}. On the other hand, for
$x \in \sph$, a quick computation shows that
$$
  \la Q_{i,j,\t} x, x \ra = (x_i^2 +x_j^2)\cos \t + \sum_{k \ne i,j} x_k^2
    =  \cos \t + \sum_{k \ne i,j} x_k^2(1- \cos \t) \ge \cos \t.
$$
Consequently, since $\cos d (x,y) = \la x, y\ra$, we obtain
$$
  d(Q_{i,j,\t}x ,x)  = \arccos \la Q_{i,j,\t} x, x \ra \le \t
$$
which shows that $Q_{i,j,\t} \in O_t$ for $0<\t \le t$. As a
result, we immediately see that the following proposition holds.

\begin{prop} \label{omega-omegaD}
For $f \in L^p(\sph)$ if $1 \le p < \infty$ and $f \in C(\sph)$ if $p = \infty$,
$$
       \o_r(f,t)_p \le \wt \o_r(f,t)_p, \qquad 1 \le p \le \infty, \quad r \in \NN.
$$
\end{prop}

The main advantage of the new modulus of smoothness is that it
reduces to forward differences in Euler angles, which live on two
dimensional circles on the sphere, and many of its properties can
be deduced from the corresponding results for trigonometric
functions of one variable.

Our new $K$-functional is defined via the differential operators $D_{i,j}$
in \eqref{partial_ij}, which can be regarded as derivatives with respect to
the Euler angles.

\begin{defn}  \label{def:K-func}
For $r \in \NN_0$ and $t \ge 0$,
\begin{equation} \label{eq:K-func-sphere}
   K_r(f,t)_p : = \inf_{g \in C^r(\sph)} \left\{ \|f - g\|_p + t^r \max_{1 \le i<j \le d}
        \|D_{i,j}^r g\|_p\right\}.
\end{equation}
\end{defn}

One usually defines the $K$-functional in $L^p$ norm by taking the minimum
over a Sobolev space, such as
$$
    W_r^p : = \{g \in L^p(\sph): \|D_{i,j}^r g\|_p < \infty,
    \quad 1 \le i\neq j \le d \}.
$$
Since we will deal with several different $K$-functionals, it is more convenient
to take the minimum over $C^r(\SS^{d-1})$, the space of $r$-th continuous
differentiable functions, which however is no less general by the density of
$C^r(\SS^{d-1})$ in the Sobolev spaces. In Section 3, we will show that
$K_r(f,t)_p$ is equivalent to $\o_r(f,t)_p$.

\subsection{Properties of the modulus of smoothness}
We will need the following elementary lemma (see Lemma 3.8.9 and Lemma 3.6.1
in \cite{DX}), in which $d\sigma$ denotes the usual Lebesgue measure on
$\SS^{d-1}$ without normalization.

\begin{lem}
Let $d$ and $m$ be positive integers. If $m \ge 2$, then
\begin{equation}  \label{IntS-B}
\int_{\SS^{d+m-1}} f(y) d\s = \int_{\BB^d} (1-\|x\|^2)^{\frac{m-2}{2}}
 \left[\int_{\SS^{m-1}} f(x, \sqrt{1-\|x\|^2}\, \xi) d\s(\xi)\right] dx,
\end{equation}
whereas if $m =1$, then
\begin{equation}  \label{IntS-Bm=1}
  \int_{\SS^{d}} f(y) d\s = \int_{\BB^d} \left[ f(x, \sqrt{1-\|x\|^2})+
         f(x, - \sqrt{1-\|x\|^2})\right]  \frac{dx}{\sqrt{1-\|x\|^2}}.
\end{equation}
\end{lem}

Because of the maximum in the definition of $\o_r(f,t)_p$, it is
more convenient, and often more useful, to state the properties on
$\sup_{|\t|\le t} \|\Delta_{i,j,\t}^r f\|_p$, some of which are
collected in the lemma below.

\begin{lem} \label{lem:Delta_ij}
Let  $r \in \NN$ and let $f\in L^p(\sph)$ with $1\leq p< \infty$, or $f\in C(\sph)$
when $p=\infty$.

\begin{enumerate}[\rm(i)]
\item For any $\l >0$,  $t\in (0, 2\pi]$, and $1\leq i<j\leq d$, we
have
$$
\sup_{|\t|\leq \l t} \|\tr^r_{i, j,\t} f\|_p\leq  (\l+1)^r
        \sup_{|\t|\leq  t} \|\tr^r_{i, j,\t} f\|_p.
$$
\item For $1\leq i\neq j\leq d$ and $\t\in [-\pi, \pi]$,
$$
  \|\tr^r_{i, j,\t} f \|_p \le 2^r \|f\|_p \quad\hbox{and} \quad
      \|\tr^r_{i, j,\t} f\|_p \le c \,|\t|^r  \|D_{i,j}^r f\|_p.
$$
\item If $f\in \Pi^d_n$ and $1\leq i<j\leq d$,  then
$$
\|\tr_{i, j, n^{-1}}^r f \|_p \sim n^{-r}\|D_{i, j}^r f\|_p
$$
\item For $1\leq i<j\leq d$ and $t\in (0, 2\pi)$,
$$
\sup_{|\t|\leq t} \|\tr^r_{i, j,\t} f\|_p^p\sim \f 1t
    \int_0^t \|\tr_{i, j,\t}^r f\|_p^p\, d\t
$$
with $\|\cdot\|_p^p$ replaced by $\|\cdot\|_\infty$ when
$p=\infty$.
\end{enumerate}
\end{lem}

\begin{proof}
Clearly we only need to consider the case of $(i,j) = (1,2)$. For
$f$ defined on $\SS^{d-1}$ we set $g_{s, y}(\phi) := f (s \cos
\phi, s \sin \phi, \sqrt{1-s^2} y)$, where $y \in \SS^{d-3}$,
$s\in [0,1]$ and $\phi\in [0, 2\pi]$.

(i) For a positive integer $n$, the well known identity
$$
 \overrightarrow{\tr}_{n \t}^r g(t) = \sum_{\nu_1 =0}^{n-1} \cdots \sum_{\nu_r=0}^{n-1}
    \overrightarrow{\tr}_\t^r g(t + \nu_1 \t + \ldots + \nu_r \t)
$$
and the connection \eqref{Delta_ij} imply immediately the inequality
$$
\sup_{|\t|\leq n t} \|\tr^r_{i, j,\t} f\|_p\leq  n^r \sup_{|\t|\leq  t} \|\tr^r_{i, j,\t} f\|_p,
$$
from which (i) follows from the monotonicity of $\sup_{|\t|\leq t} \|\tr^r_{i, j,\t} f\|_p$
in $t$ and $n = \lfloor \l \rfloor$.

(ii) By \eqref{IntS-B}), for $d >3$,
\begin{align} \label{Int_m=2}
  \int_{\sph} f(y)\, d\s(y) & = \int_{\BB^2} \int_{\SS^{d-3}}
     f(x_1, x_2, \sqrt{1-\|x\|^2} y)\, d\s(y)\,  (1-\|x\|^2)^{\f{d-4}2} dx \\
    & = \int_0^1 s(1-s^2)^{\f{d-4}2}
    \int_{\mathbb{S}^{d-3}}\int_0^{2\pi} g_{s, y} (\phi)\, d\phi\, d\s(y)\, ds, \notag
\end{align}
where the second equality follows from changing variables $(x_1,x_2) = s(\cos \phi,
\sin \phi)$. Using \eqref{Delta_ij}, the identity \eqref{Int_m=2} implies immediately
\begin{equation} \label{norm_Delta_ij}
\|\tr_{1,2,t}^r f\|_p^p =\int_0^1 s(1-s^2)^{\f{d-4}2} \int_{\mathbb{S}^{d-3}}
\left [ \int_0^{2\pi} |\overrightarrow{\tr}_t^r g_{s, y}(\phi)|^p \,
d\phi\right]\, d\s(y) \, ds.
\end{equation}
In the case of $d =3$, the formula \eqref{Int_m=2} degenerated to a form in
which the integral over $\SS^{d-3}$ is replaced by a sum of two terms,
see \eqref{IntS-Bm=1}. By
\eqref{EulerAngle1} and \eqref{partial_ij}, it is easy to see that
\begin{equation}\label{tri}
 (- 1)^rD_{1,2}^r f(s\cos\phi, s\sin\phi, \sqrt{1-s^2}y)=g^{(r)}_{s,y}(\phi)
\end{equation}
and $g_{s,y}(\phi)$ is a $2\pi$-periodic function. Hence, the desired result
follows from the corresponding result for the trigonometric functions
on $\mathbb{T}$.

(iii) If $f\in \Pi^d_n$, then $g_{s,y}(\phi)$ is a trigonometric
polynomial of degree at most $n$ in $\phi$. The classical result
of Steckin (\cite{St}) shows that for a trigonometric polynomial
$T_n$ of degree at most $n$,
$$
\|T_n^{(r)} \|_{L^p(\mathbb{T})} \sim
   h^{-r}\|\overrightarrow{\tr}^r_{h} T_n\|_{L^p(\mathbb{T})},   \qquad
   0<h\leq \pi n^{-1}
$$
with the constant of equivalence depending only on $r$. Thus, (iii) follows
by \eqref{norm_Delta_ij} and \eqref{tri}.

(iv) This again follows from \eqref{norm_Delta_ij} and the corresponding
result for trigonometric function. Indeed, by 
\cite[p.191, Lemma 7.2]{PP}, we have for $0<t\leq 2\pi$,
$$
\sup_{|\t |\leq t} \left( \int_0^{2\pi} |\overrightarrow{\tr}_\t^r
     g_{s, y}(\phi)|^p \, d\phi \right)^{\f1p}  \sim \left(\f 1t  \int_0^t\int_0^{2\pi}
       |\overrightarrow{\tr}_\t^r g_{s,y}(\phi)|^p \, d\phi\, d\t \right )^{\f1p}
$$
with the usual modification when $p=\infty$, from which (iv) follows from
\eqref{norm_Delta_ij}.
\end{proof}

\begin{prop} \label{modulus}
The modulus of smoothness $\o_r(f,t)_p$ satisfies
\begin{enumerate}[\rm(1)]
\item For $s < r$, $\o_r (f,t)_p \le 2^{r-s} \o_s (f,t)_p$.
\item For $\l > 0$, $\o_r(f,\l t)_p \le (\l +1)^r \o_r(f,t)_p$.
\item For $0<t<\f12$ and every $m>r$,
$$
    \o_r(f,t)_p \le c_m t^r \int_t^1 \f {\o_{m}(f, u)_p}{u^{r+1}}\, du.
$$
\end{enumerate}
\end{prop}

\begin{proof}
The first property follows from the identity
$$
   (I - T)^r = (I -T)^s \sum_{k=0}^{r-s} \binom{r-s}{k} (-1)^k T^k
$$
and the triangle inequality. The second one follows immediately from
(ii) of Lemma \ref{lem:Delta_ij}. The third one is the Marchaud type
inequality and it follows, by \eqref{norm_Delta_ij}, from Marchaud
inequality for the trigonometric functions, in which the additional term
$t^r \|f\|_p$ that usually appears in the right hand can be removed
upon using $\inf_{c\in \RR} \|f-c\|_p \ll  \o_r(f,\pi)$ that follows from
Lemma \ref{KeyLemma} below.
\end{proof}

Recall the distance $d(x,y) : = \arccos \la x,y\ra$ on $\sph$. It follows that
\begin{equation} \label{distant}
\|x-y\| = \sqrt{2- 2 \cos d(x,y)} = 2 \sin \f {d(x,y)}{2} \sim d(x,y).
\end{equation}

\begin{lem} \label{lem:omega_infty}
For $x,y \in \sph$,
$$
         |f(x) - f(y)|  \le c \, \o(f,d(x,y))_\infty,
$$
where $c$ depends only on dimension.
\end{lem}

\begin{proof}
We may assume that $d(x,y) \le \d_d: = 1/(2 d^2)$. Otherwise we
can select an integer $m$ such that $d(x,y) \le m \d_d < 1$, then
$m$ is finite and we can select points $x=z_0, z_1, \cdots,
z_{m}=y$ on the great circle connecting $x$ and $y$ on $\sph$ such
that $d(z_i, z_{i+1})=\f {d(x, y)}{m} \leq \d_d$ for
$i=0,1,\cdots, m-1$, and then use triangle inequality. Since
$\|x\| =1$ implies that $|x_i | \ge 1/\sqrt{d}$ for at least one
$i$, we can assume without losing generality, as $\o_r(f,t)_p$ is
independent of the order of $e_1,\ldots, e_d$, that
$x_d=\max_{1\leq j\leq d} |x_j| \ge \f 1{\sqrt{d}}$.

For $1 \le j \le d-2$, let $u_j':=(x_1,\ldots,x_j,y_{j+1},\ldots,y_{d-1})$
and $v_j := \sqrt{1-\|u_j'\|^2}$, where by the choice of $x_d$ and $\d_d$,
\begin{align*}
  \|u_j'\|^2 = &  1-(x_{j+1}^2 - y_{j+1}^2) - \ldots - (x_{d-1}^2 - y_{d-1}^2) - x_d^2 \\
   & \le 1 - \tfrac{1}{\sqrt{d} }+ 2 (d-j-2) d(x,y) \le  1 - \tfrac{1}{\sqrt{d}} + \tfrac{1}{d} < 1.
\end{align*}
 We then define $u_0=y$, $u_j =
(u_j', v_j) \in \sph$ for $1 \le j \le d-2$, and  $u_{d-1} =x$. By definition,
$u_j$ and $u_{j-1}$ differ at exactly $j$-th and $d$-th elements, so that we
can write $u_{j-1} = Q_{j,d,t_j} u_j$, where the Euler angle $t_j$ satisfies,
by Lemma \ref{lem:EulerAngle},
$$
  t_j \sim \|u_{j-1} - u_j\| / s_j,  \quad\hbox{where}  \quad s_j ^2= x_j^2 + v_j^2.
$$
Our assumption shows that
\begin{align*}
s_j^2  \ge v_j^2 & = x_d^2 + (x_{j+1}^2-y_{j+1}^2)+\ldots +
 (x_{d-1}^2 - y_{d-1}^2) \\
      & \ge d^{-1} -  \|x-y\|^2 \ge \f 1{2d},
\end{align*}
and, on the other hand, by \eqref{distant},
\begin{align*}
   \|u_j - u_{j-1}\|^2 = |x_j -y_j|^2 +  \frac{(x_j^2-y_j^2)^2} {(v_{j-1} +v_j)^2} \le
        (1+ 8d) |x_j-y_j|^2 \ll d(x,y)^2.
\end{align*}
Together the last three displayed equations imply that $t_j \ll d(x,y)$. Hence,
\begin{align*}
   |f(x)-f(y)|& \le \sum_{j=1}^{d-1} |f(Q_{j,d, t_j} u_j) -
   f(u_j)|\\
   &        \le (d-1) \o (f, c \,d(x, y))_\infty \le c \, \o (f, d(x,y))_\infty,
\end{align*}
where the last step uses (2) of Proposition \ref{modulus}.
\end{proof}

\section{Approximation on the Unit Sphere}
\setcounter{equation}{0}

In this section we show that our new modulus of smoothness and $K$-functional are
equivalent and use them to establish direct and
inverse theorem for
\begin{equation} \label{eq:bestEn}
  E_n (f)_p := \inf_{g \in \Pi_{n-1}^d} \|f - g\|_p, \quad n=1,2,...,  \quad 1 \le p \le \infty.
\end{equation}

\subsection{Preliminaries}

Recall that $\Pi_n^d$ denote the spherical polynomials of degree at most $n$.
Let $\CH_n^d$ denote the space of spherical harmonics of degree $n$, which
are the restriction of homogeneous harmonic polynomials on $\sph$. It is well
known that the reproducing kernel of the space $\CH_n^d$ in $L^2(\sph)$ is
given by the zonal harmonic
\begin{equation}\label{zonal}
    Z_{n,d}(x,y) := \f{n + \l}{\l} C_n^\l(\la x,y\ra), \qquad \l = \f{d-2}2,
\end{equation}
where $C_n^\l$ is the Gegenbauer polynomial with index $\l$, normalized
by $C_n^\l(1) = \binom{n+2\l-1}{n}$.

Let $\eta$ be a $C^\infty$-function on $[0,\infty)$ with the properties
that $\eta(x)=1$ for $0\leq x\leq 1$ and $\eta(x)=0$ for $x\ge 2$. We
define
\begin{equation}\label{Vnf}
V_nf(x)=\int_{\sph} f(y) K_n(\la x, y\ra)\, d\s(y), \quad x\in\sph, \quad n=1,2,\cdots
\end{equation}
with
$$
K_n(t)=\sum_{k=0}^{2n} \eta\left(\f kn\right) \f{k+\l}{\l}C_k^{\l}(t),
    \quad  t\in[-1,1].
$$
By now it is well known (cf. \cite[p. 316]{Rus}) that $V_n(f)$ satisfies the following properties:

\begin{lem} \label{lem:Vnf}
Let $f \in L^p$ if $1 \le p < \infty$ and $f \in C(\sph)$ if $p = \infty$. Then
\begin{enumerate}[  \rm(1)]
 \item $V_n f \in \Pi^d_{2n}$ and $V_n f = f $ for $f \in \Pi^d_n$.
 \item For $n \in \NN$, $\|V_n f\|_p \le c \|f\|_p$
 \item For $n \in \NN$,
 $$
        \|f - V_n f\|_p \le c E_n(f)_p.
 $$
\end{enumerate}
\end{lem}

More importantly, the kernel is highly localized (\cite[p. 409]{BD}); that is, for any
positive integer $\ell$,  $K_n(t)$ satisfies
\begin{equation} \label{local_kernel}
     |K_n(\cos \t)| \le  c_\ell n^{d-1} (1 + n \t)^{-\ell} =: G_n(\t), \quad \t \in [0,\pi].
\end{equation}
The following lemma plays an essential role in our study below.

\begin{lem} \label{KeyLemma}
Suppose that  $f \in L^p(\sph)$ for $1\le p< \infty$, and $G_n(t)
\equiv G_{n,\ell}(t)$ is given by \eqref{local_kernel} with $\ell
> p + d$. Then
\begin{align*}
 \int_{\sph} \int_{\sph} |f(x)-f(y)|^p G_n( d(x, y) )
   d\s(x) \, d\s(y) \leq c \,\o(f, n^{-1})_p^p.
\end{align*}
\end{lem}

\begin{proof}
Let $E^{+}_j :=\{ x\in \sph: x_j \ge \f 1{\sqrt{d}}\}$ and $E^{-}_j :=\{ x\in \sph: x_j
 \le -\f 1{\sqrt{d}}\}$ for $1 \le j \le d$. Then $\sph =
 \bigcup_{j=1}^d (E_j^{+}\cup E_j^{-})$. Hence, it is enough to show that for
 each $1\leq k\leq d$,
\begin{align} \label{KeyLem1}
 & \int_{E_k^{\pm}} \int_{\sph} |f(x)-f(y)|^p |G_n (d(x, y))| \,d\s(y) \, d\s(x)
     \le c \, \o(f, n^{-1})_p^p.
\end{align}
By symmetry, it is enough to consider $E_d^{+}$.  For $0 < \d < \pi$ and $x \in \sph$,
let $c(x,\d)$ denote the spherical cap defined by
$$
   c(x,\d):= \{y \in \sph: d(x,y) \le \d\}.
$$
We choose $\d = 1/(100d)$ and split  the integral in
\eqref{KeyLem1} into two parts:
$$
\int_{E_d^{+}}\int_{c(x, \d)}\cdots\, d\s(y)\, d\s(x) +
  \int_{ E_d^{+}}\int_{\sph \setminus c(x, \d)}\cdots\, d\s(y)\, d\s(x)
  =: A+B.
$$
We first estimate the  integral $B$:
\begin{align} \label{KeyLem-B}
B &  =  \int_{E_d^{+}} \int_{\{y\in\sph:\  d(x, y)\ge \d\}}
  |f(x)-f(y)|^p |G_n (d(x, y))| \, d\s(y) \, d\s(x) \\
& \le  c\, n^{d-1-\ell} \int_{\sph}\int_{\sph} |f(x)-f(y)|^p
    d\s(x) \, d\s(y)\notag\\
 & =c\, n^{d-1-\ell} \int_{SO(d)}\int_{\sph} |f(x)-f(Qx)|^p  d\s(x) \, dQ,  \notag
\end{align}
where the last step uses the standard realization of $S^{d-1} = SO(d)/SO(d-1)$.
Using the decomposition of $Q$ in terms of Euler angles as in \eqref{EulerAngle2},
each $Q\in SO(d)$ can be decomposed as $Q=Q_1 Q_2\cdots Q_{d(d-1)/2}$
with $Q_k=Q_{i_k, j_k, t_k}$ for some $1\leq i_k<j_k\leq d$ and $t_k \in [0, 2\pi]$.
It then follows that
 \begin{align*}
  & \int_{\sph} |f(x)-f(Qx)|^p  d\s(x) \ll \int_{\sph} |f( Q_{d(d-1)/2}x)-f( x)|^p\,d\s(x)\\
 & \qquad + \sum_{k=1}^{\f{d(d-1)}2-1} \int_{\sph} |f(Q_k\cdots Q_{d(d-1)/2}x)-
     f(Q_{k+1}\cdots Q_{d(d-1)/2} x)|^p\, d\s(x)\\
 & \ll \max_{1\leq i<j\leq d}\sup_{0<\t \le 2\pi}\int_{\sph} |f(Q_{i,j,\t} x)-f(x)|^p\,
     d\s(x) \\
     & \ll \o(f, 2\pi)_p^p\ll n^p\o (f, n^{-1})_p^p,
 \end{align*}
which, together with (\ref{KeyLem-B}), gives the desired estimate
$B\leq c\, \o(f, n^{-1})_p^p$.

It remains to estimate the integral $A$. Setting $x = (x',x_d)$
with $x_d = \sqrt{1-\|x'\|^2}$,  we deduce from \eqref{IntS-Bm=1} that
\begin{align*}
A &=\int_{E_d^{+}} \int_{c(x,\d)} |f(x)-f(y) |^p G_n(d(x,y)) \, d\s(y) \, d\s(x)\\
   & = \int_{\|x'\| \le d^*} \int_{c(x, \d)}
       |f(x)-f(y) |^p G_n(d(x,y)) \, d\s(y) \, \frac{dx'}{\sqrt{1-\|x'\|^2}},
\end{align*}
where $d^* = \sqrt{1-d^{-1}}$.  Since $x_d \ge \f{1}{\sqrt{d}}$ it
follows that for any $y=(y', y_d)\in c(x, \d)$, $y_d\ge
x_d-|y_d-x_d|\ge x_d -d(x, y) \ge x_d-\d \ge \f 1 {2\sqrt{d}}$,
which further implies, by a simple computation, that
\begin{align*}
\|x'-y'&\|\leq \|x-y\|\leq \|x'-y'\|+|x_d-y_d|=\|x'-y'\|+ \f {
|\|x'\|^2-\|y'\|^2|}{x_d+y_d}\\
&\leq (1+2\sqrt{d} ) \|x'-y'\|.\end{align*} Consequently, setting
$g(x') := f(x',\sqrt{1-\|x'\|^2})$,   using  \eqref{IntS-Bm=1}
again  and observing  that $G_n(\t_1)\sim G_n(\t_2)$ whenever
$\t_1\sim \t_2$, we obtain
\begin{align*}
  A & \le c  \int_{\|x'\| \le d^*} \int_{\|x'-y'\|\le \d} |g(x')-g(y') |^p G_n(\|x'-y'\|) \, dy' \, dx'\\
  & = c \int_{\|x'\| \le d^*} \int_{\|u\|\le \d} |g(x')-g(u+ x')|^p G_n(\|u\|) \, du \, dx'.
\end{align*}
Let $b_0(u) :=0$ and $b_j(u): = u_1 e_1+\ldots + u_j e_j$, $1 \le j \le d-1$. Since
$$
    g(x')-g(x'+u) = \sum_{j=1}^{d-1} \left( g(x'+ b_{j-1}(u)) - g(x' + b_j(u)) \right),
$$
by triangle inequality it suffices to estimate, for $1 \le j \le d-1$,
\begin{align*}
 A_j:& = \int_{\|x'\|\le d^*} \int_{\|u\|\le \d}
  |g(x'+b_{j-1}(u) ) - g(x' + b_j(u))|^p
   G_n(\|u\|) du \,dx' \notag\\
   & \le \int_{\|x'\|\le d^* + \d} \int_{\|u\|\le \d}  |g(x') - g(x' + u_j e_j)|^p
                G_n(\|u\|) du \,dx',\notag\end{align*}
where the second line follows from a change of variables $x' +
b_{j-1}(u) \mapsto x'$. By symmetry, it suffices to consider $A_1$.

Observe that for $u_1\in\RR$ and  $u=(u_1, v)\in\RR^{d-2}$,
$$
G_n (\|u\|) =n^{d-1} (1+n\|u\|)^{-\ell}\leq H_n(|u_1|) n^{d-2}
(1+n\|v\|)^{-d+1},
$$
where $H_n(s) =n (1+ns)^{-\ell+d-1}$,  and we have used the
assumption $\ell>d-1$ in the last step. This implies
\begin{align} \label{3-6-eq}
A_1 & \ll  \int_{\|x'\|\le d^* + \d}  \int_{-\d}^\d |g(x') - g(x' + u_1e_1)|^p \\
  &  \qquad  \times \left[ \int_{\{v\in \RR^{d-2}: \|v\|\leq \sqrt{\d^2-|u_1|^2}\}}
       G_n(\|(u_1, v)\|) \,dv\right]\, du_1 \,dx'\notag\\
  &  \ll \int_{\|x'\|\le d^* + \d} \int_{-\d}^\d  |g(x') - g(x' + s e_1)|^p
                H_n(|s|)  ds\,dx'.\notag
\end{align}
Set $v_1(t, x') = - x_1 + x_1 \cos t -\sqrt{1- \|x'\|^2} \sin t$.
A straightforward calculation shows that
\begin{equation}\label{1-10}
   \f {1} {\pi \sqrt{d}}\leq \f{-v_1(t,x')} t \leq 2 \ \
    \text{and}\ \ \f 1{4\sqrt{d}}\leq -\f{\partial} {\partial t} v_1(t,x') \leq 2
\end{equation}
whenever $|t|\leq \sqrt{d}\d=\d^\ast$ and  $\|x'\|\leq d^\ast+\d$.
Thus, performing a change of variable $s=v_1(t, x')$ in
(\ref{3-6-eq}) yields
\begin{align*}
A_1 & \le c \int_{\|x'\|\le d^* + \d} \int_{-\d^\ast}^{\d^\ast}
             |g(x') - g(x' +v_1(x',t) e_1)|^p
           H_n(|v_1(x',t)|) \left |\f{\partial v_1(t,x')} {\partial t}\right | dt\,dx'\\
       & \le  c' \int_{\|x'\|\le \rho} \int_{-\d^\ast}^{\d^\ast}
                 |g(x') - g(x' +v_1(x',t) e_1)|^p H_n(|t|) dt\,dx',
\end{align*}
where $\rho: = \sqrt{1-(2d)^{-1}} \ge d^\ast +\d$, and we used
(\ref{1-10}) and the monotonicity of $H_n$ in the last step. Now
observe that for $x=(x',\sqrt{1-\|x'\|^2})$ with $\|x'\|\leq
\rho$,
$$
Q_{1,d,t} x = (x' +v_1(x',t) e_1, z_d), \qquad \forall t\in [-\d^\ast,\d^*],
$$
where, using the fact that $\sin t \le t \le 1/ (8 \sqrt{d})$,
$$
z_d= \sqrt{1-\|x' +v_1(x',t) e_1\|^2} = x_1 \sin t + \sqrt{1-\|x'\|^2} \cos t
   \ge  1/(4 \sqrt{d}) >0.
$$
Thus, using \eqref{IntS-Bm=1}, we deduce that
\begin{align*}
A_1 & \ll  \int_{-\d^\ast}^{\d^\ast}\int_{\{ x\in \SS^{d-1}: x_d\ge (2d)^{-1}\}}
           |f(x) - f(Q_{1,d,t}x)|^p\, d\s(x)  H_n(|t|) dt\\
        & \ll \int_{-\d^\ast}^{\d^\ast} \o(f, |t|)_p^p H_n(|t|) dt.
\end{align*}
Hence, by (2) of Proposition \ref{modulus} and the definition of $H_n$,
\begin{align*}
A_1 &  \ll \o(f,n^{-1})_p^p \int_{0}^{\d^\ast} (1+ n t)^p H_n(t)
dt =  \o(f,n^{-1})_p^p \, n \int_0^{\d^\ast}  (1+ n t)^{-\ell+ p+d-1}\, dt  \\
      & \ll \o(f,n^{-1})_p^p \int_0^\infty  (1+s)^{-\ell+ p+d-1} ds
         \ll \o(f,n^{-1})_p^p,
\end{align*}
since $\ell > p +d$. This completes the proof.
\end{proof}

The operator $V_n f$ plays an important role in our study. Our
next lemma shows that it commutes with $\tr_{i,j,t}$.

\begin{lem} \label{lem:triV=Vtri}
For $1 \le p \le \infty$ and $1 \le i\neq j\le d$,
$$
   \tr_{i,j,t}^r V_nf =   V_n\left(\tr_{i,j,t}^r f \right ).
$$
In particular, for $t > 0$,
$$
   \o_r (f- V_n f, t)_p \le c\, \o_r (f,t)_p.
$$
\end{lem}

\begin{proof}
Recall that $T_Q f(x) = f(Q x)$ for $Q\in SO(d)$. By the definition of
$V_n f$,
\begin{align*}
T_Q V_n f (x) = & \int_{\sph} f(y) K_n(\la Qx, y \ra) d\s(y) =
   \int_{\sph} f(y) K_n(\la x, Q^{-1}y \ra) d\s(y) \\
   = &  \int_{\sph} f(Qy) K_n(\la x, y \ra) d\s(y) = V_n (T_Q
   f)(x)
\end{align*}
by the rotation invariance of $d \s(y)$, which gives the stated
result as $\tr_{i,j,t}^r = (I -T_{Q_{i,j,t}})^r$. By (2) of  Lemma
\ref{lem:Vnf},
$$
   \|\tr_{i,j,t}^r ( f - V_n f)\|_p =  \|\tr_{i,j,t}^r f - V_n  \tr_{i,j,t}^r  f\|_p
      \le (1+ c) \|\tr_{i,j,t}^r f\|_p,
$$
from which the stated inequality follows.
\end{proof}

\subsection{Direct and inverse theorems for best approximation}

We start with direct and inverse theorem characterized by our new modulus of
smoothness.

\begin{thm} \label{thm:best}
For $f \in L^p$ if $1 \le p < \infty$ and $f \in C(\sph)$ if $p =\infty$, we have
\begin{equation} \label{Jackson}
       E_n(f)_p \le c \, \o_r(f,n^{-1})_p, \qquad 1 \le p \le \infty.
\end{equation}
On the other hand,
\begin{equation} \label{inverse}
       \o_r(f,n^{-1})_p \le c\, n^{-r} \sum_{k=1}^n k^{r-1} E_{k-1}(f)_p, \qquad 1 \le p \le \infty,
\end{equation}
where $\o_r(f,t)_p$ and $E_n(f)_p$ are defined in
\eqref{eq:modulus} and \eqref{eq:bestEn}, respectively.
\end{thm}

\begin{proof}
When $r=1$ and $1\leq p< \infty$, we use Lemma \ref{lem:Vnf},
H\"older's inequality and the fact that $\int_{\SS^{d-1}} |
K_n(\la x, y \ra) |d\s(y) \le c$ for all $x \in \sph$ to obtain
\begin{align*}
 E_n(f)_p & \le \|f-V_{\lfloor \frac n 2 \rfloor} f\|_p \\
   &  \ll \left(\int_{\sph}\int_{\sph} |f(x)-f(y)|^p
    |K_{\lfloor \frac n 2 \rfloor} (\la x, y \ra)|\, d\s(x)\, d\s(y)\right)^{\f 1p},
\end{align*}
from which \eqref{Jackson} for $r =1$ follows from Lemma \ref{KeyLemma}.
For $r =1$ and $p = \infty$, we use Lemma \ref{lem:omega_infty} and $V_nf$ to
conclude
\begin{align*}
 E_n(f)_\infty \le  \|f-V_{\lfloor \f n 2 \rfloor} f\|_\infty &
  \le \int_{\sph} |f(x)-f(y)| |K_{\lfloor \f n 2 \rfloor}(\la x,  y\ra)|\,d\s(y)\\
&  \ll \int_{\sph}\o (f, d(x,y))_\infty |K_{\lfloor \f n 2 \rfloor}(\la x, y \ra)|\, d\s(y)\\
&  \ll \o(f, n^{-1})_\infty  \int_{\sph} (1+nd(x,y)) |K_{\lfloor \f n 2 \rfloor}(\la x, y \ra)|\,d\s(y) \\
& \ll \o (f, n^{-1})_\infty,
\end{align*}
where the last inequality follows from \eqref{local_kernel} and the fact that
$\la x, y\ra = \cos d(x,y)$, just like the estimate of $A_j$ in the previous
proof.

Fo $r > 1$, we follow the induction procedure
 on $r$ (\cite[p.106-107]{Dai2}, \cite[p.191-192]{DD1}) using
$V_nf$ in Lemma \ref{lem:Vnf} and Lemma \ref{lem:triV=Vtri}.
Assume that we have proven \eqref{Jackson} for some positive
integer $r \ge 1$. Let $g = f-V_{\lfloor \frac n 2 \rfloor}f$. It
suffices to show that $\|g\|_p \le c\, \o_{r+1}(f,n^{-1})_p$. The
definition of $V_n f$ implies that $V_{\lfloor \frac n 4
\rfloor}g=0$, so that
\begin{align*}
    \|g\|_{p} = \|g-V_{\lfloor \frac n 4 \rfloor}g\|_{p} \le
    c E_{\lfloor \frac n 4 \rfloor}(g)_{p}\le c_1 \o_r(g,n^{-1})_{p}.
\end{align*}
On the other hand, using (2) and (3) of the Proposition \ref{modulus}, we
obtain, for any $m \in \NN$,
\begin{align*} 
\o_{r} (g, t)_p & \le c_r t^r \int_{t}^{2^m t} \f {\o_{r+1} (g,u)_p}{u^{r+1}}\, du
       + c_r 2^{r+1}t^r \|g\|_p\int_{2^m t}^1 u^{-r-1}\, du\\
  & \le c_1(m,r) \o_{r+1}(g, t)_p + c_2(r) 2^{-mr}\|g\|_p, \notag
\end{align*}
where $c_2(r)$ is independent of $m$. Choosing $m$ so that $4^{-1}
 \le c_1 c_2(r) 2^{-mr} < 2^{-1}$, we deduce from these two equations that
$$
    \|g\|_p \le c\, \o_{r+1} (g,n^{-1})_p \le c\, \o_{r+1} (f,n^{-1})_p,
$$
where the last step follows from Lemma \ref{lem:triV=Vtri}.
This completes the proof of \eqref{Jackson}.

The proof of \eqref{inverse} follows the standard approach for deriving it
from the Bernstein inequality (see, for example, \cite[p. 208]{DL}). Upon using (ii) of
Lemma \ref{lem:Delta_ij}, it reduces to the Bernetein type inequality
$$
   \|D_{i,j}^r P\|_p \le c n^r \|P\|_p, \qquad P \in \Pi^d_n, \quad 1 \le i< j \le d,
$$
which, however, immediately follows from (ii) and (iii) of
Lemma \ref{lem:Delta_ij}.
\end{proof}

As a corollary of Theorem \ref{Jackson}, we have the following:

\begin{cor} \label{cor:BestApp}
For $ 0 < \a < r$ and $ f \in L^p(\sph)$ if $1 \le p < \infty$ and $f \in C(\sph)$
if $p = \infty$,
$$
   E_n(f)_p \sim n^{-\a} \qquad \hbox{and} \qquad \o_r(f,t)_p \sim t^{-\a}
$$
are equivalent.
\end{cor}

\subsection{Equivalence of modulus of smoothness and $K$-functional}

\begin{thm} \label{thm:equivalent}
Let $r \in \NN$ and let $f \in L^p$ if $1 \le p < \infty$ and $f \in C(\sph)$
 if $p=\infty$.  For $0 < t < 1$,
$$
   \o_r(f,t)_p \sim K_r(f,t)_p, \qquad 1 \le p \le \infty.
$$
\end{thm}

\begin{proof}
By (ii) of Lemma \ref{lem:Delta_ij} and the triangle inequality,
$$
   \| \tr_{i,j,\t}^r f\|_p \le  \| \tr_{i,j,\t}^r (f - g)\|_p +  \| \tr_{i,j,\t}^r g\|_p
    \ll \| f -g \|_p + \t^r  \| D_{i,j}^r g\|_p
$$
from which $\o_r(f,t)_p \ll K_r(f,t)_p$ follows. On the other hand, for $t> 0$ set
$n = \lfloor \frac 1 t \rfloor$, then by
Lemma \ref{lem:Vnf}, \eqref{Jackson} and (iii) of Lemma \ref{Delta_ij}
\begin{align*}
K_r(f,t)_p & \le \|f - V_n f\| + t^r \max_{1 \le i<j\le d} \|D_{i,j}^r V_n f\|_p \\
  & \ll  \o_r(f, n^{-1})_p + t^r n^r \max_{1\le i<j\le d} \|\tr_{i,j,n^{-1}}^r V_n f\|_p \\
  & \ll \o_r(f,n^{-1})_p \ll  \o_r(f,t)_p
\end{align*}
where the last step follows from (i) of Lemma \ref{Delta_ij}.
\end{proof}

The proof of the above theorem, together with Lemma
\ref{lem:triV=Vtri} and (iii) of  Lemma \ref{lem:Delta_ij},
yields a realization of the $K$-functional.

\begin{cor} \label{cor:realization}
Under the assumption of Theorem \ref{thm:equivalent},
$$
    K_r(f, n^{-1})_p \sim \|f - V_n f\|_p +
             n^{-r} \max_{1\le i<j \le d} \|D_{i,j}^r V_nf\|_p.
$$
\end{cor}

\subsection{Comparison with other moduli of smoothness}
We want to compare our modulus of smoothness $\o_r(f,t)_p$ with
$\o_r^*(f,t)_p$ defined in \eqref{omega*} and $\wt \o_r(f,t)_p$ defined in
\eqref{omegaD}. By the equivalence of modulus of smoothness and
$K$-functional, we can work with $K_r^*(f,t)_p$ given in \eqref{K-func*} and
$K_r(f,t)_p$ in \eqref{eq:K-func-sphere}. Thus, we need to deal with the
equivalence
$$
    \max_{1 \le i<j\le d} \|D_{i,j}^r g\|_p \sim \|\Delta_0^r g\|_p,
$$
where $\Delta_0$ is the Laplace-Beltrami operator.

Let $\CH_n^d$ denote the space of spherical harmonics of degree $n$ in
$d$-variables. The operator $\Delta_0$ has $\CH_n^d$ as its space
of eigenfunctions. For any $\a \in \RR$, we have
\begin{equation} \label{eigen}
   (-\Delta_0)^\a Y =  (n(n+d-1))^\a Y, \qquad Y \in \CH_n^d.
\end{equation}

\begin{lem} \label{commute}
For $1 \le i\neq  j \le d$ and $\a \in \RR$
$$
     D_{i,j} ( - \Delta_0)^\a = (- \Delta_0)^\a D_{i,j}.
$$
\end{lem}

\begin{proof}
Because of the density of the polynomials, we only need to
establish this commutativity for spherical polynomials, which can
be further decomposed in terms of spherical harmonics. Thus, it
suffices to work with spherical harmonics. By \eqref{eigen}, we
only need to show that $D_{i,j} \CH_n^d \subset \CH_n^d$. Let
$\Delta$ be the usual Laplacian operator. A straightforward
computation shows that
$$
  \partial_i^2 D_{i,j} = 2 \partial_1 \partial_2 + D_{i,j} \partial_i^2,
  \qquad  \partial_k^2 D_{i,j} = D_{i,j} \partial_k^2, \quad k \ne i, \, k \ne j,
$$
from which it follows readily that
$$
     \Delta D_{i,j} = D_{i,j} \Delta, \qquad 1 \le i,j \le d.
$$
This  implies that $D_{i,j} Y \in \ker \Delta$ if  $Y \in \ker \Delta$, which
shows $D_{i,j} \CH_n^d \subset \CH_n^d$.
\end{proof}

\begin{lem} \label{lem:normDij}
Let $f \in  \bigcup_{n=0}^\infty \Pi_n^d$. For $r \in \NN$,
\begin{equation} \label{normDij}
    \max_{1\le i < j \le d} \|D_{i,j} f\|_p \sim  \|(-\Delta_0)^{1/2} f \|_p,
           \qquad 1 < p < \infty,
\end{equation}
and, for $r \in \NN$,
\begin{equation}\label{normDij-r}
       \max_{1\le i < j \le d} \|D_{i,j}^r f\|_p \le c  \|(-\Delta_0)^{r/2} f \|_p, \quad
        1 < p < \infty.
\end{equation}
\end{lem}

\begin{proof}
The key ingredient of the proof is the following result proved in \cite{DDH}\
\begin{equation} \label{DDH}
     \|(- \Delta_0)^{1/2} f \|_p \sim \| \grad f\|_p, \qquad 1 < p < \infty,
\end{equation}
where $\grad f$ is defined by
$$
   \grad f := \nabla F \vert_{\sph}, \quad\hbox{with}\quad
       F(x): = f\left(\f {x}{\|x\|} \right), \quad  x \in \RR^d \setminus \{0\},
$$
in which $\nabla F:= (\partial_1 F, \ldots, \partial_d F)$ and $\partial_j : =
\frac{\partial}{\partial x_j}$. The norm $\|\grad f\|_p$ is taken over the
Euclidean norm of $\grad f$; that is, $\|\grad f\|_p =
\|\la \grad f, \grad f\ra^{1/2}\|_p$.

We relate $\|\grad f\|_p$ to $\|D_{i,j} f\|_p$. A straightforward computation shows
that
\begin{equation} \label{partial-Dij}
   \frac{\partial}{\partial x_j} \left[ f\left(\frac{x}{\|x\|}\right)\right]_{\|x\|=1}
      = \partial_j f - x_j \sum_{i=1}^d x_i \partial_i f =
          - \sum_{\{i:\  \ 1\leq i\neq j\leq d\}} x_i D_{i,j} f,
\end{equation}
where we have used $1 = \sum_{i=1}^d x_i^2$ in both equations. The
first equation of \eqref{partial-Dij} implies that
$$
  | D_{i,j} f | = | \la x_i e_j - x_j e_i ,  \grad f \ra | \le \la \grad f, \grad f\ra^{1/2}
$$
on $\sph$, from which follows immediately that $\|D_{i,j} f\|_p  \le \|\grad f\|_p$.
On the other hand, the second equation of \eqref{partial-Dij} implies that
$$
    \la \grad f, \grad f\ra = \sum_{j=1}^d \left(\sum_{i\neq j} x_i D_{i,j} f\right)^2
       \le  \sum_{j=1}^d \sum_{i\neq j}  (D_{i,j} f)^2
        \le d^2 \max_{1 \le i< j \le d} (D_{i,j} f)^2,
$$
so that we have
$$
      \|D_{i,j} f\|_p \le  \|\grad f\|_p \le d
       \max_{1 \le i < j\le d} \|D_{i,j} f\|_p,
$$
which proves, upon using \eqref{DDH}, the equivalence in \eqref{normDij}.
Furthermore, by the commutativity of $D_{i,j}$ and $(-\Delta_0)^\a$ in
Lemma \ref{commute}, we have
$$
      \| D_{i,j}^r  f\|_p \le  c \| (-\Delta_0)^{1/2} D_{i,j}^{r-1} f\|_p
          = c  \|  D_{i,j}^{r-1} (-\Delta_0)^{1/2} f\|_p
$$
from which the inequality \eqref{normDij-r} follows from induction on $r$.
\end{proof}

\begin{rem} \label{rem:normDij}
The decomposition \eqref{Laplace-Betrami2} implies immediately that
\begin{equation} \label{DeltaD_r=2}
   \|(-\Delta_0) f \|_p \le  \f{d(d-1)}2 \max_{1\le i< j \le d} \|D_{i,j}^2 f\|_p, \qquad
        1 \le p \le \infty.
\end{equation}
Together with \eqref{normDij} and \eqref{normDij-r}, we see that
\begin{equation} \label{conjecture}
  \|(-\Delta_0)^{r/2} f \|_p  \sim  \max_{1 \le i < j \le d} \|D_{i,j}^r f\|_p, \qquad
        1 < p < \infty,
\end{equation}
holds for $r =1$ and $r=2$. However, we do not know if
\eqref{conjecture} is true for all $r$.
\end{rem}

We can now state and prove our main result in this subsection:

\begin{thm} \label{thm:K-funcEqu}
Let $f \in L^p(\sph)$, $1< p < \infty$. For $r \in \NN$ and $0 < t <1$,
\begin{equation} \label{eq:K=K1}
         K_r(f,t)_p \le c \,  K^*_r(f,t)_p, \qquad 1 < p < \infty.
\end{equation}
Furthermore,  for $r = 1$ or $2$,
\begin{equation} \label{eq:K=K2}
         K_r(f,t)_p \sim \, K^*_r(f,t)_p, \qquad 1 < p < \infty.
\end{equation}
\end{thm}

\begin{proof}
We only need to prove the inequality with $t = n^{-1}$. Furthermore, by
Corollary \ref{cor:realization}, we only need to work with polynomials
$V_n f$, for which the stated results are immediate consequences of
Lemma \ref{lem:normDij} and Remark \ref{rem:normDij}.
\end{proof}

For $1< p < \infty$, it is shown in \cite{DDH} that $\wt
\o_r(f,t)_p \sim K_r^*(f,t)_p$ so that $\wt \o_r(f,t)_p \sim
\o_r^*(f,t)_p$ for $1 < p < \infty$. As a corollary of Proposition
\ref{omega-omegaD}, Theorems \ref{thm:equivalent} and \ref{thm:K-funcEqu},
we can state the following equivalence:

\begin{cor} \label{cor:moduliEqu}
Let $f \in L^p(\sph)$ with $1< p < \infty$. For $r \in \NN$ and $0
< t <1$,
\begin{equation} \label{eq:omega=omega1}
         \o_r(f,t)_p \le \wt \o_r(f,t)_p \sim \o_r^*(f,t)_p, \qquad 1 < p < \infty.
\end{equation}
Furthermore,  for $r = 1$ or $2$,
\begin{equation} \label{eq:omega=omega2}
        \o_r(f,t)_p \sim \wt \o_r(f,t)_p \sim \o_r^*(f,t)_p, \qquad 1 < p < \infty.
\end{equation}
\end{cor}

According to the inequality \eqref{eq:omega=omega1}, our new
modulus of smoothness $\o_r(f,t)_p$ is at least as good as
$\o_r^*(f,t)_p$ for $1 < p < \infty$ and $r \ge 1$, and they are
equivalent when $r =1,2$.
 We do not know if the equivalence holds
for $r \ge  3$. In the case of $r =2$, the inequality
\eqref{DeltaD_r=2} shows that $ \o^\ast_2(f,t)_p \le c \,
\o_2(f,t)_p$ also holds for $p=1$ and $p = \infty$. A recent
example of \cite{Di3} shows that the equivalence
\eqref{eq:omega=omega2} fails at the endpoints $p=1, \infty$.


\section{Weighted Approximation on the Unit Sphere}
\setcounter{equation}{0}

In this section we consider approximation in the weighted $L^p$
space on the sphere. The main result establishes the analogue of
the Jackson estimate for the doubling weight using our new modulus
of smoothness. Such a result has been established in \cite[p. 94]{Dai2}
using the weighted version of $\wt \o_r(f,t)_p$, following the
lead of \cite{MT1,MT2} for weighted approximation on the interval.
We shall follow the approach in \cite{Dai2} closely.

\subsection{Definition of modulus of smoothness for doubling weight}

A non-negative integrable  function $w$ on $\sph$ is called  a
{\it doubling weight}  if there exists a constant $L>0$, called
the doubling constant,  such that for any $x\in\sph$ and $t>0$
$$
  \int_{c(x,2t)}w(y)\, d\s(y) \leq L \int_{c(x,t)} w(y)\, d\s(y).
$$
Many of the weight functions on $\sph$ that appear in analysis
satisfy the doubling condition, including all weights of the form
\begin{equation}\label{1-3}
 h_{\a, \mathbf{v}}(x) = \prod_{j=1}^m |x\cdot v_j|^{\a_j}, \qquad
  \a_j>0, \quad v_j\in \sph,
\end{equation}
as shown in \cite[(5.3)]{Dai}, which contains reflection invariant
weight functions introduced by Dunkl (see \cite[Chapt. 5]{DX}). Throughout
this section, we assume that $w$ is a  doubling weight with the
doubling constant $L$. We further define
\begin{equation}\label{eq:wn}
w_n(x) := n^{d-1} \int_{c(x,\f1n)} w(y)\, d\s(y),\qquad
n=1,2,\ldots, \   \
   x\in\sph
\end{equation}
and set $w_0(x) := w_1(x)$. Then $w_n$ is again a doubling weight
with doubling constant  comparable to $L$. Moreover, it satisfies
the following inequality
\begin{equation}\label{wn-wn}
  w_n(x)\leq L (1+nd(x,y))^s w_n(y),\quad   s={\log L}/{\log 2},
    \quad n=0,1,\ldots.
\end{equation}

We denote by $L^p(w)$ the weighted Lebesgue space endowed with the
norm
\begin{equation}\label{eq:norm-wn}
\|f\|_{p,w}:=\left(\int_{\sph} |f(y)|^pw(y) d\s (y)\right)^{1/p}
\end{equation}
with the usual change when $p=\infty$.
For $f \in L^p(w)$, our weighted $r$-th moduli of smoothness are
defined by
\begin{equation}\label{weight-modulus}
\o_r(f, t)_{p,w_n}: =\max_{1\leq i <j \leq d} \sup_{|\t|\leq t}\|
     \tr_{i,j, \t}^r f\|_{p,w_n},  \qquad  0< p \leq \infty,
\end{equation}
and the corresponding weighted $r$th order $K$-functional is
defined by
\begin{equation}\label{weight-Kfunc}
 K_r(f, t)_{p,w_n} :=\inf_{g \in C^r(\sph)}
   \left\{\|f-g\|_{p,w_n}+ t^r \max_{1\leq i<j\leq d}
\|D_{i,j}^r g\|_{p,w_n} \right\}.
\end{equation}
These definitions are analogues of those defined in \cite[p. 181]{MT2}
and \cite[p. 91]{Dai2}. They are used to study the weighted best
approximation defined by
\begin{equation}\label{weight-best}
E_k(f)_{p,w_n} :=\inf_{g\in\Pi_{k-1}^d} \|f-g\|_{p,w_n}, \qquad 1
\le p \le \infty.
\end{equation}

The direct and the inverse theorems for $E_k(f)_{p,w_n}$ were
established in \cite{Dai2} using the weighted modulus of
smoothness
$$
\wt{\o}_r (f, n^{-1})_{p,w_n}: = \sup_{Q\in O_{n^{-1}}} \|\tr_Q^r
f\|_{p,w_n}.
$$
Our development is parallel and follows along the same line. After
establishing the properties of the modulus of smoothness, most of
our proof will be similar to the unweighted case in Section 3, so
that we can be brief.

\subsection{Properties of modulus of smoothness}
It was shown in \cite[Corollary 3.4]{Dai} that if $f\in\Pi_n^d$
then $\|f\|_{p,w}\sim\|f\|_{p,w_n}$ for all $0< p<\infty$, with
the constant of equivalence depending only on $L$ and $d$. An
important tool for our study is the Marcinkiewicz-Zygmund
inequality. Let $\b > 0$. A subset $\Lambda$ of $\sph$ is said to
be maximal $\b$-separated if $\sph= \bigcup_{\o \in\Lambda}
c(\o,\b)$ and $\min \{d(\o,\o'): \o,\o'\in\Lambda, \o\neq \o' \}
\ge \b$. The following result is a simple consequence of
\cite[Corollary 3.3]{Dai}.

\begin{lem}\label{lem-3-2-ver2}
There exists a positive  number  $\ve$ depending only on $d$ and
$L$ such that  for any maximal $\frac{\d}{n}$-separated subset
$\Lambda$ of $\sph$ with $0<\d \leq \ve$, $f\in \Pi_n^d$ and
$0<p<\infty$, we have
\begin{align*}
\|f\|^p_{p,w} \sim   \sum_{\o \in\Lambda} \l_\o\min_{x\in c(\o, \f
\d n)}
    |f(x)|^p \sim  \sum_{\o \in\Lambda}\l_\o  \max_{x\in c(\o, \f \d n)}|f(x)|^p,
\end{align*}
where $\l_\o =\int_{c(\o, \d n^{-1})} w(x)\, d\s(x)$, and the
constants of equivalence depend only on $L$, $d$ and $p$.
\end{lem}

We start with the following analogue of Lemma \ref{lem:Delta_ij}:

\begin{lem}\label{lem-3-4-ver2}
Let  $r \in \NN$ and $f\in C^r(\sph)$.
\begin{enumerate}[\rm(i)]
\item If  $0<|t|\leq cn^{-1}$ and $1\leq p<\infty$
then
$$
 \|\tr_{i,j,t}^r f\|_{p,w_n} \ll |t|^r \|D_{i,j}^r f\|_{p,w_n},\qquad  1\leq i<j\leq d.
$$
\item Let  $f\in \Pi_n^d$, $1 \le p< \infty$ and let  $\ve$ be as in
Lemma \ref{lem-3-2-ver2}. If  $0<|t|\leq \f \ve{nr} $, then
$$
\|\tr_{i,j, t}^r f\|_{p,w}\sim |t|^r \|D_{i,j}^r f\|_{p,w}, \qquad
1\leq i<j\leq d,
$$
\end{enumerate}
where $w_n$ and $\|\cdot \|_{w_n}$ are defined in \eqref{eq:wn} and \eqref{eq:norm-wn},
respectively.
\end{lem}

\begin{proof}
(i) Let $F_{i,j}(t,x): = f(Q_{i,j,t}x)$.  Note that for any
$t_1,t_2
> 0$ and $x \in \sph$,
$$
   F_{i,j}(t_1+t_2,x) = f(Q_{i,j,t_1+t_2}x) =
   F_{i,j}(t_1,Q_{i,j,t_2}x),
$$
it follows  by the definition of the $\tr_{i,j,t}$ that
\begin{align}\label{3-2-eq}
\tr_{i,j,t}^r f(x) & = \int_0^{t} \cdots\int_0^{t}
\frac{\partial^r}{\partial t^r}
 F_{i,j}( t_1+\cdots+t_r, x)\, dt_1\cdots dt_r \\
  & = \int_0^{t} \cdots\int_0^{t}   \frac{\partial^r}{\partial t^r}
 F_{i,j}\left(0, Q_{i,j, t_1+\cdots+t_r}x\right)\, dt_1\cdots dt_r,
\end{align}
which implies, by Minkowski's inequality, that
\begin{align*}
\|\tr_{i,j,t}^r f\|_{p,w_n}\leq \int_{-|t|}^{|t|}
\cdots\int_{-|t|}^{|t|}\left \|\frac{\partial^r}{\partial t^r}
 F_{i,j}\left(0, Q_{i,j, t_1+\cdots+t_r}x\right)\right\|_{p,w_n}\, dt_1\cdots dt_r.
\end{align*}
Since, by \eqref{wn-wn}, $w_n(y) \sim w_n(x)$ whenever $d(x,y) \le
c n^{-1}$ and $d(Q_{i,j,t}x, x) \le t$, it follows from the
rotation  invariance of $d\s$ and \eqref{partial_ij} that for
$0<|t|\leq n^{-1}$,
\begin{align*}
\|\tr_{i,j,t}^r f\|_{p,w_n}  \ll \int_{-|t|}^{|t|}
\cdots\int_{-|t|}^{|t|}
   \left \| D_{i,j}^r f \right\|_{p,w_n}\, dt_1\cdots dt_r \ll |t|^r \|D_{i,j}^r f\|_{p,w_n}.
\end{align*}

(ii) Let $\Lambda$ be a maximal $\f{\ve}{ n}$-separated subset of
$\sph$ with $\ve$ being the same constant as in Lemma
\ref{lem-3-2-ver2}. Using (\ref{3-2-eq}),  for any $\o \in\Lambda$
and $0<|t|\leq \f \ve{ nr}$, we have
\begin{align*} 
|\tr_{i,j,t}^r (f)(\o)| &=\left |\int_0^{t} \cdots\int_0^{t}
\frac{\partial^r}{\partial t^r}
 F_{i,j}\left(0, Q_{i,j, t_1+\cdots+t_r} \o\right)\, dt_1\cdots dt_r\right | \notag\\
& \leq |t|^r \max_{0\leq u\leq r|t|} \left |
\frac{\partial^r}{\partial t^r}
 F_{i,j}(0,Q_{i,j,u}\o)\right |\leq |t|^r \max_{y\in c(\o, \ve/ n)}
\left | D_{i,j}^r f(y)\right|.
\end{align*}
Thus using Lemma \ref{lem-3-2-ver2} and setting $\l_\o
=\int_{c(\o, \ve n^{-1})} w(x)\, d\s(x)$,  we obtain
\begin{align*}
\|\tr_{i,j, t}^r f\|_{p,w}^p & \sim \sum_{\o\in\Lambda} \l_\o
|\tr_{i,j,t}^r f(\o)|^p \leq c |t|^{rp}  \sum_{\o\in\Lambda} \l_\o
\max_{y\in c(\o, {\ve}/{n})}
   \left | D_{i,j}^r f (y)\right|^p.
\end{align*}
However,  as shown in the proof of Lemma \ref{commute}, $D_{i,j}^r
f \in \Pi_n^d$, so that the right hand side of the above
expression is, by Lemma \ref{lem-3-2-ver2}, equivalent to
$|t|^{rp} \| D_{i,j}^r f \|_{p,w}^p$. Thus, we have established
the desired upper estimate  $\|\tr_{i,j, t}^r f\|_{p,w}\ll |t|^r
\|D_{i,j}^r f\|_{p,w}$.

The lower estimate can be carried out along the same line. In
fact, using (\ref{3-2-eq}), for  $\o\in\Lambda$ and $|t|\leq \f
\ve{nr}$, we have
$$
|\tr_{i,j,t}^r (f)(\o)|\ge |t|^r \min_{y\in c(\o, \ve/n)} \left|
 \frac{\partial^r}{\partial t^r} F_{i,j}(0, x)\right |
     =  |t|^r \min_{y\in c(\o, \ve/n)}  \left |D_{i,j}^r f(y)\right |.
$$
Since  $\tr_{i,j,t}^r (f)\in\Pi_n^d$ and $D_{i,j}^r f \in\Pi_n^d$,
it follows by Lemma \ref{lem-3-2-ver2} that
\begin{align*}
\|\tr_{i,j, t}^r f\|^p_{p,w} & \sim  \sum_{\o\in\Lambda}\l_\o
    |\tr_{i,j,t}^r f(\o)|^p \\
   &    \ge c |t|^{rp}  \sum_{\o \in\Lambda} \l_\o \min_{y\in c(\o,n^{-1}\ve)}
  \left |D_{i,j}^r f(y)\right|^p  \sim |t|^{rp} \left \|D_{i,j}^r f \right\|_{p,w}^p
\end{align*}
This gives the desired lower estimate and completes the proof.
\end{proof}

\begin{lem}\label{2-2-lem}
Let $r \in \NN$ and let $f\in L^p(\sph)$ if $1\leq p<\infty$ and
$f \in C(\sph)$ if $p = \infty$.
\begin{enumerate}[\rm(i)]
\item For $\l \ge 1$ and $t>0$,
\begin{equation*}
\o_r(f,\l t)_{p,w_n}\leq c(p, r,w) (1+n \l t)^{sr/p} \l^{r}\o_r
(f,t)_{p,w_n}.
\end{equation*}

\item Let $m, n$ be positive integers, For $ 0 < t \le \f 1{2^m n}$,
\begin{equation*}
\o_r(f,t)_{p, w_{_{n}}} \leq   c_1(m,r)
\o_{r+1}(f,t)_{p,w_n}+c_2(r)
     \d^{mr} \|f\|_{p, w_n},
\end{equation*}
where $\d= \f { (L2^s)^{1/p}}  {(L2^s)^{1/p}+1}\in (0,1)$,
$c_1(m,r)>0$ depends only on $m$, $r$ and $L$, and $c_2(r)>0$
depends only on $r$ and $L$.
\end{enumerate}
\end{lem}

The analogue of this lemma using the weighted version of the
modulus $\wt \o_r (f,t)_{p,w_n}$ was proved in Lemma 2.2 and Lemma
4.1 of \cite{Dai2}. The proof there carries over to our new
modulus of smoothness with obvious modification.

\subsection{Weighted approximation on the sphere}
Our main result is the Jackson estimate in the following theorem:

\begin{thm}\label{thm-2-1}
Let $ f\in L^p(\sph)$ when  $1\leq p<\infty$ or $f\in C(\sph)$
when $p=\infty$. Then
$$
        E_n(f)_{p,w_n} \le c \, \o_r(f,n^{-1})_{p, w_n},
$$
where $E_n(f)_{p,w_n}$ and $\o_r(f,n^{-1})_{p, w_n}$ and defined in
\eqref{weight-best} and \eqref{weight-modulus}, respectively.
\end{thm}

The inverse theorem in terms of $\o_r(f,t)_{p,w_n}$ follows from
the one given in terms of $\wt \o_r(f,t)_{p,w_n}$, the weighted
version of the modulus given in \eqref{omegaD}, in \cite[p. 94]{Dai2},
since $\o_r(f,t)_{p,w_n} \le \wt \o_r(f,t)_{p,w_n}$ as shown in
Proposition \ref{omega-omegaD}. We can also state the following
theorem:

\begin{thm}
Let $ f\in L^p(\sph)$ when $1\leq p<\infty$ or $f\in C(\sph)$ when
$p=\infty$. Then
\begin{align*}
 \o_r (f, n^{-1})_{p, w_n} \sim K_r(f, n^{-1})_{p,w_n}.
\end{align*}
Furthermore, a realization of the $K$-functional is given by
$$
  K_r(f, n^{-1})_{p,w_n} \sim \|f-V_n f\|_{p,w_n}
  + n^{-r} \max_{1\leq i<j \leq d} \|D_{i,j}^r (V_n f)\|_{p,w_n}.
$$
\end{thm}

These two theorems are analogues of results in Subsections 3.2 and
3.3. Their proofs are also similar, using the properties of the
modulus of smoothness given above and the following two lemmas.
The first one, proved in Lemma 2.5 of \cite[p. 97]{Dai2}, is as follows.

\begin{lem}
For $0\leq k \leq 4n$, $f \in L^p(\sph)$ if $1 \le p \le \infty$,
and $f \in C(\sph)$ if $p = \infty$,
\begin{equation*}
  \|V_nf\|_{p, w_k}\leq c \|f\|_{p, w_k},
  \quad \hbox{and} \quad
\|f-V_nf\|_{p,w_k}\leq c E_n(f)_{p,w_k},
\end{equation*}
 where $c>0$ depends only on $d$ and $L$.
\end{lem}

The second lemma is the analogue of Lemma \ref{KeyLemma}. Let
$$
G_{n}(\t) \equiv G_{n,\ell}(\t) := n^{d-1}(1+n|\t|)^{-\ell}\qquad
\text{with \quad $\ell> d+s+\f sp$}.
$$

\begin{lem}\label{lem-2-3}
Suppose $f\in L^p(\sph)$ for $1\leq p<\infty$. Then
\begin{align*}
 \int_{\sph} & \int_{\sph} |f(x)-f(y)|^p G_n( d(x, y))w_n(x)
  \, d\s(x) \, d\s(y) \leq c \, \o(f, n^{-1})_{p,w_n}^p.
\end{align*}
\end{lem}

\begin{proof}
The proof is similar to that  of Lemma \ref{KeyLemma}. We only
list the necessary modification for the weighted cases. First we
need to replace $d\s(x)$ and $dx'$ by $w_n(x)d\s(x)$ and
$\bar{w}_n(x')dx'$, where $\bar{w}_n(x')=w_n (x',
\sqrt{1-\|x'\|^2})$. Second, when making the change of variable
$x'+b_{j-1}(u) \mapsto x'$, we need to use the estimate
$\bar{w}_n(x'-b_{j-1}(u)) \ll (1+n\|u\|)^s \bar{w}_n (x')$, which
follows from \eqref{wn-wn}. Third, Lemma \ref{2-2-lem} (i) and
(\ref{wn-wn}) have to be used several times in the proof, and it
is often necessary to replace $G_{n,\ell}$ by $G_{n, \ell-s}$ when
(\ref{wn-wn}) is used.

 We also refer to the proof  of
Lemma 3.1 in \cite{Dai2} for details.

\end{proof}

\part{Approximation on the Unit Ball}

This part is organized as follows. In Section 6 we derive a pair of new
modulus of smoothness and $K$-functional on the ball from the results
on  the sphere. In Section 7, we study another pair of new modulus of
smoothness and $K$-functional, which are extensions of those defined
by Ditzian-Totik on $[-1,1]$ to the ball. Finally, in Section 8, we discuss
extensions of our result on the unit ball to $W_\mu$ with $\mu$ being
a nonnegative real number.

\section{Approximation on the Unit Ball, Part I}
\setcounter{equation}{0}

We consider approximation on the unit ball $\BB^d$ and we
often deal with the weighted $L^p$ spaces $L^p(\BB^d, W_\mu)$ for $1 \le p
< \infty$, where the weight function is defined by
\begin{equation}\label{weightW}
   W_\mu(x) : = (1 - \|x\|^2)^{\mu -1/2}, \qquad \mu \ge 0.
\end{equation}
For $1 \le p < \infty$ we denote by $\|f\|_{p,\mu}$ the norm for
$L^p(\BB^d, W_\mu)$,
\begin{equation}
  \|f\|_{p,\mu}: = \left( \int_{\BB^d} |f(x)|^p W_\mu(x) dx \right)^{1/p},
\end{equation}
and $\| f\|_{\infty,\mu} : = \|f\|_\infty$ for $f \in C(\BB^d)$. When we need to
emphasis that the norm is taken over $\BB^d$, we write $\|f\|_{p,\mu}
= \|f\|_{L^p(\BB^d,W_\mu)}$.

\subsection{Preliminaries}
There is a close relation between orthogonal structure on the sphere
and on the ball, so much so that a satisfactory theory for the best
approximation on the ball, including modulus of smoothness and its
equivalent $K$-functional, can be established accordingly
\cite[Sect. 4]{X05a} and \cite[Sect. 3]{X05b}.

For $f \in L^p(\BB^d, W_\mu)$, the modulus of smoothness in \cite{X05b} is
defined by
\begin{equation} \label{moduousB-1}
 \o_r^*(f,t)_{p,\mu} : = \sup_{|\t| \le t} \| (I - T_\t^\mu)^{r/2} f\|_{p,\mu},
\quad 1 \le p \le \infty,
\end{equation}
where $T_\t^\mu$ is the generalized translation operator of the
orthogonal expansion, which can be written explicitly as an
integral operator (\cite[Theorem 3.6]{X05b}). A $K$-functional
$K^\ast_r(f,t)_{p,\mu}$ that is equivalent to this modulus of
smoothness is defined by
\begin{equation} \label{K-funcB-1}
 K_r^*(f,t)_{p,\mu} : = \inf_{g} \left\{ \|f -g \|_{p,\mu} +  t^r
    \|(-\CD_\mu)^{r/2} g \|_{p,\mu}\right\}, \quad 1 \le p \le \infty,
\end{equation}
where $\CD_\mu$ is the second order differential operator
\begin{equation}\label{D-mu}
  \CD_\mu:=  \sum_{i=1}^d (1-x_i^2)  \partial^2_i - 2
   \sum_{1\le i < j \le d}  x_i x_j \partial_i \partial_j  -
     (d+2 \mu)\sum_{i=1}^d x_i \partial_i,
\end{equation}
which has orthogonal polynomials with respect to $W_\mu$ as
eigenfunctions, see \eqref{D-eigen}. Both $ \o^*_r(f;t)_{p,\mu}$
and $ K^*_r (f;t)_{p,\mu}$ satisfy all the usual properties of
moduli of smoothness and $K$-functionals, and they can be used to
prove the direct and inverse theorems for
\begin{equation}\label{best-Ball}
E_n(f)_{p, \mu}: = \inf_{ g  \in \Pi_{n-1}^d} \| f - g \|_{p,\mu}.
\end{equation}
The approach in \cite{X05a} is based on treating both $T_\t^\mu$ and
$\CD_\mu$ as multiplier operators of the orthogonal expansions, and the
results can be deduced from weighted counterparts on the unit sphere.

We shall define a new modulus of smoothness in the case of $\mu =
\frac{m-1}{2}$ and $m\in \NN$. The reason that we consider such
values of $\mu$ lies in a close relation between the orthogonal
structure on $\SS^{d+m-1}$ and the one on $\BB^d$, which was
explored in \cite[ ]{X01}.
%

Given a function $f$ on $\BB^d$, we will frequently need to regard it as
a projection onto $\BB^d$ of a function $F$, defined on $\SS^{d+m-1}$ by
\begin{equation} \label{eq:F}
  F(x,x') : = f(x), \qquad (x,x') \in \SS^{d+m-1},\quad x \in \BB^d, \quad x' \in \BB^{m}.
\end{equation}
Under such an extension of $f$,  the equations \eqref{IntS-B} and
\eqref{IntS-Bm=1} become, for example,
\begin{equation}  \label{IntS-B-wt}
\int_{\SS^{d+m-1}} F(y) d\s(y) = \s_m \int_{\BB^d} f(x)
(1-\|x\|^2)^{\frac{m-2}{2}} dx,
\end{equation}
where $\s_m$ denotes the surface area of $\SS^{m-1}$ for $m \ge 2$
and $\s_1= 2$.

Let $\CV_n^d(W_\mu)$ denote the space of orthogonal polynomials
of degree $n$ with respect to the weight function $W_\mu$ on $\BB^d$.
The elements of $\CV_n^d(W_\mu)$ satisfies (cf. \cite[p. 38]{DX})
\begin{equation} \label{D-eigen}
\CD_\mu P =  - n(n+ d+ 2\mu-1 ) P \quad\hbox{for all} \quad P \in \CV_n^d(W_\mu).
\end{equation}
We denote by $P_n^\mu(x,y)$ the reproducing kernel of $\CV_n^d(W_\mu)$
in $L^2(B^d, W_\mu)$. It is shown in \cite[Theorem 2.6]{X01} that
\begin{equation}\label{reprodB-S}
    P_n^\mu (x,y) = \int_{\SS^{m-1}} Z_{n,d+m} \left (  \la x, y \ra +
        \sqrt{1-\|y\|^2} \, \la x', \xi \rangle \right ) d\s (\xi)
\end{equation}
for any $x, y \in \BB^d$ and  $(x,x') \in \SS^{d+m-1}$,  where
$Z_{n,d}(t)$ is the zonal harmonic defined in \eqref{zonal} and
$\mu = \frac{m-1}{2}$. For
$\eta$ being a $C^\infty$-function on $[0,\infty)$ that satisfies
the  properties as defined in Section 3.1, we define an operator
\begin{equation}\label{Vnf_mu}
 V_n^\mu f(x) := a_\mu \int_{\BB^d} f(y) K_n^\mu(x,y) W_\mu(y) dy, \quad x \in \BB^d,
\end{equation}
where $a_\mu$ is the normalization constant of $W_\mu$ and
\begin{equation}\label{kernelB}
   K_n^\mu(x,y): = \sum_{k=0}^{2n} \eta \left(\frac{k}{n}\right) P_k^\mu(x,y).
\end{equation}
The operator $V_n^\mu$ is an analogue of the operator $V_n$  in
\eqref{Vnf}, and it shares the same properties satisfied by $V_n
f$. In particular, the kernel $K_n^\mu$ is highly localized \cite[Theorem 4.2]{PX}
and an analogue of Proposition \ref{Vnf} holds for
$V_n^\mu f$ (\cite[p. 16]{X05a}):

\begin{lem} \label{lem:Vnf-mu}
Let $f \in L^p(\BB^d, W_\mu)$ if $1 \le p < \infty$ and $f \in C(\BB^d)$
if $p = \infty$. Then
\begin{enumerate}[  \rm(1)]
 \item $V_n^\mu f \in \Pi^d_{2n}$ and $V_n^\mu f = f $ for $f \in \Pi^d_n$.
 \item For $n \in \NN$, $\|V_n^\mu f\|_p \le c \|f\|_p$
 \item For $n \in \NN$,
 $$
        \|f - V_n^\mu f\|_{p,\mu} \le c E_n(f)_{p,\mu}.
 $$
\end{enumerate}
\end{lem}

The following result shows a further connection between $V_n F$ and
$V_n^\mu f$.

\begin{lem} \label{lem:VnF}
Let $V_n$ denote the operator defined in \eqref{Vnf} on $\SS^{d+m-1}$.
For $x \in \BB^d$, $(x,x') \in \SS^{d+m-1}$ and $F$ in \eqref{eq:F},
$$
    (V_n F) (x,x') = (V_n^\mu f) (x), \quad\hbox{where}\quad \mu = \frac{m-1}{2}.
$$
\end{lem}

\begin{proof}
By the definition of $V_n$ and \eqref{IntS-B-wt},
\begin{align*}
 (V_n F) (x,x') & = \s_{d+m} \int_{\SS^{d+m-1}} F(y) K_n(\la (x,x'), y\ra) d\s(y) \\
   &  =  \s_{d+m} \int_{\BB^d} f(v)\int_{\SS^{m-1}} K_n\bigl( \la x, v \ra +
   \sqrt{1-\|v\|^2} \la x', \xi \ra\bigr) d\s(\xi) W_\mu(v) dv \\
    & = \s_{d+m} \int_{\BB^d} f(v) K_n^\mu(x, v) W_\mu(v) dv = (V_n^\mu f)(x),
\end{align*}
where the third step follows from \eqref{reprodB-S} and the definitions of $K_n$
and $K_n^\mu$.
\end{proof}

\subsection{Modulus of smoothness and best approximation}
For a given function $f \in L^p(B^d, W_{\frac{m-1}{2}})$, the extension
$F$ in \eqref{eq:F} is an element of $L^p(\SS^{d+m-1})$ according to
\eqref{IntS-B-wt}. This relation can be used to define a modulus of
smoothness on the unit ball.

We denote by $\wt f$ the extension of $f$ in \eqref{eq:F} in the case
of $m=1$; that is,
\begin{equation} \label{wt-f}
  \wt f(x,x_{d+1}) = f(x), \qquad (x,x_{d+1}) \in \RR^{d+1}, \quad x \in \BB^d.
\end{equation}
Recall that $\tr_{i,j,\t} = \tr_{Q_{i,j,\t}}$ and $Q_{i,j,\t}$ is the rotation in angle $\t$
in the $(x_i,x_j)$-plane.

\begin{defn} \label{defn:modulus-B}
Let $\mu = \frac{m-1}{2}$ and $m\in \NN$. Let $f \in L^p(\BB^d, W_\mu)$
if $1 \le p < \infty$ and $f \in C(\BB^d)$ if $p = \infty$. For $r \in \NN$ and $t>0$
\begin{align}  \label{omegaB-def}
  \o_r(f,t)_{p,\mu} :=  \sup_{|\t|\le t} &
  \left\{ \max_{1 \le i<j \le d}  \| \tr_{i,j,\t}^r  f\|_{L^p(\BB^{d},W_ {\mu})}, \right. \\
& \quad
 \left. \max_{1 \le i \le d}  \| \tr_{i,d+1,\t}^r \wt f\|_{L^p(\BB^{d+1},W_ {\mu-1/2})} \right\},
\notag
\end{align}
where, for $m=1$, $ \| \tr_{i,d+1,\t}^r \wt f\|_{L^p(\BB^{d+1},W_ {\mu-1/2})}$
is replaced by  $\| \tr_{i,d+1,\t}^r \wt f\|_{L^p(\SS^d)}$.
\end{defn}

Several remarks are in order. First, the second  term in the right
hand side of \eqref{omegaB-def} is necessary, as  for any radial function $f$ and
$1\leq i<j\leq d$, $\tr_{i,j,\t}^r  f=0$. The first term is also necessary, as
will be shown in our examples in Section 10 (see the discussion after
Example \ref{ex-ball-4}).  Also,  the second term in
\eqref{omegaB-def} can be made more explicit by, recalling \eqref{Delta_ij},
$$
 \tr_{i,d+1,\t}^r \wt f(x,x_{d+1}) = \overrightarrow{\tr}_\t^r
   f \left (x_1,\ldots,x_{i-1}, x_i\cos (\cdot) - x_{d+1} \sin (\cdot),
     x_{i+1},\ldots,x_d \right),
$$
with the forward difference in the right hand side being evaluated
at $0$.
Second, when $\mu = 1/2$, or $m =2$, we have the unweighted case,
whereas if $m =1$ then $W_{\mu -1/2}$ becomes singular and
the following limit holds:
$$
 \lim_{\mu \to 0+} \| \tr_{i,d+1,\t}^r \wt f\|_{L^p(\BB^{d+1},W_ {\mu-1/2})}
   = \| \tr_{i,d+1,\t}^r \wt f\|_{L^p(\SS^d)},
$$
which follows from the limit relation that, for a generic function $f$,
\begin{align}\label{lim-mu=0}
& \lim_{\mu \to 0+} c_\mu \int_{\BB^{d+1}} f(x) (1-\|x\|^2)^{\mu-1} dx \\
& \qquad
 = \lim_{\mu \to 0+} c_\mu \int_0^1 \int_{\SS^d} f(s x')  d\s(x') s^d (1-s^2)^{\mu-1} ds
   =  \int_{\SS^d}   f (x') d\s(x'),  \notag
\end{align}
where $c_\mu = \s_{d+1}(\int_{\BB^{d+1}} (1-\|x\|^2)^{\mu-1}dx)^{-1}$.
Third, in \eqref{omegaB-def}, we have used the
notation $W_\mu$ for both weight function on $\BB^d$ and
$\BB^{d+1}$ which implies that $x$ in the definition of $W_\mu$ is
assumed to be in the appropriate set accordingly. Finally, just as we remarked
after Definition \ref{def:modulus}, the moduli $\o_r(f,t)_{p,\mu}$ is not
rotationally invariant and it relies on the standard basis $e_1,\ldots, e_d$
but independent of the order of $e_1,\ldots, e_d$.

We can also define $\o_r(f,t)_p$ in an equivalent but more compact form:
\begin{align}  \label{omegaB-def2}
 \o_r(f,t)_{p,\mu}: = \sup_{|\t|\le t} \max_{1 \le i<j \le d+1}
      \| \tr_{i,j,\t}^r \wt f\|_{L^p(\BB^{d+1},W_ {\mu-1/2})}.
\end{align}
Indeed, if $1 \le i<j \le d$, then $\tr_{i,j,\t}^r \wt
f(x,x_{d+1}) = \tr_{i,j,\t}^r f(x)$ by the definition of
$Q_{i,j,\t} x$ in \eqref{EulerAngle1}; consequently,
$$
     \| \tr_{i,j,\t}^r \wt f\|_{L^p(\BB^{d+1},W_{\mu-1/2})} = c
              \| \tr_{i,j,\t}^r  f\|_{L^p(\BB^d,W_{\mu})},
$$
which follows from, for a generic function $f$ and $\l > -1$,
\begin{align} \label{B_d+1B_d}
 \int_{\BB^{d+1}} \wt f(y) (1-\|y\|^2)^\l dy & =
    \int_{\BB^d}  f(x) \int_{-\sqrt{1-\|x\|^2}}^{\sqrt{1-\|x\|^2}}
         (1-\|x\|^2 - u^2)^{\l} du dx \\
  & =   c \int_{\BB^d}  f(x) (1-\|x\|^2)^{\l +1/2} dx, \notag
\end{align}
where $c = \int_{-1}^1 (1-t^2)^{\l}\,dt$. Thus, \eqref{omegaB-def}
and \eqref{omegaB-def2} are  equivalent.

To emphasis the dependence on the dimension, we shall write the modulus
of smoothness on the sphere as $\o_r(f,t)_p = \o_r(f,t)_{L^p(\SS^{d-1})}$ in
the following.

\begin{lem} \label{lem:moduliBS}
Let $\mu =\frac{m-1}{2}$ and $m \in \NN$. Let $f \in L^p(\BB^d,W_{\mu})$
if $1 \le p < \infty$ and $f\in C(\BB^d)$ if $p  = \infty$, and let $F$ be defined
as in \eqref{eq:F}. Then
$$
   \o_r (f, t)_{L^p(\BB^d, W_\mu)} \sim  \o_r(F,t)_{L^p(\SS^{d+m-1})}.
$$
\end{lem}

\begin{proof}
If $1 \le i<j \le d$, then $\tr_{i,j,\t}^r F (x,x') =
\tr_{i,j,\t}^r f(x)$ by \eqref{EulerAngle1} and, hence, for $m \ge
1$, it follows by \eqref{IntS-B-wt} and \eqref{B_d+1B_d} that
\begin{align*}
   \int_{\SS^{d+m-1}}  |\tr_{i,j,\t}^r F(y)|^p d\s(y)   = \s_m
    \int_{\BB^d} |\tr_{i,j,\t}^r f(x)|^p (1-\|x\|^2)^{\frac{m-2}{2}} dx.
\end{align*}
If $1 \le i \le d$ and $ d+1 \le j \le d+m$, then it follows from
\eqref{Delta_ij} that $\tr_{i,j,\t}^r F(x,x') =  \tr_{i, d+1,\t}^r
\wt f(x, x_j)$, where $x \in \BB^d$, so that for $m \ge 2$,
$$
 \int_{\SS^{d+m-1}}  |\tr_{i,j,\t}^r F(y)|^p d\s(y)   = \s_{m-1}
    \int_{\BB^{d+1}} |\tr_{i,d+1,\t}^r \wt f(x)|^p (1-\|x\|^2)^{\frac{m-3}{2}} dx
$$
by \eqref{IntS-B} and \eqref{IntS-Bm=1}, whereas there is nothing
to prove for $m=1$ by the modification in the definition of
$\o_r(f,t)_{L^p(\BB^d, W_\mu)}$ in that case.
\end{proof}

\begin{thm} \label{thm:JacksonB}
Let $\mu = \frac{m-1}{2}$ and $m \in \NN$. For $f \in L^p(\BB^d,W_\mu)$
if $1 \le p < \infty$ and $f\in C(\BB^d)$ if $p  =\infty$,
\begin{equation}\label{JacksonB}
   E_n (f)_{p,\mu} \le c\, \o_r(f,n^{-1})_{p,\mu}, \qquad 1 \le  p \le \infty;
\end{equation}
on the other hand,
\begin{equation} \label{inverseB}
\o_r (f,n^{-1})_{p,\mu} \le c\, n^{-r} \sum_{k=1}^n k^{r-1}
E_{k}(f)_{p,\mu}
\end{equation}
where $ E_n (f)_{p,\mu}$ and $\o_r (f,t)_{p,\mu}$ are defined in \eqref{best-Ball} and
\eqref{omegaB-def}, respectively.
\end{thm}

\begin{proof}
Let $F$ be defined as in \eqref{eq:F}. By Lemma \ref{lem:VnF},
$(V_n^\mu f)(x)= (V_n F)(x,x')$, so that by \eqref{IntS-B-wt} and the Jackson
estimate for $F$ in \eqref{Jackson},
\begin{align*}
  \|V_n^\mu f - f\|_{p,\mu}^p & = \int_{\BB^d} |V_n^\mu f(x) - f(x)|^p W_\mu(x)dx \\
    & = c \int_{\SS^{d+m-1}} |V_n F(y)-F(y)|^p d \s(y) \\
    & \le c \,\o_r(F, n^{-1})_{L^p(\SS^{d+m-1})}
       \le  c\, \o_r(f, n^{-1})_{p,\mu},
\end{align*}
which proves \eqref{JacksonB}. The inverse theorem follows likewise from
\begin{align*}
  E_n(F)_{L^p(\SS^{d+m-1})} & \le  c \,
     \|V_{\lfloor \frac{n}{2} \rfloor} F - F\|_{L^p(\SS^{d+m-1})} \\
      &  = c \|V_{\lfloor \frac{n}{2} \rfloor}^\mu f - f\|_{p, \mu}
         \le  c E_{\lfloor \frac{n}{2} \rfloor}(f)_{p,\mu}
\end{align*}
and the inverse theorem for $F$ in \eqref{inverse}.
\end{proof}

\subsection{Equivalent $K$-functional and comparison}
Recall the derivatives $D_{i,j}$ defined in \eqref{partial_ij}. We use them to
define a $K$-functional.

\begin{defn} \label{defn:K-funcB}
Let $\mu = \frac{m-1}{2}$ and $m \in \NN$. Let $f \in L^p(\BB^d, W_\mu)$,
if $1 \le p < \infty$ and $f \in C(\BB^d)$ if $p = \infty$. For $r \in \NN$ and $t>0$,
\begin{align}  \label{K-functB-def}
 &   K_r(f,t)_{p,\mu}
:= \inf_{g\in C^r(\BB^d)} \Big\{  \|f-g\|_{L^p(\BB^d, W_\mu)}   \\
  & \qquad
     + t^r   \max_{1 \le i<j \le d}  \| D_{i,j}^r  g\|_{L^p(\BB^{d},W_ {\mu})} +
   t^r \max_{1 \le i \le d}  \| D_{i,d+1}^r  \wt g\|_{L^p(\BB^{d+1},W_ {\mu-1/2})} \Big \}.
\notag
\end{align}
where if $m =1$, then $\| D_{i,d+1}^r  \wt g\|_{L^p(\BB^{d+1},W_ {\mu-1/2})}$
is replaced by $\| D_{i,d+1}^r  \wt g\|_{L^p(\SS^d)}$.
\end{defn}

Although $\wt g(x,x_{d+1})  = g(x)$ is a constant in $x_{d+1}$
variable so that $\partial_{d+1} \wt g(x,x_{d+1}) =0$, we cannot
replace $D_{i,d+1}^r$ by $(x_i \partial_i)^r$, since $D_{i,d+1}^r
\ne (x_i \partial_i)^r$ if $r > 1$. Observe also,  that if
$f(x)=f_0(\|x\|)$ is a radial function and  $1\leq i<j\leq d$,
then $D_{i,j} (fg)(x)=f_0(\|x\|) D_{i,j}g(x);$ in particular, this
implies that  $D_{i,j}f=0$ for any radial function $f$.

In the case of $m \ge 2$, we can also define the $K$-functional in
an equivalent but more compact form
\begin{align*}
K_r(f,t)_{p, \mu} = \inf_{g\in C^r(\BB^d)} \Big\{  \|f-g\|_{L^p(\BB^d, W_\mu)}
 + t^r \max_{1 \le i<j \le d+1} \|D_{i,j}^r  \wt g \|_{L^p(\BB^{d+1},W_ {\mu-1/2})}
 \Big \}. \notag
\end{align*}
The equivalence of the two definitions follows from \eqref{B_d+1B_d}.
Just as in the case of modulus of smoothness,  these $K$-functionals are
related to the $K$-functional $K_r(f,t)_p =  K_r(f,t)_{L^p(\SS^{d-1})}$ defined
on the sphere.

\begin{lem} \label{lem:K-funcBS}
Let $\mu =\frac{m-1}{2}$ and $m \in \NN$. Let $f \in L^p(\BB^d,W_{\mu})$
if $1 \le p < \infty$ and $f\in C(\BB^d)$ if $p  = \infty$, and let $F$ be defined
as in \eqref{eq:F}. Then
$$
 K_r(f,t)_{p,\mu}\equiv  K_r (f, t)_{L^p(\BB^d, W_\mu)}
                  \sim  K_r(F,t)_{L^p(\SS^{d+m-1})}.
$$
\end{lem}
\begin{proof}
The estimate $K_r(F,t)_{L^p(\SS^{d+m-1})}\leq cK_r (f,
t)_{L^p(\BB^d, W_\mu)}$ follows directly from the definition, and
the fact that  for any $g\in C^r(B^d)$
$$\|D_{i,j}^r \wt g \|_{L^p(\BB^{d+1},W_ {\mu-1/2})} =c\|D_{i,j}^r
G \|_{L^p(\SS^{d+m-1})},\   \  1\leq i<j\leq d+1,$$ where
$G(x,x')=g(x)$ for $x\in B^d$ and $(x,x')\in\SS^{d+m-1}$.

To show the inverse inequality,  we observe that on account of
Lemma \ref{lem:VnF}, \eqref{IntS-B}, and \eqref{IntS-Bm=1},
 $\|f-V_n^\mu
f\|_{L^p(\BB^d, W_\mu)}=c \|F-V_n F\|_{L^p(\SS^{d+m-1})}$, and
$$\|D_{i,j}^r \wt {V_n^\mu f} \|_{L^p(\BB^{d+1},W_ {\mu-1/2})} =c\|D_{i,j}^r
V_n(F) \|_{L^p(\SS^{d+m-1})},\   \  1\leq i<j\leq d+1.$$The
inverse inequality $K_r (f, t)_{L^p(\BB^d, W_\mu)}\leq c
K_r(F,t)_{L^p(\SS^{d+m-1})}$ then follows by choosing $g=V_n^\mu
f$ with $n\sim \f1t$ in Definition \ref{defn:K-funcB}, and using
Corollary \ref{cor:realization}.
\end{proof}

By Lemmas \ref{lem:moduliBS} and \ref{lem:K-funcBS} and the
equivalence in Theorem \ref{thm:equivalent}, we further arrive at
the following:

\begin{thm}\label{thm:equivalenceB}
Let $r\in \NN$ and let $f \in L^p(\BB^d,W_\mu)$ if $1 \le p <
\infty$ and $f \in C(\BB^d)$ if $p  = \infty$. Then, for $0 < t <
1$,
$$
   \o_r (f , t)_{p,\mu} \sim K_r (f , t)_{p,\mu},   \qquad 1 \le p \le \infty.
$$
\end{thm}

Next we compare the moduli of smoothness $\o_r(f,t)_{p,\mu}$ with
$ \o^\ast_r(f,t)_{p,\mu}$ defined in \eqref{moduousB-1}. By
Theorem \ref{thm:equivalenceB} and \cite[Theorem 3.11]{X05b}, it
is enough to compare the $K$-functional $K_r(f,t)_{p,\mu}$ with $
K^*_r(f,t)_{p,\mu}$ defined in \eqref{K-funcB-1}. We start with an
observation on the differential operator $\CD_\mu$ defined in
\eqref{D-mu}. To emphasize the dependence on the dimension, we
shall use the notation $\Delta_{0,d} = \Delta_0$ for the
Laplace-Betrami operator on $\sph$.

\begin{lem}
Let $\mu = \frac{m-1}{2}$ and $m \in \NN$. Let $F$ be defined as in
\eqref{eq:F}. Then
$$
\Delta_{0,d+m} F (x,x') = \CD_\mu f (x), \qquad x \in \BB^d,
        \quad (x,x') \in \SS^{d+m-1}.
$$
\end{lem}

In fact, this follows immediately from comparing the expressions
\eqref{Laplace-Betrami3} and \eqref{D-mu}. Since both are
multiplier operators, their fractional powers are also defined and
equal. Thus, as a consequence of Lemma \ref{lem:K-funcBS} and
Corollary \ref{cor:realization}, we see that the comparison of the
$K$-functionals, Theorem \ref{thm:K-funcEqu}, and the comparison
of the moduli of smoothness,  Corollary \ref{cor:moduliEqu}, on
the sphere carry over to the comparison on the ball.

\begin{thm} \label{thm:moduliEqu-B}
Let $\mu = \frac{m-1}{2}$ and $m \in \NN$ and let $f \in L^p(\sph)$, $1< p < \infty$.
For $r \in \NN$ and $0 < t <1$,
\begin{equation} \label{eq:omega=omega1B}
         \o_r(f,t)_{p,\mu} \le c \, \o_r^*(f,t)_{p,\mu}, \qquad 1 < p < \infty.
\end{equation}
Furthermore,  for $r = 1$ or $2$,
\begin{equation} \label{eq:omega=omega2B}
        \o_r(f,t)_{p,\mu} \sim \o_r^*(f,t)_{p,\mu}, \qquad 1 < p < \infty.
\end{equation}
\end{thm}

An equivalent result can be stated for $K$-functionals.


\subsection{The moduli of smoothness on $[-1,1]$}
\label{subsection:d=1}

When $d =1$, the ball becomes the interval $B^1 = [-1,1]$. It turns out
that our modulus of smoothness appears to be new even in this case.
For $\mu = \frac{m-1}{2}$ and $m \in \NN$, the definition in
\eqref{omegaB-def} becomes, written out explicitly,
\begin{equation}\label{modulus_d=1}
\o_r(f,t)_{p,\mu} := \sup_{|\t| \le t}
  \left(c_\mu \int_{B^2} \left |\overrightarrow{\tr}_\t ^r
    f(x_1 \cos (\cdot)  + x_2 \sin (\cdot)) \right|^p  W_{\mu-\f12}(x) dx \right)^{1/p}
\end{equation}
for $1 \le p < \infty$ with the usual modification for $p =
\infty$, where $c_\mu^{-1}  = \int_{B^2} W_{\mu-\f12}(x)dx$. The
difference  $\overrightarrow{\tr}_\t ^r$ in this definition can be
evaluated at any fixed point $t_0\in [0,2\pi]$. More precisely,
 $\overrightarrow{\tr}_\t ^r f(x_1 \cos (\cdot)  + x_2 \sin (\cdot))=
    \overrightarrow{\tr}_\t ^r g_{x_1,x_2}(t_0)$  for a fixed
$t_0\in [0,2\pi]$, where $g_{x_1,x_2}(\t)=f(x_1 \cos \t  + x_2
\sin \t) $. Clearly, the definition is independent of the choice
of $t_0$, and makes sense for all real $\mu$ such that $\mu > 0$,
whereas for $\mu =0$ the integral is taken over $\SS^1$ upon using
the limit \eqref{lim-mu=0}.

This modulus of smoothness is computable, as shown by the
following example.

\begin{exam}
For $g_\a(x) = (1-x)^\a$, $\a > 0$ and $x\in [-1,1]$, $\mu > 0$ and $1 \le p
\le \infty$,
\begin{equation}\label{exam:d=1}
 \o_2(g_\a,t)_{p,\mu} \sim  \begin{cases}
       t^{2 \a + \frac{2\mu+1}{p}}, & - \frac{2\mu+1}{2 p} < \a < 1 - \frac{2\mu+1}{2p},\\
       t^{2} |\log t |^{1/p}, &  \a = 1 - \frac{2\mu+1}{2p}, \quad p \ne \infty,\\
       t^{2},  &  \a > 1 - \frac{2\mu+1}{2p}.
      \end{cases}
\end{equation}
\end{exam}

This is proved later in Lemma \ref{lem-4-3}. Notice that for $\mu = \frac{m-1}{2}$,
this modulus of smoothness is the restriction of the modulus of smoothness from
the sphere.

In this setting, several moduli of smoothness were defined and studied in the literature; we refer to the discussion in \cite[Chapter 13]{Di-To}.  In particular,
$\o_r^* (f, t)_{p,\mu}$ was studied by Butzer and his school and by Potapov.
The most successful one has been the Ditzian-Totik modulus of smoothness
\cite{Di-To}, which we now recall.

Let $\varphi(x) = \sqrt{1-x^2}$ and let $w_\mu(x):= (1-x^2)^{\mu-1/2}$
on $[-1,1]$. The Ditzian-Totik $K$-functional
 with respect to the weight $w_\mu$ is defined by
\begin{equation}\label{K-func-DT}
     \wh K_{r} (f,t)_{p,\mu} : = \inf_g  \left \{ \|f-g\|_{p,\mu} +
         t^r \|\varphi^{r} g^{(r)}\|_{p,\mu} \right \}.
\end{equation}
This K-functional  is equivalent to a modulus of smoothness
$\wh \o_r(f,t)_{p,\mu}$, called the Ditzian-Totik modulus of smoothness
and usually denoted by $\o_\varphi^r(f,t)_{p,\mu}$ (see \cite{Di-To}):
\begin{equation}\label{Eq-DT}
  \wh K_r(f,t)_{p,\mu}\sim \wh \o_r (f,t)_{p,\mu},\quad
     1\leq p\leq \infty, \quad 0<t<t_r.
\end{equation}
In the unweighted case (i.e. in the case of $\mu=\f12$), the
Ditzian-Totik modulus of smoothness  is defined  by, as in
\eqref{1-2omegaDT},
\begin{equation}\label{modu-DT}
  \wh \o_r (f,t)_{p} :=
  \sup_{0 < h \le t} \| \wh \tr_{h \varphi}^r f\|_{p},
   \qquad 1 \le p \le \infty,
\end{equation}
where $\wh \Delta_h^r$ is the $r$-th symmetric difference defined by
$$
 \wh \Delta^r_{\t \varphi} f(x) = \sum_{k=0}^r (-1)^k \binom{r}{k}
     f\left(x+ (\tfrac{r}{2} - k) \t \varphi(x)\right),
$$
in which we define $\wh \tr^r_{\t \varphi} f(x) =0$ whenever $x + r h \varphi(x)/2$
or $x - r h \varphi(x)/2$ is not in $(-1,1)$. In the weighted case
with $\mu > 1/2$, the Ditzian-Totik modulus of smoothness is more
complicated and defined by, for $1 \le p < \infty$,
\begin{align}\label{modu-DT-weight}
& \wh \o_r (f,t)_{p,\mu} :=
   \sup_{0 < h \le t} \|\wh \tr_{h \varphi}^r f\|_{L^p(I_{rh}; w_\mu)} \\
  & \quad + \sup_{0 < h \le 12r^2 t^2} \| \overleftarrow \tr_h^r f\|_{L^p(J_{1,rt}; w_\mu)}
  + \sup_{0 < h \le 12 r^2t^2} \| \overrightarrow \tr_h^r f\|_{L^p(J_{-1,rt}; w_\mu)},
    \notag
\end{align}
where the norms are taken over the intervals indicated with
$$
   I _{t} := [-1+ 2 t^2, 1- 2 t^2], \quad J_{1,t}:= [1- 12 t^2, 1], \quad
     J_{-1,t}: = [-1, -1+ 12 t^2],
$$
and $\wh \o_r (f,t)_{\infty,\mu}\equiv \wh \o_r (f,t)_{\infty} $
is defined as $\sup_{x\in[ -1,1]} |\wh\tr_{h \varphi}^r f(x)|$.

One important property of $\wh \o_r(f,t)_p$ is the following equivalence
established in \cite[(2.1.4), (2.2.5)]{Di-To}:
\begin{equation}\label{4-28-June4}
 \wh \o_r (f,t)_{p}^p\sim \f 1t \int_0^t \| \wh \tr_{h \varphi}^r f\|_{p}^p\, dh,\
    \   \  1\leq p<\infty,
 \end{equation}
with the usual change when $p=\infty$. In the weighted case, the right hand
side needs to be replaced by a sum of three integrals on the respective
intervals \cite[(6.19)]{Di-To}.

The success of the $\wh \o_r(f,t)_{p}$ lies in the fact that it is computable
and can be used to establish both the direct and  inverse theorems for
algebraic polynomial approximation on $[-1,1]$. The definition of
$\wh \o_r(f,t)_{p,\mu}$ for the weight $w_\mu$ is more complicated and will
be discussed in Subsection 7.8. For even more general weight,
see the book \cite{Di-To}.

The connection between our modulus of smoothness and that of
Ditizian-Totik is given in the following theorem.

\begin{thm}\label{muduli_d=1}
Let $\mu = \frac{m-1}{2}$, $m\in \NN$ and $r\in\NN$. Let $f \in
L^p([-1,1],w_\mu)$ if $1 \leq  p < \infty$, and $f\in C[-1,1]$ if
$p=\infty$. Assume further that $r$ is odd if $p=\infty$. Then
\begin{equation} \label{omega-omegaDT}
      \o_r(f,t)_{p,\mu} \le c\, \wh \o_r(f,t)_{p,\mu}
         +c\, t^r \|f\|_{p,\mu},   \quad  0<t\leq t_r,
\end{equation}
where the term $t^r\|f\|_{p,\mu}$ can be dropped when $r=1$.
\end{thm}

\begin{proof}
By Theorem \ref{thm:equivalenceB} and  the equivalence \eqref{Eq-DT},
it suffices to prove the inequality for the corresponding $K$-functionals:
$$
K_r(f,t)_{p,\mu}\leq c \wh K_r(f,t)_{p,\mu}+c \, t^r \|f\|_{p,\mu}, \quad
    1\leq p\leq\infty
$$
with the additional assumption $r$ is odd when $p=\infty$. This
inequality, together with  the equivalence $K_1(f,t)_{p,\mu}\sim
\wh K_{1}(f,t)_{p,\mu}$, is given in Theorem \ref{thm:whK=K} of
the next section.
\end{proof}

It is worth to point out that \eqref{4-28-June4} and
\eqref{omega-omegaDT} imply that
\begin{align} \label{tri-wh_tri}
 & \int_{\BB^2} \left |\overrightarrow{\tr}_\t ^r
    f(x_1 \cos (\cdot)  + x_2 \sin (\cdot)) \right|^p \frac{dx_1 dx_2}{\sqrt{1-x_1^2-x_2^2}} \\
 & \qquad\qquad\qquad
     \le c \frac{1}{t} \int_0^t \|\wh \tr_{h \varphi}^r f \|_p^p d h
    + c\, t^{rp} \|f\|_p^p, \notag
\end{align}
with the usual modification when $p = \infty$. This highly non-trivial inequality
will play a pivotal role in Subsection \ref{subsection-DT}.

We do not know if the reversed inequality of \eqref{omega-omegaDT}
holds when $r\ge 2$. We note, however, that the order of
$\o_r(g_\a,t)_{p,\mu}$ given in \eqref{exam:d=1} coincides, when
$\mu = 1/2$,  with the computation in \cite[p. 34]{Di-To} for $\wh
\o_r(g_\a,t)_{p}$ and $1 \le p \le \infty$  except when $\a =1$. For
$\a =1$, $g_\a(t)$ is a linear polynomial so that $\wh \o_r(g_\a,t)_{p}
 =0$, whereas $\o_r(g_\a,t)_{p, \frac12}$ is non-zero.


\section{Approximation on the Unit Ball, Part II}
\setcounter{equation}{0}

In this section, we  introduce  another pair of modulus of
smoothness and $K$-functional on the ball  that are in  analogy
with  those of Ditzian and Totik on $[-1,1]$, and  utilize them to
study best approximation on the unit ball. Both the direct and the
inverse theorems are established.

\subsection{A new $K$-functional and comparison} \label{subset:lem}
Let $\varphi(x): = \sqrt{1-\|x\|^2}$ for  $ x \in \BB^d$. Recall
that  the Ditzian-Totik K-functional $\wh K_{r}(f,t)_{p,\mu}$ on
$[-1,1]$ is  defined in \eqref{K-func-DT}. We now  define its
higher dimensional  analogue on the ball $\BB^d$.

\begin{defn} \label{defn:K-funcB2}
Let $f \in L^p(\BB^d, W_\mu)$ if $1 \le p < \infty$ and $f \in
C(\BB^d)$ if $p = \infty$. For $r \in \NN$ and $t>0$, define
\begin{align*}  
  \wh K_r(f,t)_{p,\mu}
   := \inf_{g\in C^r(\BB^d)} \Big\{ \|f-g\|_{p,\mu}
        + t^r   \max_{1 \le i<j \le d}  \| D_{i,j}^r  g\|_{p,\mu} +
      t^r \max_{1 \le i \le d}  \| \varphi^r \partial_i^r g\|_{p,\mu} \Big \}.
\end{align*}
\end{defn}

We establish a connection between  $\wh K_r(f,t)_{p,\mu}$ and
the $K$-functionals $K_r(f,t)_{p,\mu}$ defined in \eqref{defn:K-funcB}.
The result plays a crucial role in our development in this section. Recall,
in particular, that the proof of Theorem \ref{muduli_d=1} relies on
the theorem below.

\begin{thm} \label{thm:whK=K}
Let $\mu = \frac{m-1}{2}$ and $m \in \NN$. Let $f \in L^p(\BB^d,W_\mu)$
if $1 \le p < \infty$, and $f\in C(\BB^d)$ if $p=\infty$. We further assume
that $r$ is odd when $p=\infty$.  Then
\begin{align} \label{whK=K_r=1}
   \wh  K_1(f,t)_{p,\mu} \sim K_1(f,t)_{p,\mu},
\end{align}
and for $r > 1$, there is a $t_r > 0$ such that
\begin{align} \label{whK=K}
       K_r(f,t)_{p,\mu} \le c \, \wh K_r(f,t)_{p,\mu}
             + c \, t^r \|f\|_{p,\mu},  \qquad 0<t<t_r.
\end{align}
\end{thm}

The proof of Theorem \ref{thm:whK=K} relies on several lemmas. The
first one contains two Landau type inequalities. In the case of no
weight function and $r$ is even, this lemma appeared in \cite[p.
135]{Di-To} with $\|f\|_p$ in place of $\|\varphi^r f \|_{p}$ in
the right hand side of the inequalities. The proof of the general
case follows along the same line, but there are enough
modifications that we decide to include a proof.

\begin{lem} \label{weightHardy}
Let $\mu \ge 0$ and $r \in \NN$. Assume $f$ defined on $[-1,1]$ satisfies
$\varphi^r f^{(r)} \in L^p[-1,1]$.
\begin{enumerate}[\rm (i)]
\item If  $1\leq p<\infty$ and $1\leq i\leq \f r2$ or $p=\infty$ and
$1\leq i<\f r2$, then
\begin{equation}\label{weightHardy1}
     \| \varphi^{r-2i} f^{(r-i)}\|_{p,\mu}
        \le c_1 \|\varphi^{r} f^{(r)}\|_{p,\mu} +
        c_2\|\varphi^r f\|_{p,\mu}.
\end{equation}
\item If $r$ is even, set $\delta_r:= 0$ and assume $1\leq i\leq \f r2$,
 $1\leq p<\infty$;
If $r$ is odd, set $\d_r : =1$ and assume $1\leq i \leq \f{r+1}2$, $1\leq p\leq \infty$.
Then
\begin{equation}\label{weightHardy2}
   \|\varphi^{\delta_r} f^{(i)}\|_{p,\mu} \le c_1\|\varphi^{r} f^{(r)}\|_{p,\mu} +
                 c_2  \|\varphi^r f\|_{p,\mu}.
\end{equation}
\end{enumerate}
\end{lem}

\begin{proof}First, we show that for $1 \le p \le \infty$ and $1 \le i \le r-1$,
\begin{align} \label{weighted-Landau}
     \| f^{(i)}\|_{p,\mu}  \le  c \left( \|f^{(r)}\|_{p,\mu} +
         \|f\|_{p,\mu} \right).
\end{align}
In the case when  there is no weight function, this inequality is
well known. We only need to establish it for $1 \le p < \infty$.
We derive it from the following result in \cite[p. 109]{KZ},
\begin{align*} 
 \int_{0}^\infty  x^{\a} |g^{(i)}(x)|^p dx \le c \left(\int_{0}^\infty
        x^\a |g(x)|^p dx  \right )^{1-\frac{i}{r}}
   \left(\int_{0}^\infty  x^{\a} |g^{(r)}(x)|^p dx \right)^{\frac{i}{r}}
\end{align*}
for all $\a > -1$, which implies, by the elementary inequality $|a
b |\le \frac{|a|^p}{p}  + \frac{|b|^q}{q}$ for $\frac{1}{p} +
\frac{1}{q} =1$, that for $0\leq i\leq r$,
\begin{align} \label{KZ-Landau}
  \int_{0}^\infty  x^{\a} |g^{(i)}(x)|^p dx  \le
     c  \int_{0}^\infty    x^\a |g(x)|^p dx + c\int_{0}^\infty  x^{\a} |g^{(r)}(x)|^p dx.
\end{align}
For $f$ defined on $[-1,1]$, we write $f = f_1 + f_2 = f \psi + f
(1-\psi)$, where $\psi$ is a $C^\infty$ function on $\RR$ such
that $\psi(x) = 1$ for $x \le -1/2$ and $\psi (x) =0$ for $x \ge
1/2$. It then follows by \eqref{KZ-Landau}  that for $0\leq i\leq
r$
$$\|f_j^{(i)}\|_{p,\mu}\leq c \|f_j^{(r)}\|_{p,\mu} + c \|f
_j\|_{p,\mu}\leq c\|f_j^{(r)}\|_{p,\mu} + c \|f\|_{p,\mu},\   \
j=1,2.$$ Thus, the proof of \eqref{weighted-Landau} is  reduced to
showing
\begin{equation}\label{weighted-Landau-j}
\|f_j^{(r)}\|_{p,\mu}\leq c \|f\|_{p,\mu}+c \|f^{(r)}\|_{p,\mu},\
\  j=1,2.\end{equation} To see this, we observe that $\psi'$ is
supported in $[-\f12, \f12]$. Thus, by the Leibnitz rule, we
obtain
\begin{align*} \|f_j^{(r)}\|_{p,\mu}&\leq c \|f^{(r)}\|_{p,\mu}+ c
\max_{0\leq i\leq
r-1}\|f^{(i)}\|_{L^p[-\f12,\f12]}\\
& \leq c \|f^{(r)}\|_{p,\mu}+c\|f^{(r)}\|_{L^p[-\f12,\f12]}+
c\|f\|_{L^p[-\f12,\f12]}\\
&\leq c \|f^{(r)}\|_{p,\mu} + c \|f\|_{p,\mu},\end{align*} where
the second step uses the unweighted version of
\eqref{weighted-Landau}. This proves the desired inequality
\eqref{weighted-Landau-j}, and hence completes the proof of
\eqref{weighted-Landau}.

Now we return to the proof of
\eqref{weightHardy1} and \eqref{weightHardy2}. For $f =  f_1 +
f_2$ decomposed as above,  using \eqref{weighted-Landau-j} and the
fact that $\varphi^{rp} w_\mu= w_{rp/2 + \mu}$, we deduce
\begin{align} \label{weighted-Landau2}
   \|\varphi^r f_j^{(r)}\|_{p,\mu}  \le  c \left( \|\varphi^r f^{(r)}\|_{p,\mu} +
       \|\varphi^r f\|_{p,\mu} \right),\   \  j=1,2.
\end{align}
              Thus,  we can work with
$f_j$, $j=1,2$ instead of $f$. However, by symmetry, we can
assume, without loss of generality, that $f$ is supported in
$[-1,\f12]$. We claim that if $j\in\ZZ_+$ and $g\in C^j[-1,1]$ is
supported in $[-1, \f12]$ then
\begin{equation}\label{5-10-eq}
\|\varphi^ag\|_{p,\mu}\leq c \|\varphi ^{a+2j}g^{(j)}\|_{p,\mu}
\end{equation}
whenever $\mu+\f12+\f {ap}2>0$ and $1\leq p<\infty$ or $a>0$ and
$p=\infty$. Clearly, once (\ref{5-10-eq}) is proved, then
(\ref{weightHardy1}) follows by setting $g=f^{(r-i)}$, $a=r-2i$
and $j=i$, whereas \eqref{weightHardy2} follows by setting
$g=f^{(i)}$, $a=\d_r$ and $j=r-i$.

Clearly, for the proof of  the claim (\ref{5-10-eq}),  it suffices
to consider the case of  $j=1$. For $p=\infty$, we use the
inequality $|g(x)|\leq |\int_x^{\f 12} g'(t)\, dt|$ to obtain, for
$-1\leq x \leq \f12$,
\begin{align*} |\varphi^a(x) g(x)| &\leq c
(1+x)^{\f a2} \int_x^{\f12} |g'(t)|\, dt \\
& \leq c \|\varphi^{a+2}
g'\|_\infty(1+x)^{\f a2}\int_{x}^{\f12} (1+t)^{-\f a2-1}\, dt
\leq c \|\varphi^{a+2} g'\|_\infty,
\end{align*}
where we used the
assumption $a>0$ in the last step. This proves (\ref{5-10-eq}) for
$p=\infty$.

Next, we show  (\ref{5-10-eq}) for $1\leq p<\infty$. Again, we
only  need to consider $j=1$.
 Our main tool is the following  Hardy
inequality: for $1 \le p < \infty$ and $\b > 0$,
\begin{align} \label{Hardy}
 \left( \int_0^\infty \left(\int_x^\infty |f (y)| dy) \right)^p x^{\b-1} dx \right)^{1/p}
  \le \frac{p}{\b}  \left( \int_0^\infty  |y f(y)|^p y^{\b-1} dy \right)^{1/p}.
\end{align}

Setting  $G(y) = g( \frac{3}{2} y -1)$ with $y\in [0,1]$, we
 obtain
 \begin{align*}
 \|\varphi^a g\|_{p,\mu}^p &\leq c \int_{-1}^{\f12}
 |g(x)|^p(1+x)^{\mu-\f12+\f {pa}2}\,dx=c\int_{0}^1 |G(y)|^p
      y^{\mu-\f12+\f {pa}2}\, dy\\
 &\leq c \int_0^1\Bigl ( \int_y^1 |G'(x)|\, dx\Bigr)^p
 y^{\mu-\f12+\f {pa}2}\, dy\leq c \int_0^1 |yG'(y)|^p y^{\mu-\f12+\f {pa}2}\, dy \\
 &=c \int_{-1}^{\f12} |g'(x)|^p (x+1)^{\mu-\f12+\f {pa}2+p}\, dx
 \leq c \|\varphi^{a+2}g'\|_{p,\mu}^p
 \end{align*}
 proving the claim (\ref{5-10-eq}) for $j=1$. This completes the
 proof.
\end{proof}


\begin{lem} \label{lem:D12F}
Let $\wt f$ be defined as in \eqref{wt-f}. Then
$$
D_{1,d+1}^r \wt f(x,x_{d+1}) = \sum_{j=1}^r p_{j,r}(x_1,x_{d+1})
\partial_1^j f(x), \quad    x \in \BB^d,\  (x, x_{d+1})\in\BB^{d+1},
$$
where $p_{r,r}(x_1,x_{d+1})= x_{d+1}^r$ and
\begin{align}
p_{j,2r}(x_1,x_{d+1})  & = \sum_{\max\{0,j-r\} \le \nu \le j/2}
a_{\nu,j}^{(2r)}
       x_1^{j-2\nu} x_{d+1}^{2\nu},    \label{p_jrEven}\\
p_{j,2r -1}(x_1,x_{d+1}) & = \sum_{\max\{0,j-r\}  \le \nu \le
(j-1)/2}
     a_{\nu,j}^{(2r-1)} x_1^{j-1-2\nu} x_{d+1}^{2\nu+1}  \label{p_jrOdd}
\end{align}
for $1 \le j \le  2 r -1$ and $1 \le j \le  2 r -2$, respectively,
and $a_{\nu,j}^{(r)}$ are absolute constants.
\end{lem}

\begin{proof}
Recall that $\wt f(x,x_{d+1}) = f(x)$, so that $\partial_{d+1} \wt
f (x,x_{d+1}) =0$. Starting from
\begin{align*}
    D_{1,d+1}^{r+1} \wt f(x_1,x_{d+1}) = (x_{d+1} \partial_1 - x_1 \partial_{d+1})
        \sum_{j=1}^r p_{j,r}(x_1,x_{d+1}) \partial_1^j f(x),
\end{align*}
a simple computation shows that $p_{j,r}$ satisfies the recurrence
relation
\begin{equation} \label{p_jrRecur}
  p_{j,r+1} = x_{d+1} p_{j-1,r} + (x_{d+1} \partial_1 - x_1 \partial_{d+1}) p_{j,r},
   \quad 1 \le j \le r,
\end{equation}
where we define $p_{0,r} : =0$, and $p_{r+1,r+1}= x_{d+1}
p_{r,r}$. Since $p_{1,1} =x_{d+1}$, we see that $p_{r,r} =
x_{d+1}^r$ by induction. The general case also follows by
induction: assuming $p_{j,r}$ takes the stated from, we apply
\eqref{p_jrRecur} twice to get $p_{j,r+2}$ and verify that they
are of the form \eqref{p_jrEven} and \eqref{p_jrOdd}.
\end{proof}

We will also need the following integral formula, which is a simple
consequence of a change of variables.

\begin{lem}
For $1 \le m \le d-1$,
\begin{equation}\label{B-B}
  \int_{\BB^d} f(x)dx  = \int_{\BB^{d-m}} \left[\int_{\BB^m}
       f(\sqrt{1-\|v\|^2}u, v) du \right] (1-\|v\|^2)^{\frac{m}2} dv
\end{equation}
\end{lem}

We are now in a position to prove Theorem \ref{thm:whK=K}.

\medskip\noindent
{\it Proof of Theorem \ref{thm:whK=K}.} We give the proof for
the case $m \ge 2$ only. The proof for the case $m=1$ follows
along the same line. The only difference in this case is that we
need to replace the integral over $\BB^{d+1}$   by the one over
$\SS^d$ according to Definition \ref{defn:K-funcB}, and use
\eqref{IntS-Bm=1} instead of \eqref{B-B}.

By definition, we need to compare
$\| D_{i,d+1}^r \wt g\|_{L^p(\BB^{d+1},W_{\mu-1/2})}$ with $\| \varphi^r
\partial_i^r g\|_{p,\mu}$, where $\|\cdot\|_{p,\mu}\equiv
\|\cdot\|_{L^p(\BB^d, W_\mu)}$. More precisely, we need to show
\begin{equation}\label{5-15-Jun4}
\| D_{i,d+1} \wt g\|_{L^p(\BB^{d+1},W_{\mu-1/2})}\sim \| \varphi
\partial_i g\|_{p,\mu},\   \  1\leq i\leq d\end{equation}
and for $r\ge 2$
\begin{equation}\label{5-16-Jun4}
\| D_{i,d+1}^r \wt g\|_{L^p(\BB^{d+1},W_{\mu-1/2})}\leq c \|
\varphi^r
\partial_i^r g\|_{p,\mu}+c\|g\|_{p,\mu},\   \  1\leq i\leq d.\end{equation}

If $r =1$, then by \eqref{Delta_ij} $D_{1,d+1} \wt g (x,x_{d+1}) =
 x_{d+1} \partial_1 g(x)$. Hence, by \eqref{B-B},
\begin{align*}
  \|D_{i,d+1}\wt g\|_{L^p(\BB^{d+1},W_{\mu-1/2})}^p
  &   = \int_{\BB^{d+1}} | x_{d+1} \partial_1 g(x) |^p
  (1-\|x\|^2 - x_{d+1}^2)^{\mu-1} d(x,x_{d+1}) \\
  &   = c \int_{\BB^d} |\varphi(x) \partial_1 g(x)|^p (1-\|x\|^2)^{\mu-1/2} dx
        =  c \| \varphi \partial_1 g\|_{p, \mu}^p,
\end{align*}
where $c = \int_{-1}^1 |s|^p(1-s^2)^{\mu-1}ds$. The above argument
with slight modification  works equally well  for $p=\infty$. This
proves \eqref{5-15-Jun4}.

Next, we show \eqref{5-16-Jun4} for $r\ge 2$. By symmetry,  we
only need to consider the case $i =1$. We start with  the case of
even $r=2\ell$ with $\ell\in\NN$. In this case, $1\leq p<\infty$,
and by \eqref{p_jrEven}, we have
$$
   |D_{1,d+1}^{2\ell} \wt g(x,x_{d+1})| \le c\sum_{j=1}^{2\ell}
          \max_{ \max\{0, j-\ell\} \le \nu \le j/2}
             \left| x_1^{j-2\nu} x_{d+1}^{2\nu} \partial_1^j g(x)
             \right|.
$$
This implies
\begin{equation}\label{5-17-June4}\|D_{1,d+1}^{2\ell}
\wt g\|_{L^p(\BB^{d+1}, W_{\mu-1/2})}\leq c\sum_{j=1}^{\ell}
          \max_{ 0 \le \nu \le j/2} I_{j,\nu}+c\sum_{j=\ell+1}^{2\ell}
          \max_{ j-\ell \le \nu \le j/2} I_{j,\nu},
\end{equation}
with
\begin{align*}
I_{j,\nu} & :=      \int_{\BB^{d+1}} \left | x_1^{j-2\nu} x_{d+1}^{2\nu}
    \partial_1^j g(x) \right |^p  (1-\|x\|^2 - x_{d+1}^2)^{\mu-1} d(x, x_{d+1}) \\
 & =  c  \int_{\BB^d} \left | x_1^{j-2\nu} \varphi^{2\nu}(x) \partial_1^j g(x) \right |^p
                  (1-\|x\|^2)^{\mu-1/2} dx,
\end{align*}
where the last equation follows by  \eqref{B-B}. Let $x = (x_1, x') \in \BB^d$.
Using \eqref{B-B} again, and setting $g_{x'}(t) =g(t\varphi(x'),x')$, we
see that
\begin{align*}
  I_{j,\nu} &  =  c  \int_{\BB^{d-1}} \int_{-1}^1
        \left | t^{j-2\nu} \varphi^j(x') \varphi^{2\nu}(t) \partial_1^j
            g(\varphi(x') t, x') \right |^p
        (1-t^2)^{\mu-1/2} dt (1-\|x'\|^2)^{\mu} dx' \\
    & \le  c  \int_{\BB^{d-1}} \left [\int_{-1}^1
         | \varphi^{2\nu}(t) g_{x'}^{(j)}(t)  |^p
              (1-t^2)^{\mu-1/2} dt \right] (1-\|x'\|^2)^{\mu} dx',
\end{align*}
where the inequality is resulted from $|t^{j-2\nu}| \le 1$.

If $1 \le j \le \ell=\f r2$ and  $\nu\ge 0 $, then  $\varphi^{2\nu}(t) \le
1$, so that we can apply \eqref{weightHardy2} in Lemma
\ref{weightHardy} to conclude that
\begin{align*}
   I_{j,\nu} & \le  c  \int_{\BB^{d-1}}\left[ \int_{-1}^1
        \left |\varphi^{2\ell}(t)\frac{d^{2\ell}}{dt^{2\ell}}
              \left[ g(\varphi(x') t, x')\right] \right |^p
                       (1-t^2)^{\mu-1/2} dt\right] (1-\|x'\|^2)^{\mu} dx' \\
      & \qquad\qquad + c \int_{\BB^{d-1}} \int_{-1}^1 \left | g(\varphi(x') t, x')\right |^p
                       (1-t^2)^{\mu-1/2} dt (1-\|x'\|^2)^{\mu} dx' \\
      & = c  \int_{\BB^d} \left |\varphi^{2\ell}(x) \partial_1^{2\ell} g(x)\right |^p
                       (1-\|x\|^2)^{\mu-1/2} dx + c \|g\|_{p,\mu}^p.
\end{align*}

 If $\ell +1\leq  j \le 2\ell$, and  $\nu \ge j-\ell$, then  $\varphi^{2\nu}(t)
\le \varphi^{2j-2\ell}(t)$, so that we can apply
\eqref{weightHardy1} in Lemma \ref{weightHardy} with $i = 2\ell
-j$ and $r=2\ell$ to the integral over $t$, which leads exactly as
in the previous case to $I_{j,\nu} \le c  \|\varphi^{2\ell}
\partial_1^{2\ell} g\|_{p,\mu}^p+ c \|g\|_{p,\mu}^p$.

Putting these together, and using (\ref{5-17-June4}),  we have
established the desired result for the case of even $r=2\ell$. The
proof for  the case of odd $r$ follows along the same
line. This completes the proof. \qed

\begin{rem}
In the case of $r =2$ and $1 < p < \infty$, the reversed
inequality \eqref{whK=K} holds:
\begin{align} \label{whK=K-r=2}
   c_1 \wh K_2(f,t)_{p,\mu} \le  K_2(f,t)_{p,\mu} \le c_2 \, \wh K_2(f,t)_{p,\mu}
    + c \, t^2 \|f\|_{p,\mu},  \qquad 1 < p < \infty.
\end{align}
This will be proved in a more general setting in Theorem
\ref{thm:whK-K*}.
\end{rem}

We  do not know if the first inequality in \eqref{whK=K-r=2}
holds for $r  >2$, but they have to be close as both direct and
inverse theorems hold using either $K$-functional.

\subsection{Direct and inverse theorem by $K$-functional}
\label{5-4-subsection}

For the K-functional given in Definition \ref{defn:K-funcB2}, we
establish  both the direct and the inverse inequalities.

\begin{thm} \label{thm:JacksonB2}
Let $\mu = \frac{m-1}{2}$,  $m \in \NN$ and $r\in \NN$.  Let  $f
\in L^p(\BB^d,W_\mu)$ if $1 \le p < \infty$, and $f\in C(\BB^d)$
if $p=\infty$. Then
\begin{equation}\label{JacksonB2}
   E_n (f)_{p,\mu} \le c\, \wh
   K_r(f,n^{-1})_{p,\mu}+cn^{-r}\|f\|_{p,\mu}
\end{equation}
and
\begin{equation} \label{inverseB2}
\wh K_r (f,n^{-1})_{p,\mu} \leq  c\, n^{-r} \sum_{k=1}^n k^{r-1}
E_{k}(f)_{p,\mu}.
\end{equation}
Furthermore, the additional term $n^{-r}\|f\|_{p,\mu}$ on the
right hand side of \eqref{JacksonB2} can be dropped when $r=1$.
\end{thm}

\begin{proof}
When $1\leq p<\infty$ and $r\in\NN$ or $p=\infty$ and $r$ is
odd, the Jackson type estimate \eqref{JacksonB2} follows
immediately from \eqref{JacksonB} and Theorem \ref{thm:whK=K}.
Thus, it remains to prove \eqref{JacksonB2} for even   $r=2\ell$
and $p=\infty$. Since we already proved \eqref{JacksonB2}  for
$\wh K_{2\ell+1}(f,t)_\infty$, it suffices to show the inequality
\begin{equation}\label{5-20-eq}
   \wh K_{2\ell+1} (f, t)_\infty \leq c \wh K_{2\ell}(f,t)_\infty.
\end{equation}
For $d=1$, \eqref{5-20-eq} has already been proved in \cite[p. 38]{Di-To};
whereas in the case of $d\ge 2$ it is a consequence of the
following inequalities
$$
   \|D_{i,j}^{r+1} g\|_{\infty} \leq c \|D_{i,j}^{r} g\|_{\infty} \quad
    \text{and} \quad  \|\vi^{r+1} \partial_i^{r+1} g\|_{\infty}
   \leq c \|\vi^{r} \partial_i^{r} g\|_{\infty}
$$
which can be deduced directly from the corresponding results  for
functions of one variable; see, for example, \eqref{5-10-eq}.

The inverse estimate \eqref{inverseB2} follows as usual from the
Bernstein inequalities: For $1\leq p\leq \infty$ and $P\in\Pi_n^d$,
\begin{align}\label{Bernstein}
 \max_{1\leq i<j\leq d}
 \| D_{i,j}^r P\|_{p,\mu} \le c n^r \|P \|_{p,\mu} \quad \text{and} \quad
 \max_{1\leq i\leq d}
 \|\varphi^r  \partial_i^r P \|_{p,\mu} \le c n^r \| P \|_{p,\mu}.
\end{align}
The second inequality in \eqref{Bernstein} has already been
established in \cite[Theorem 8.2]{Dai}, so we just need to show
the first inequality.  Without loss of generality, we may assume
$(i,j)=(1,2)$. We then have, for $1\leq p<\infty$,
\begin{align}
\|D_{1,2}^r f\|_{p,\mu}^p &=\int_{\BB^{d-2}} \left [\int_{\BB^2}
|D_{1,2}^r f(\vi(u)x_1, \vi(u)x_2,u)|^p (1-x_1^2-x_2^2)^{\mu-\f12}\,
dx_1dx_2\right] \notag \\
& \qquad\qquad \times (1-\|u\|^2)^{\mu+\f12}\, du\notag \\
&=\int_{\BB^{d-2}}\int_{0}^1\left[\int_0^{2\pi} |D_{1,2}^r
f(\vi(u)\rho \cos\t, \vi(u)\rho\sin\t, u)|^p\, d\t\right]\notag\\
& \qquad \qquad \times
(1-\rho^2)^{\mu-\f12}\rho\, d\rho (1-\|u\|^2)^{\mu+\f12}\, du\notag\\
&=\int_{\BB^{d-2}}\int_{0}^1\left[\int_0^{2\pi} |f_{u,
\rho}^{(r)}(\t)|^p\, d\t\right] (1-\rho^2)^{\mu-\f12}\rho\, d\rho
(1-\|u\|^2)^{\mu+\f12}\, du\label{5-24-jun}\\
&\leq c n^{rp}\int_{\BB^{d-2}}\int_{0}^1\left[\int_0^{2\pi} |f_{u,
\rho}(\t)|^p\, d\t\right] (1-\rho^2)^{\mu-\f12}\rho\, d\rho
(1-\|u\|^2)^{\mu+\f12}\, du\notag\\
&=cn^{rp} \|f\|_{p,\mu}^p,\notag
\end{align}
where $f_{u,\rho}(\t) =f(\vi(u)\rho \cos\t, \vi(u)\rho \sin\t,u)$, and
the inequality step uses the usual Bernstein inequality for trigonometric
polynomials. Using \eqref{partial_ij2}, the same argument works for
$p=\infty$. This completes the proof of the inverse estimate.
\end{proof}

\begin{rem}
Since  the Bernstein inequality \eqref{Bernstein} is proved for
all $\mu>-\f12$, the inverse estimate \eqref{inverseB2} holds for
all $\mu>-\f12$ as well.
\end{rem}

\subsection{Analogue of Ditzian-Totik modulus of smoothness on $\BB^d$}
\label{subsection-DT} Recall the definition of the Ditzian-Totik
modulus of smoothness in \eqref{modu-DT}. We define an analogue on
the ball $\BB^d$. Since the definition for the weighted space has
an additional complication, we consider only the unweighted case,
that is, the case $W_{1/2}(x)dx = dx$, in this section. Let $e_i$
be the $i$-th coordinate vector of $\RR^d$ and let  $\wh\Delta_{h
e_i}^r$ be the $r$-th central difference in the direction of
$e_i$, more precisely,
$$
\wh\Delta_{h e_i}^r f (x) = \sum_{k=0}^r (-1)^k \binom{r}{k}
f\left(x+ (\tfrac{r}2 -k)h e_i\right).
$$
As in the case of $[-1,1]$, we assume that $\wh\Delta_{h e_i}^r$
is zero if either  of the points $x \pm r \frac{h}{2} e_i$ does
not belong to $B^d$. We  write $L^p(\BB^d)$, $\|f\|_p$ and $\wh
K_r(f,t)_p$ for $L^p(\BB^d, W_{1/2})$,
$\|f\|_{L^p(\BB^d,W_{1/2})}$ and $\wh K_r(f,t)_{p,1/2}$
respectively. The modulus of smoothness $\wh \o_r(f,t)_{p}$ in
\eqref{modu-DT} for the case $d=1$ suggests the following
definition:

\begin{defn} \label{defn:modulusB2}
Let $f \in L^p(\BB^d)$ if $1 \le p < \infty$ and $f \in C(\BB^d)$
if $p = \infty$. For $r \in \NN$ and $t>0$,
\begin{equation}\label{modulusB2}
  \wh \o_r(f,t)_{p} = \sup_{0 < |h| \le t} \left\{ \max_{1 \le i<j \le d}
      \| \tr_{i,j,h}^r  f\|_{p},
      \max_{1 \le i \le d}  \|\wh \tr_{h\varphi  e_i}^r f\|_{p} \right \}.
\end{equation}
\end{defn}

As in the case of Definitions \eqref{def:modulus} and \eqref{defn:modulus-B},
the new moduli are not rotationally invariant, they depend on the standard
basis $e_1,\ldots,e_d$ of $\RR^d$ but independent of the order of this basis.
In the case of $d=1$, there is no Euler angle and the definition becomes
exactly the one in \eqref{modu-DT}.

Much of the properties of the modulus of smoothness $\wh{\o}_r(f,t)_{p}$
follows from the corresponding properties of the moduli of smoothness on
the sphere and on $[-1,1]$. For example, we have the following lemma.

\begin{lem} \label{ball2-modulus}
Let $f \in L^p(\BB^d)$ for $1 \le p < \infty$ and $f \in C(\BB^d)$ for $p = \infty$.
\begin{enumerate}[\rm(1)]
\item For $0 < t < t_0$, $\wh \o_{r+1} (f,t)_{p} \le c\, \wh \o_r (f,t)_{p}$.
\item For $\l > 0$, $\wh \o_r(f,\l t)_{p} \le c (\l +1)^r \wh \o_r(f,t)_{p}$.
\item For $0<t<\f12$ and every $m>r$,
$$
   \wh \o_r(f,t)_{p} \le c_m \left( t^r \int_t^1 \f {\wh \o_{m}(f, u)_{p}}
         {u^{r+1}}\, du + t^r \|f\|_{p} \right).
$$
\item For $0 < t < t_0$, $\wh \o_{r} (f,t)_{p} \le c\, \|f\|_p$.
\end{enumerate}
\end{lem}

\begin{proof}
For $1 \le p < \infty$ and $\tr^r_{i,j,\t}f$, we use the integral formula
$$
   \|\tr_{i,j,\t}^r f\|_{p}^p  =  \int_0^1 s^{d-1} \int_{\SS^{d-1}}
        \left| \tr_{i,j,\t}^r f(s x') \right|^p d\s(x')d s
$$
and apply Proposition 2.5. For $\wh \tr^r_{\t \varphi e_i}f$, we use
\eqref{B-B} with $m=1$ and the fact that if $x = (\varphi(u)s,u)$
then $\varphi(x) = \varphi(s) \varphi(u)$ to conclude that
\begin{align*}
 \|\varphi^r \wh \tr^r_{\t \varphi e_i}f\|_{p}^p
& = \int_{BB^d} \left| \varphi^r(x) \wh \tr^r_{\t \varphi(x) e_i}f(x) \right|^p dx \\
& = \int_{\BB^{d-1}} \int_{-1}^1 \left|\varphi^r(u)\varphi^r(s)
   \wh \tr^r_{\t \varphi(s)\varphi(u) e_i}f(\varphi(u)s, u) \right|^p ds \varphi(u) du \\
& = \int_{\BB^{d-1}} (\varphi(u))^{rp+1} \left[ \int_{-1}^1
\left|\varphi^r(s)
        \wh \tr^r_{\t \varphi(s) e_i} f_u (s) \right|^p ds \right]
           du,
\end{align*}
where $f_u(s) = f(\varphi(u) s, u)$ and apply the result of one
variable in \cite[p. 38, p. 43]{Di-To} to the inner integral as
well as the equivalence \eqref{4-28-June4}.
\end{proof}

Next we establish the direct and the inverse theorems in $\wh \o_r(f,t)_p$, one
of the central result in this section.

\begin{thm} \label{thm:JacksonB2-Mod}
Let $f \in L^p(\BB^d)$  if $1 \le p < \infty$, and  $f \in C(\BB^d)$
if $p = \infty$. Then for $r \in\NN$
\begin{equation}\label{JacksonB2-2}
   E_n (f)_{p} \le c\, \wh\o_r(f,n^{-1})_p+n^{-r}\|f\|_{p}
\end{equation}
and
\begin{equation} \label{inverseB2-2}
\wh\o_r(f,n^{-1})_p \leq  c\, n^{-r} \sum_{k=1}^n k^{r-1} E_{k}(f)_{p}.
\end{equation}
Furthermore, the additional term $n^{-r}\|f\|_p$ on the
right hand side of \eqref{JacksonB2-2} can be dropped when $r=1$.
\end{thm}

\begin{proof}
We start with the proof of the Jackson type inequality \eqref{JacksonB2-2}.
By \eqref{JacksonB}, it suffices to show that for the modulus $\o_r(f,t)_p$
given in Definition \ref{defn:modulus-B},
\begin{equation}
  \o_r(f, n^{-1})_p \leq c\, \wh\o_r(f, n^{-1})_p+c n^{-r}\|f\|_p.
\end{equation}
However, using Definitions \ref{defn:modulus-B} and \ref{defn:modulusB2},
this amounts  to showing that for $1\leq i\leq d$
\begin{equation}\label{5-21-Jun4}
\sup_{|\t|\leq t}\|\tr_{i, d+1, \t} ^r \wt f \|_{L^p(\BB^{d+1},
W_0)}\leq c \, \wh\o_r(f,t)_p +c t^r \|f\|_p,
\end{equation}
where $\wt f(x, x_{d+1})=f(x)$ for $x\in\BB^d$ and $(x,x_{d+1})\in \BB^{d+1}$.
By symmetry, we only need to consider $i=1$. Set
$$
f_v(s) =f(\varphi(v)s, v),\   \   v\in \BB^{d-1},\   s\in [-1,1],
$$
where $\varphi(v)=\sqrt{1-\|v\|^2}$. We can then write, by \eqref{B-B},
\begin{align*}
&\|\tr_{1, d+1, \t}^r \wt f \|^p_{L^p(\BB^{d+1}, W_0)}\\
&=\int_{\BB^{d+1}} |\overrightarrow{\triangle}_\t^r f(x_1\cos
(\cdot)+x_{d+1}\sin (\cdot), x_2,\cdots, x_d)|^p\f
{dx}{\sqrt{1-\|x\|^2}}\\
&=\int_{\BB^{d-1}}\left[ \int_{\BB^2}
|\overrightarrow{\triangle}_\t^r f_v(x_1\cos (\cdot)+x_{d+1}\sin
(\cdot)|^p\f  {dx_1 dx_{d+1}}{\sqrt{1-x_1^2-x_{d+1}^2}}\right]
\varphi(v)\, dv.
\end{align*}
Applying \eqref{tri-wh_tri} to the inner integral,  the last expression
is bounded by, for $|\t|\leq t$,
\begin{align*}
& c \f 1t \int_0^t \int_{\BB^{d-1}} \left[\int_{-1}^1 |\wh
\triangle _{h\varphi(s)}^r f_v(s)|^p\, ds \right]\, \varphi(v)\,
dv\, dh+t^{rp}
\int_{\BB^{d-1}}\int_{-1}^1 |f_v(s)|^p\, ds \, \varphi(v)\, dv\\
& =c\f 1t \int_0^t \int_{\BB^{d-1}}\int_{-1}^1 |\wh
{\triangle}^r_{h\varphi(\varphi(v)s,v) e_1}f(\varphi(v)s,v)|^p\,
ds
\varphi(v)\, dv\, dh+c t^{rp} \|f\|^p_p\\
&=c\f 1t \int_0^t\int_{\BB^d} |\wh {\triangle}^r_{h\varphi(x)
e_1}f(x)|^p\, dx\, dh+c t^{rp} \|f\|^p_p\leq c \, \wh\o_r(f,t)^p_p +c
t^{rp} \|f\|^p_p.
\end{align*}
For $r =1$ the additional term $t^{rp} \|f\|^p_p$ can be dropped because
of Theorem \ref{muduli_d=1}. Obviously, the above argument with
slight modification works equally well for the case $p=\infty$.
This proves the Jackson  inequality (\ref{5-21-Jun4}).

Finally, the inverse estimate \eqref{inverseB2-2} follows by
\eqref{inverseB2} and  the inequality $\wh\o_r(f, t)_{p,\mu}
\allowbreak \leq c \wh K_r (f, t)_{p,\mu}$, which will be given in Theorem
\ref{thm:5-5} in the next subection.
\end{proof}

\subsection{Equivalence of $\wh \o_r(f,t)_p$ and $\wh K_r(f,t)_p$}
\label{5-7-subsection}

As a consequence of Theorem \ref{thm:JacksonB2-Mod}, we can deduce
the equivalence of the modulus of smoothness $ \wh\o^r(f,t)_p$ and the
K-functional $ \wh K_r(f,t)_p$:

\begin{thm}\label{thm:5-5}
Let   $f \in L^p(\BB^d)$  if $1 \le p < \infty$, and  $ f \in
C(\BB^d)$ if  $p = \infty$.
Then for $r\in\NN$ and $0<t<t_r$,
$$
c^{-1}\wh\o^r(f,t)_p\le  \wh K_r(f,t)_p\leq c\, \wh\o^r(f,t)_p+c\,t^r\|f\|_p.
$$
Furthermore, the term $t^r\|f\|_p$ on the right side can be dropped
when $r=1$.
\end{thm}

For the proof of Theorem \ref{thm:5-5}, we need the following lemma.
\begin{lem}\label{lem-5-9-Jun}
 For $1\leq p\leq \infty$ and $f\in \Pi_n^d$, we have
\begin{equation}\label{5-21-jun} n^{-r}\|D^r_{i,j} f\|_{p, \mu} \sim \sup_{|\t|\leq
n^{-1}} \|\triangle_{i,j,\t}^r f\|_{p,\mu},\   \ 1\leq i<j\leq d,
\end{equation}
and
\begin{equation}\label{5-22-jun}
n^{-rp}\|\varphi^r
\partial_i^r f\|^p_{p} \sim n\int_0^{n^{-1}}
\|\wh\triangle_{h\varphi e_i}^r f\|_{p}^p\, dh,\ \  1\leq i\leq d,
\end{equation}
with the usual change when $p=\infty$.
\end{lem}

\begin{proof}
The relation \eqref{5-21-jun} follows directly from
\eqref{5-24-jun} and the corresponding inequality for
trigonometric inequality (see, for instance, \cite{St}). The
relation \eqref{5-22-jun} can be proved similarly. In fact,
setting $i=1$ and $f_u(s)=f(\vi(u)s,u)$, we have
\begin{align}
 n^{-rp}\|\vi^r\partial_1^r f\|_{p}^p & =n^{-rp}\int_{\BB^{d-1}}\left[
\int_{-1}^1 |\vi^r(s)f_u^{(r)}(s)|^p\, ds \right] \vi(u)\, du\notag \\
&\sim n\int_{\BB^{d-1}}\left[ \int_0^{n^{-1}}\int_{-1}^1
|\wh\triangle_{h\vi(s)}^r f_u(s) |^p\, ds\, dh\right] \vi(u)\, du \\
&= n\int_0^{n^{-1}} \|\wh\triangle_{h\varphi e_1}^r f\|_{p}^p\,
dh,\notag
\end{align}
where we have used the equivalence of one variable in \cite[p. 191]{HL}
and \eqref{4-28-June4}.
\end{proof}

\noindent
{\it Proof of Theorem \ref{thm:5-5}.}
We start with the proof of the inequality
\begin{equation}\label{5-23-jun}
    \wh{\o}_r (f,t)_p \leq c \wh{K}_r(f, t)_p,\quad 0<t<t_r.
\end{equation}
Let $g_t \in  C^r(\BB^d)$ be chosen such that
$$
   \|f- g_t\|_{p} \le 2  \wh K_r(f,t)_{p},
      \quad t^r\max_{1\leq i<j\leq d}
      \|D_{i,j}^r g_t \|_{p}\le 2 \wh K_r(f,t)_{p},
$$
and
$$
 t^r\max_{1\leq i\leq d}\| \varphi^r \partial_i^r g_t\|_{p} \le 2 \wh K_r(f,t)_{p}.
$$
From the definition of $\wh \omega_r(f,t)_{p}$ and (4) of Lemma
\ref{ball2-modulus} it follows that
$$
   \wh \omega_r(f,t)_{p} \le  \wh \omega_r(f -g_t ,t)_{p} +
    \wh  \omega_r(g_t,t)_{p}\leq c \wh K_r(f,t)_{p}+\wh  \omega_r(g_t,t)_{p}.
$$
Consequently, for the proof of the inequality of \eqref{5-23-jun},
it suffices to show that for  $g \in C^r(\BB^d)$,
\begin{equation}\label{eq2}
 \|\tr^r_{i,j, \t} g\|_{p} \le  c\, \t^r\| D^r_{i,j}g\|_{p}
       \quad  \hbox{and}\quad
  \|\wh \tr^r_{\t \varphi e_i} g\|_{p} \le c\, \t^r
   \| \varphi^r \partial_i^r g\|_{p}.
\end{equation}
First we consider $\wh \Delta^r_{\t \varphi e_i} f$, for which we
will need the corresponding result for $[-1,1]$. It is known
\cite[(2.4.4)]{Di-To} that there exists $t_r\in (0,1)$ such
that for $0<h<t_r$,
\begin{align} \label{1d-est}
    \|\wh \tr^r_{h\varphi } g_t\|_{L^p[-1,1]} \le c \, h^r
   \|\varphi^r g_t^{(r)}\|_{L^p[-1,1]}.
\end{align}
For $p=\infty$, the proof of \eqref{eq2} follows from the usual
relation between forwarded differences and derivatives.  For
$1 \le p < \infty$, we only need to consider the case of $i =1$.
Using \eqref{B-B} with $d$ replaced by $d-1$, we obtain by
\eqref{1d-est} that
\begin{align*}
  \left \|\wh \tr^r_{\t \varphi e_1} g \right \|_{p}^p
 & = \int_{\BB^{d-1}} \int_{-1}^1 \left|\wh \tr^r_{\t \varphi(y) \varphi(s)e_1}
     g\left(\varphi(y) s, y \right)\right |^p \, ds\, \varphi(y) dy\\
   &=  \int_{\BB^{d-1}} \int_{-1}^1 |\wh \tr^r_{\t  \varphi(s)}
     g_y(s) |^p \, ds\, \varphi(y) dy\\
& \le c \int_{\BB^{d-1}}  \t^{r p}  \int_{-1}^1\left| \varphi^r(s)
    \frac{d^r}{d s^r} \left[g\left(\varphi(y) s, y \right) \right] \right |^p
     \, ds \, \varphi(y) dy \\
&  =   c\, \t^{r p} \int_{\BB^d} \left| \varphi^r(x) \partial_1^r
g(x) \right |^p  dx
     = c\, \t^{rp} \|\varphi^r \partial_1^r g\|_{p}^p,
\end{align*}
where $g_y(s)=g(\varphi(y)s, y)$. This proves the second
inequality of \eqref{eq2}.

Next, we consider $\Delta^r_{i,j,\t} g$, for which we will need
the corresponding result for trigonometric functions. Let $h$ be a
$2\pi$ periodic function in $L^p[0,2\pi]$ and let $\|h\|_p: =
\left(\int_0^{2\pi} |h(\t)|^p d\t \right)^{1/p}$ in the rest of
this proof. Then it is known (see, for example, \cite{DL}) that
\begin{equation}\label{trig-est}
      \|\overrightarrow \tr^r_h h\|_p \le c h^r \|h^{(r)}\|_p.
\end{equation}
We consider only the case of $(i,j) = (1,2)$. By \eqref{Delta_ij},
\begin{align*}
 & \| \tr^r_{1,2, \t} g\|_{p}^p =  \int_{\BB^{d-2}} \int_{\BB^2}
     \left| \tr^r_{1,2,\t} g\left(v,\varphi(v) u \right) \right |^p
        \varphi(v)^{d-2}  dv \, du  \\
  & =  \int_{\BB^{d-2}} \int_0^1 \rho \int_{0}^{2\pi} \left| \overrightarrow
       \tr^r_{\t} g\left(\rho \cos t, \rho \sin t,\varphi(\rho)
       u \right) \right |^p
         dt   \varphi(\rho)^{d-2} d\rho \,  du.
\end{align*}
Setting $g_{\rho, u} (t)=g\left(\rho \cos t, \rho \sin
t,\varphi(\rho) u \right)$, we deduce from \eqref{trig-est} that
\begin{align*}
  \|\tr^r_{1,2,\t}f \|_{p}^p
  & \le c \, \t^{rp} \int_{\BB^{d-2}} \int_0^1 \rho \int_{0}^{2\pi}
       | g_{\rho,u}^{(r)}(t) |^p
                 dt   \varphi(\rho)^{d-2} d\rho du \\
   &  = c \, \t^{rp} \int_{\BB^d}  \left| D^r_{1,2} g (x) \right |^p
                  dx =c\, \t^{rp}\|D_{1,2}^r f\|_{p}^p,
  \end{align*}
which proves the first  inequality of \eqref{eq2}. Consequently,
we have proved the inequality \eqref{5-23-jun}.

We now prove the reversed inequality
\begin{equation}\label{K=inverse}
  \wh K_r(f,t)_{p}\leq c\, \wh \o_r(f,t)_{p}+c\, t^r\|f\|_p.
\end{equation}
Setting  $n = \lfloor \frac{1}{t} \rfloor$, we have
$$
  \wh K_r(f,t)_p \le \|f - V_n^\mu f\|_p
      + t^r\max_{1 \le i < j \le d} \|D_{i,j}^r V_n^\mu f\|_p +
       t^r\max_{1\le i \le d} \| \varphi^r \partial_i^r V_n^\mu f\|_p.
$$
The first term is bounded by $c\, \wh \o_r(f,t)_p+c n^{-r}\|f\|_p$ by
\eqref{JacksonB2-2}. For the second term, we use \eqref{5-21-jun}
to obtain
\begin{align*}
t^r\max_{1 \le i < j \le d} \|D_{i,j}^r V_n^\mu f\|_p &
    \leq c\, \wh\o_r(V_n^\mu f, n^{-1})_p
    \leq c \,\wh\o_r (V_n^\mu f -f, n^{-1})_p +c\, \wh \o_r (f, n^{-1})_p\\
& \le c \|f-V_n^\mu f \|_p +c \, \wh \o_r (f, n^{-1})_p
   \leq c \, \wh\o_r (f, t)_p+c t^r \|f\|_p.
\end{align*}
The third term can be treated similarly, using Lemma \ref{lem-5-9-Jun}.
This completes the proof of \eqref{K=inverse}. \qed

\subsection{Analogue of Ditzian-Totik modulus of smoothness with
weight} \label{5-5-subsection}

For $w_\mu(x) = (1-t^2)^{\mu-1/2}$, the Ditzian-Totik modulus of
smoothness $\wh \o_r(f,t)_{p,\mu}$ is defined in
\eqref{modu-DT-weight}, with two additional terms of forward and
backward differences close to the boundary, which are shown to be
necessary in \cite[p. 56]{Di-To}.

For the unit ball, we can define the modulus of smoothness with
respect to $W_\mu$ for $\mu > 1/2$ in an analogous way. For this
purpose, we first need to define the analogues of $I_t$ and
$J_{\pm 1,t}$. For $x \in \BB^d$ and $1 \le i \le d$, we define
$\wh x_i :=
 (x_1,\ldots,x_{i-1}, x_{i+1},\ldots, x_d)$. For $ 1 \le i \le d$, we define
$$
  \II_{i,t} : = \left \{x \in \BB^d: \frac{x_i}{\sqrt{1-\|\wh x_i\|^2}}
  \in I_t \right\},\quad
 \JJ_{\pm 1,i,t} : = \left \{x \in \BB^d: \frac{x_i}{\sqrt{1-\|\wh x_i\|^2}} \in J_{\pm 1,t}
  \right\}.
$$

\begin{defn}
For $\mu > 1/2$ and $1 \le p < \infty$, define
\begin{align*}
   \wh \o_r (f,t)_{p,\mu} := & \sup_{|\t| \le t} \max_{1 \le i< j \le d}
      \|  \tr_{i,j,\t}^r f\|_{p,\mu}  +  \sup_{0 < h \le t}
     \max_{1 \le i \le d}  \|\wh \tr_{h \varphi e_i}^r f\|_{L^p(\II_{i,rh},W_\mu)} \\
  & + \sup_{0 < h \le 12r^2t^2} \max_{1\le i \le d} \left(
   \| \overleftarrow \tr_{h e_i}^r f\|_{L^p(\JJ_{1,i,rt},W_\mu)} +
   \| \overrightarrow \tr_{h e_i}^r f\|_{L^p(\JJ_{-1,i,rt},W_\mu)} \right),
    \notag
\end{align*}
with the usual change when $p=\infty$. 
\end{defn}

The direct and the inverse theorems hold for this modulus of
smoothness.

\begin{thm} \label{thm:JacksonB2-Mod-weight}
Let $\mu = \frac{m-1}{2}$ and $m > 2$, let $f \in L^p(\BB^d, W_\mu)$
if $1 \le p < \infty$. Then for $r \in\NN$
\begin{equation}\label{JacksonB2-2-weight}
   E_n (f)_{p,\mu} \le c\, \wh\o_r(f,n^{-1})_{p,\mu} +n^{-r}\|f\|_{p,\mu}
\end{equation}
and
\begin{equation} \label{inverseB2-2-weight}
\wh\o_r(f,n^{-1})_{p,\mu} \leq  c\, n^{-r} \sum_{k=1}^n k^{r-1} E_{k}(f)_{p,\mu}.
\end{equation}
Furthermore, the additional term $n^{-r}\|f\|_{p,\mu}$ on the
right hand side of \eqref{JacksonB2-2-weight}  can be dropped when
$r=1$.
\end{thm}

\begin{proof}
The proof of \eqref{JacksonB2-2-weight} follows along the line of
Theorem \ref{thm:JacksonB2-Mod} and we shall be brief. Since
\eqref{JacksonB} is established for $W_\mu$ with $\mu =
\frac{m-1}{2}$, we again come down to showing that for $1\leq
i\leq d$
\begin{equation}\label{Jackson-B-weight-1}
 \sup_{|\t|\leq t}\|\tr_{i, d+1, \t} ^r \wt f \|_{L^p(\BB^{d+1},W_{\mu-1/2})}
    \leq c \, \wh\o_r(f,t)_{p,\mu} +c t^r \|f\|_{p,\mu},
\end{equation}
where $\wt f(x, x_{d+1})=f(x)$ for $x\in\BB^d$ and $(x,x_{d+1})\in
\BB^{d+1}$. The major difference is that instead of \eqref{tri-wh_tri}, we have
\begin{align*}
 & \int_{B^2} \left |\overrightarrow{\tr}_\t ^r
  f(x_1 \cos (\cdot)  + x_2 \sin (\cdot)) \right|^p (1-x_1^2-x_2^2)^{\mu-1} dx \\
 & \qquad
  \ll  \frac{1}{t} \int_0^t \|\wh \tr_{h \varphi}^r f \|_{L^p(I_{rh},w_\mu)}^p d h
  + \frac{1}{t^2} \int_0^{t^2} \|\overleftarrow \tr_h^r f \|_{L^p(J_{1,1,rt},w_\mu)}^p d h\\
 & \qquad \  +
  \frac{1}{t^2} \int_0^{t^2} \|\overrightarrow \tr_h^r f \|_{L^p(J_{-1,1,rt},w_\mu)}^p d h
    +  t^{rp} \|\varphi^r f\|_{p,\mu}^p, \notag
\end{align*}
which follows from \eqref{omega-omegaDT} and \cite[p. 57]{Di-To}.
Now, it is easy to see that if $x = (\varphi(v) s, v)$, then $s
\in I_{t}$ is equivalent to $x \in \II_{1,t}$ and $s\in
J_{\pm1,t}$ is equivalent to $x \in \JJ_{\pm 1,1,t}$, from which
we can carry out the computation and establish
\eqref{Jackson-B-weight-1} exactly as in the proof of Theorem
\ref{thm:JacksonB2-Mod}.

The proof of \eqref{inverseB2-2-weight} follows again by \eqref{inverseB2}
and Theorem \ref{thm:7-14} below.
\end{proof}

Two remarks are in order. First, it is worth to pointe out that,
in the case of $d =1$, this theorem did not appear in
\cite{Di-To}, which gave the Jackson estimate for the weighted
approximation in terms of the {\it main-part modulus of
smoothness}. The result was later proved in \cite[p. 556]{LuTo} and, for
$1 < p < \infty$, in \cite[Corollary 7.3]{DD2}. Second, in the
case of $p = \infty$, the norm in $\wh \o_r(f,t)_{\infty,\mu}$ is
taken as $\|f\|_{\infty,\mu} = \|W_\mu f\|_{\infty}$, which is not
what the norms in $\o_r(f,t)_{\infty,\mu}$ or $\wh
K_r(f,t)_{\infty,\mu}$ are taken. Thus, we exclude the case of $p
=\infty$ in the above theorem and the theorem below.

\begin{thm}\label{thm:7-14}
Let $\mu = \frac{m-1}{2}$, $m >1$, and let $f \in L^p(\BB^d, W_\mu)$
if $1 \le p < \infty$. Then for $r\in\NN$ and $0<t<t_r$,
\begin{equation} \label{thm7-14-1}
  c^{-1}\wh\o^r(f,t)_{p,\mu}\le  \wh K_r(f,t)_{p,\mu}
                         \leq c\, \wh\o^r(f,t)_p+c\,t^r\|f\|_{p,\mu}.
\end{equation}
Furthermore, the term $t^r\|f\|_p$ on the right side can be dropped
when $r=1$.
\end{thm}

\begin{proof}
The proof of this theorem follows along the same line as that  of
the proof of Theorem \ref{thm:5-5} and we only need to point out
the difference. For the left hand inequality
$\wh\o^r(f,t)_{p,\mu}\le c\, \wh K_r(f,t)_{p,\mu}$ of
\eqref{thm7-14-1}, the counterpart on $\tr_{i,j,\t}^r f$ follows
as in the unweighted case without further complication, so that
the essential part is to show that
\begin{align*}
 & \|\wh \tr^r_{\t \varphi e_i} g\|_{L^p(\II_{i,r\t},W_\mu)} \le c \|f\|_{p,\mu},
 \quad   \|\wh \tr^r_{\t \varphi e_i} g\|_{L^p(\II_{i,r\t},W_\mu)}  \le c\, \t^r
   \| \varphi^r \partial_i^r g\|_{p,\mu}, \\
& \sup_{|\t|\leq 12 r^2t^2}\|\overleftarrow \tr^r_{\t e_i} g\|_{L^p(\JJ_{1,i,rt},W_\mu)} \le c \|f\|_{p,\mu}, \\
&  \sup_{|\t|\leq 12 r^2t^2}
  \|\overleftarrow \tr^r_{\t e_i} g\|_{L^p(\JJ_{1,i,rt},W_\mu)}  \le c\, t^r
   \| \varphi^r \partial_i^r g\|_{p,\mu}
\end{align*}
and a similar inequality for $\overrightarrow \tr^r_{\t e_i} g$
with $g \in C^r(\BB^d)$. As in the proof of Theorem \ref{thm:5-5},
the proof of these inequalities reduces to the corresponding
inequalities in one variable, and the weighted version of
\eqref{1d-est}, which however follow from the results given in
\cite[p. 58]{Di-To}.

For the right hand inequality of \eqref{thm7-14-1}, we can follow the proof
of Theorem \ref{thm:5-5} verbatim once we establish the relation,
for $f \in \Pi_n^d$,
\begin{align*} 
 & n^{-rp}\|\varphi^r \partial_i^r f\|^p_{p,\mu}  \sim n
\int_0^{n^{-1}}
  \|\wh\triangle_{h\varphi e_i}^r f\|_{L^p(\II_{i,rh}, W_\mu)}^p\, dh \\
& \qquad +  n^2 \int_0^{n^{-2}} \|\overleftarrow \tr_h^r f
\|_{L^p(J_{1,i,rn^{-1}},W_\mu)}^p d h
  + n^2 \int_0^{n^{-2}} \|\overrightarrow \tr_h^r f \|_{L^p(J_{-1,i,rn^{-1}},W_\mu)}^p
   d h,
\end{align*}
which is the analogue of \eqref{5-22-jun}. The proof of this relation follows
as that of \eqref{5-22-jun} from the corresponding result in one variable,
and the equivalence in one variable follows from \cite[p. 57]{Di-To} and
\cite[p.193]{HL}.
\end{proof}

It should be mentioned that \cite{Di-To} considers far more
general weight functions than $w_\mu$ in the case of $d =1$, but
we can only deal with $W_\mu$ as our results depend on Subsection
\ref{5-4-subsection}, in which weight is $W_\mu$ with $\mu =
\frac{m-1}{2}$. On the other hand, it is possible to consider
doubling weights and establish the results as in Section 5.

We note that it is more involved to derive properties for the weighted
modulus of smoothness, which requires us to verify that the corresponding
results hold for the weighted $L^p$ space on $[-1,1]$. Such results are
stated mostly for weighted main-part modulus of smoothness in \cite{Di-To}
and a close look at the proof in \cite{Di-To} indicates that the weighted case
requires caution and perhaps further work. Since the result is not needed in
this paper, we shall not pursue it here.

\section{The Weighted $L^p(\BB^d,W_\mu)$ Space with $\mu \ne (m-1)/2$}
\setcounter{equation}{0}

The results that we obtained in the previous sections are
established for the space $L^p(\BB^d, W_\mu)$ with $\mu
=\frac{m-1}{2}$. The definitions of the moduli of smoothness and
the $K$-functionals, however, make sense for all $\mu \ge 0$. A
natural question is if our results can be extended to the case of
$L^p(\BB^d, W_\mu)$ with $\mu \ne \frac{m-1}{2}$. This, however,
appears to be a difficult problem. Below we give a positive result
for the case of $r=2$.

\subsection{Decomposition of $\CD_\mu$}
Recall the second differential operator $D_\mu$ given in
\eqref{D-mu} and the operators $D^2_{i,j}$, $1 \le i < j \le d$
defined in \eqref{partial_ij}. We further define
$$
   D^2_{i,i} := [W_\mu(x)]^{-1} \partial_i \left[(1-\|x\|^2) W_\mu(x) \right]
        \partial_i, \qquad 1 \le i \le d.
$$
It turns out that $D_\mu$ can be decomposed as a sum of second
order differential operators.

\begin{prop}
The differential operator $\CD_\mu$ can be decomposed as
\begin{equation}\label{decomp}
  \CD_\mu = \sum_{i=1}^d D^2_{i,i} +  \sum_{1\le i < j \le d} D^2_{i,j}
      =  \sum_{1\le i \le j \le d} D^2_{i,j}.
\end{equation}
\end{prop}

The proof is a straightforward computation. In the case of $d=2$,
$D^2_{1,2}$ is simply the second partial derivative with respect
to $\t$ in the polar coordinates. In this case, it is tempting to
write the decomposition entirely in terms of polar coordinate $(r,
\t)$ but it does not seem to offer further structure.

The decomposition \eqref{decomp} implies immediately that
$\|\CD_\mu g\|_{p,\mu}$ is bounded by the sum of $\|D^2_{i,j}
g\|_{p,\mu}$ for all $g$ for which the norms involved are finite.
More importantly, however,  the reversed inequality holds. For
this, we relate $\CD_\mu$ with a differential operator,
$\CD_{\mu,T}$ on the simplex $T^d : = \{x \in \RR^d: 1-x_1 -\ldots
- x_d \ge 0, x_i\ge 0, 1\leq i\leq d\}$ and use a result for
$\CD_{\mu,T}$. Let $\BB_+^d: = \{x \in \BB^d: x_i \ge 0, 1 \le i
\le d\}$ and let
\begin{equation}\label{eq:psi}
   \psi: (u_1,\ldots,u_d) \in T^d \mapsto (\sqrt{u_1},\ldots, \sqrt{u_d}) \in \BB_+^d.
\end{equation}
This change of variables leads immediately to the relation
\begin{equation}\label{IntegralBT}
  \int_{\BB^d_+} f(x_1,\ldots,x_d) dx =  \frac{1}{2^d}
     \int_{T^d} f(\sqrt{u_1},\ldots,\sqrt{u_d}) \frac{du}{\sqrt{u_1\cdots u_d}}.
\end{equation}
In particular, it maps the weight function $W_\mu$ to the weight
function
\begin{equation}\label{weightT}
  W_\mu^T (x) = x_1^{-1/2} \ldots x_d^{-1/2} (1-|x|_1)^{\mu-1/2},
     \quad   x \in T^d,
\end{equation}
where $|x|_1 = x_1 + \ldots +x_d$.  Furthermore, the mapping
\eqref{eq:psi} sends the differential operator $\CD_\mu$ to
\begin{align} \label{D-muT}
   \CD_{\mu, T}: = \sum_{i=1}^d x_i(1-x_i)  \partial^2_i -
      2 \sum_{1 \le i < j \le d} x_ix_j  \partial_i \partial_j
       + \sum_{i=1}^d \left ( \tfrac12 - \left(\mu+ \tfrac{d+1}{2}\right)
          x_i \right) \partial_i,
\end{align}
and $\CD_{\mu, T}$ has orthogonal polynomials with respect to
$W_\mu^T$ on $T^d$ as eigenfunctions. Much of the analysis on
$\BB^d$ or $T^d$ can be carried over to the other domain through
this connection (see, for example, \cite{X06}). It is known that
$\CD_{\mu,T}$ satisfies a decomposition \cite{BSX,Di},
$$
  \CD_{\mu,T} = \sum_{i=1}^d U_{i,i}^T +  \sum_{1\le i < j \le d} U_{i,j}^T
       = \sum_{1 \le i \le j \le d} U_{i,j}^T,
$$
where, with $\partial_{i,j} : = \partial_i - \partial_j$,
\begin{align*}
 U_{i,i}^T  & = [W_\mu^T(x)]^{-1} \partial_i \left[x_i (1-|x|)
    W_\mu^T(x) \right]   \partial_i, \qquad 1 \le i \le d, \\
  U_{i,j}^T & = [W_\mu^T(x)]^{-1} \partial_{i,j} \left[x_i x_j
    W_\mu^T(x) \right]  \partial_{i,j}, \qquad 1 \le i < j \le d.
\end{align*}
Let $\|\cdot \|^T_{p,\mu}$ denote the norm of $L^p(W_\mu^T) =
L^p(T^d, W_\mu^T)$. In fact, the decomposition of $\CD_\mu$ can
also be derived from the mapping \eqref{eq:psi}.

\begin{lem}
For $g \in C^2(\BB^d)$,  $1 \le i \le j  \le d$ and $u\in T^d$,
\begin{equation} \label{U:T-B}
  \CD_\mu g(\psi(u)) = 4 \CD_{\mu,T}(g\circ \psi)(u) \quad \hbox{and}\quad
     D^2_{i,j} g(\psi(u)) = 4 U_{i,j}^T (g \circ \psi)(u).
\end{equation}
\end{lem}

\begin{proof}
Under the change of variables $x_i = \sqrt{u_i}$, $1 \le i \le d$
so that $(g \circ \psi)(u) = g(\sqrt{u_1},\ldots,
\sqrt{u_d})=g(x)$, the relation for $\CD_\mu$ follows from a
straightforward computation and so is the case $D^2_{i,i}$, since
$x_i = \sqrt{u_i}$ implies that $\partial_{x_i} = 2
\sqrt{u_i}\partial_{u_i}$. We now consider the case of $D^2_{i,j}$
with $i < j$. First we note that
$$
   D^2_{i,j} = [W_\mu(x)]^{-1}  \left( x_i \partial_{x_j} - x_j  \partial_{x_i}\right)
      W_\mu(x) \left( x_i \partial_{x_j} - x_j \partial_{x_i}\right).
$$
In fact, the above identity holds if any differentiable radial
function is in place of $W_\mu$. Setting $x_i = \sqrt{u_i}$, we
see easily that  $x_i \partial_{x_j} - x_j \partial_{x_i} = 2
\sqrt{u_i u_j} (\partial_{u_j} - \partial_{u_i})$. Consequently,
it follows that
\begin{align*}
  D^2_{i,j}g (x) & = 4 (1-|u|_1)^{-\mu+\frac12} \sqrt{u_iu_j}
       \left(\partial_{u_j} - \partial_{u_i}\right)  \sqrt{u_i u_j} \\
  & \qquad \qquad \qquad \times (1-|u|_1)^{\mu-\frac12}
          \left(\partial_{u_j} - \partial_{u_i}\right)(g\circ \psi)(u) \\
    & = 4  \left[W_\mu^T(u)\right]^{-1}   \left(\partial_{u_j} - \partial_{u_i}\right)
         \left[u_i u_j W_\mu^T(u) \right]
            \left(\partial_{u_j} - \partial_{u_i}\right) = 4 U_{i,j}^T (g\circ \psi)(u),
\end{align*}
which verifies \eqref{U:T-B}.
\end{proof}

\subsection{Differential operators and $K$-functional}
The following result was established in \cite{DHW} recently: for
$f \in C^2(T^d)$,
\begin{equation} \label{decompT}
 \|\CD_{\mu,T}  f\|_{L^p(T^d,W_\mu^T)} \sim \sum_{1 \le i \le j \le d}
       \|U_{i,j}^T f \|_{L^p(T^d,W_\mu^T)},   \quad 1 < p < \infty.
\end{equation}
With the connection in the previous subsection, it immediately
leads us to an analogous result for $\CD_\mu$ on $\BB^d$.

\begin{thm} \label{thm:decom}
For $g \in C^2(\BB^d)$,
$$
   \|\CD_\mu  g\|_{p,\mu} \sim \sum_{1 \le i \le j \le d} \|D^2_{i,j}g \|_{p,\mu},
     \quad 1< p < \infty.
$$
\end{thm}

\begin{proof}
By \eqref{decomp}, it suffices to prove that
\begin{equation}\label{8-8-D}
\|D^2_{i,j} g\|_{p,\mu} \le c \| \CD_\mu g\|_{p,\mu}\end{equation}
for $1 \le i,j \le d$. This is the same as in the proof in
\cite{DHW}, which amounts to showing that for $f\in C^2(T^d)$ and
$1\leq i<j\leq d$,
\begin{equation} \label{estimateT}
   \|U_{i,j}^T f \|_{L^p(T^d,W_\mu^T)}
   \le c  \|\CD_{\mu,T}  f\|_{L^p(T^d,W_\mu^T)},
   \quad 1 < p < \infty.
\end{equation}
We shall deduce \eqref{8-8-D} from \eqref{estimateT} and the
change of variables \eqref{eq:psi}. Now, for $\ve=(\va_1, \cdots,
\va_d)\in\{-1,1\}^d$, we define $g_\ve$ by $g_\ve(x) = g(x_1
\ve_1,\ldots, x_d \ve_d)$, $x\in \BB^d$. We claim that
\begin{equation} \label{U:comm}
    (D^2_{i,j} g)_\ve (x) =  (D^2_{i,j} g_\ve) (x), \qquad 1 \le i \le j \le
    d,\   \   x\in \BB^d.
\end{equation}
 Indeed,  since $\partial_i g_\ve (x) = \ve_i (\partial_i g)_\ve(x)$, this can be
verified via a straightforward computation. Using
\eqref{IntegralBT}, \eqref{U:T-B} and \eqref{U:comm}, we obtain
\begin{align*}
   \int_{\BB^d} |D^2_{i,j} g (x)|^p W_\mu(x) dx   &=  \sum_{\ve \in \{-1,1\}^d}
        \int_{\BB^d_+}  |(D_{i,j}^2 g)_\ve (x)|^p W_\mu(x) dx  \\
   &   =  \sum_{\ve \in \{-1,1\}^d}
        \int_{\BB^d_+}  |D_{i,j}^2 g_\ve (x)|^p W_\mu(x) dx \\
    & =  \frac{4^p}{2^d} \sum_{\ve \in \{-1,1\}^d}
             \int_{T^d} | U_{i,j}^T (g_\ve \circ \psi)(u)|^p W_\mu^T(u) du.
\end{align*}
Consequently, using \eqref{estimateT}, followed by using
\eqref{U:T-B} and \eqref{U:comm} again, we conclude that
\begin{align*}
   \int_{\BB^d} |D^2_{i,j} g (x)|^p W_\mu(x) dx & \le
       c \frac{4^p}{2^d} \sum_{\ve \in \{-1,1\}^d} \int_{T^d}
           | \CD_{\mu,T} (g_\ve \circ \psi) (u)|^p W_\mu^T(u) du \\
      &  \le   c \frac{1}{2^d} \sum_{\ve \in \{-1,1\}^d} \int_{\BB_+^d}
           | (\CD_{\mu} g)_\ve  (x)|^p W_\mu(x) dx \\
       & =   c \int_{\BB^d}  |\CD_{\mu} g (x)|^p W_\mu(x) dx
\end{align*}
for $1 < p < \infty$. This completes the proof.
\end{proof}

The differential operators $D^2_{i,j}$ for $i < j$ are second
order derivatives with respect to the Euler angles, whereas
$D^2_{i,i}$ does not have such simple interpretation. Our next
result shows that $\|D_{i,i}^2 g\|_{p,\mu}$ can be further
reduced. Recall $\varphi (x)= \sqrt{1-\|x\|^2}$.

\begin{thm} \label{thm:phi^2}
For $1 < p < \infty$, $1 \le i \le d$ and $g \in C^2(\BB^d)$,
$$
  c_1 \| \varphi^2 \partial_{i}^2 g\|_{p,\mu} \le \|D^2_{i,i} g\|_{p,\mu}
      \le  c_2 \| \varphi^2 \partial_{i}^2 g\|_{p,\mu} + c_2 \|g\|_{p,\mu}.
$$
\end{thm}

\begin{proof}
It is enough to consider $D^2_{1,1}$.  We make a change of
variables, $x \mapsto (s, y)$, where $y = (y_2,\ldots, y_d)$,  by
setting $x_1 = \sqrt{1-\|y\|^2} s$ and $x_i = y_i$ for $i
=2,\ldots, d$. It follows immediately that $\varphi(x) =
\varphi(y) \varphi(s)$. Furthermore, a quick computation shows
that
$$
   D^2_{1,1} g(x) = \Lambda_s g_y (s), \quad
   \Lambda_s: = (1-s^2)\frac{d^2}{d s^2} - (2\mu+1) s
   \frac{d}{ds},
$$
where $g_y(s)=g(s \varphi(y) , y)$.  Let $w_{\mu}(s) =
(1-s^2)^{\mu -1/2}$. It is easy to see then that $\Lambda_s$ can
be written as
$$
     \Lambda_s = w_{\mu}(s)^{-1} \frac{d}{ds}
              \left[ (1-s^2) w_{\mu}(s)\right] \frac{d}{ds}.
$$
It is known that the differential operator $\Lambda_s$ satisfies
(\cite{DHW})
\begin{equation}\label{Dt}
c_1 \|\varphi^2 g''\|_{L^p(w_{\mu})} \le  \| \Lambda_s
g\|_{L^p(w_{\mu})}
    \le  c_2 \|\varphi^2 g''\|_{L^p(w_{\mu})} + c_2 \|g\|_{L^p(w_{\mu})}
\end{equation}
for $1 < p < \infty$, where the norm is taken over $[-1,1]$. By
\eqref{B-B}, applying \eqref{Dt} we obtain
\begin{align*}
  \|D^2_{1,1}g\|_{p,\mu}^p  & =
    \int_{\BB^{d-1}}  \int_{-1}^1 \left| \Lambda_s g_y (s)\right|^p
         w_\mu(s) ds W_{\mu + \frac12}(y) dy \\
& \ge c_1^p \int_{\BB^{d-1}}  \int_{-1}^1 \left | \varphi^2(s)
\frac{d^2}{ds^2}
         \left[ g \left(s \varphi(y), y\right)\right] \right |^p w_\mu(s) ds
           W_{\mu + \frac12}(y) dy \\
& =  c_1^p \int_{\BB^{d-1}}  \int_{-1}^1 \left |
\varphi^2(s)^2\varphi^2(y)
          (\partial_1^2 g) \left (s \varphi(y), y \right )  \right |^p w_\mu(s) ds
           W_{\mu + \frac12}(y) dy \\
& = c_1^p \int_{\BB^d} \left | \varphi^2(x)
          \partial_1^2 g(x) \right |^p W_{\mu}(x) dx
             = c_1^p \|\varphi^2 \partial_1^2 g\|_{p,\mu},
\end{align*}
where in the last step we used  \eqref{B-B} again. This proves the
the lower bound. The upper bound is proved likewise.
\end{proof}

As a consequence of this theorem, we can replace $\|\CD_\mu
g\|_{p,\mu}$ in $K$-functional $K_2^*(f;t)_{p,\mu}$, defined in
\eqref{K-funcB-1} by the sum of $\|\varphi^2 \partial_i^2
g\|_{p,\mu}$ and $\|D_{i,j}^2 g\|_{p,\mu}$, at least for $1 < p <
\infty$, which leads to a comparison between $K_r^*(f,t)_{p,\mu}$
and $\wh K_r(f,t)_{p,\mu}$ defined in Definition
\ref{defn:K-funcB2}. Indeed, from Theorems \ref{thm:decom} and
\ref{thm:phi^2}, and the triangle inequality, we obtain the
following:

\begin{thm} \label{thm:whK-K*}
For $f \in L^p(\BB^d,W_\mu)$,  $1 < p < \infty$, and $0 < t \le
1$,
\begin{equation}\label{whK-K*}
 c_1\wh K_2(f;t)_{p,\mu} \le  K_2^*(f;t)_{p,\mu} \le  c_2 \wh K_2(f;t)_{p,\mu} + c_2
             t^2 \|f\|_{p,\mu}.
\end{equation}
\end{thm}

Recall that both direct and inverse theorems for best polynomial
approximation are established for $K_r^*(f,t)_{p,\mu}$, the above
shows that the same can be stated  for $\wh K_2(f,t)_{p,\mu}$.

In the case of $\mu = \frac{m-1}{2}$, $K_2^*(f,t)_{p,\mu}$ is
equivalent to $K_2(f,t)_{p,\mu}$ for $1 < p < \infty$ by the
$K$-functional counterpart of Theorem \ref{thm:moduliEqu-B}, which
gives the inequality \eqref{whK=K-r=2}.


\part{Computational Examples}
In this part we give examples of functions for which the asymptotic orders
of our new moduli of smoothness and best approximation by polynomials
are explicitly determined. The first section contains a lemma, upon which
most of the computations of our examples are based, and one of its applications.
The examples for moduli of smoothness on the sphere are given in Section 9
and examples on the ball are given in Section 10.


\section{Main lemma for computing moduli of smoothness}
\setcounter{equation}{0}

One of the advantages of our new moduli of smoothness lies in the fact
that the divided difference in Euler angles can be reduced to the forwarded
difference for trigonometric functions (cf. \eqref{Delta_ij}), which are
classical and well studied. Our claim that the new moduli of smoothness
are computable is based on this fact. Below we present a lemma that
gives the asymptotic order of the modulus of smoothness for a simple
trigonometric function, upon which most of our examples in the following
two sections are based.

\begin{lem}\label{lem-0-1}
Assume that  $ |a| <1$,  $1\leq p\leq \infty$ and $\a\neq 0$. If
$\a\neq \f12$, then there exists $\d_\a\in (0,1)$ depending only
on $\a$ such that for $|\t|\leq \d_\a$,
\begin{align} \label{ex-basic}
  & \left( \int_0^{2\pi}   \left| \overrightarrow\Delta_\t^2
  \bigl(1-a \cos (\phi + \{\cdot\})\bigr)^\a \right|^p d\phi  \right)^{1/p} \\
   &  \   \sim |a| \t^{2}
     \begin{cases}  (1- |a| +|a|\t^2)^{(\a -1)+\frac{1}{2p}},
      & \a < 1-\frac1{2p}, \\
           \Bigl|\log (1-|a| +|a|\t^2)\Bigr|, & \a =1-\frac1{2p}, \\
            1, & \a > 1-\frac1{2p},
      \end{cases}  \notag
\end{align}
with the usual change of the $L^p$-norm when $p=\infty$, where the
constants of equivalence are independent of $a$. If $\a=\f12$, then
the upper estimates in \eqref{ex-basic} remain true, whereas the lower
estimate holds under the additional condition $36\t^2 \leq 1-|a|$ when
$p >1$, and the lower bound becomes $c \t^2$ when $p=1$.
\end{lem}

\begin{proof} Without loss of
generality,  we may  assume $a
>0$, since $\cos (\pi+ \phi) = - \cos \phi$. We shall prove
the lemma for $p<\infty$ only. The case $p=\infty$ can be treated
similarly, and in fact, is simpler. Let us set $h_\a(\phi) = (1- a
\cos \phi)^\a$. We will use the fact that $\overrightarrow
\tr_\t^2 h_\a(\phi) = (1/2) \t^2 h_\a''(\phi+\xi)$ for some $\xi$
between $0$ and $2 \t$, and
\begin{align} \label{exam2-1}
  h_\a''(\phi+\xi )=
     & \a(\a-1) (1- a \cos (\phi+ \xi))^{\a-2} a^2 \sin^2 (\phi+ \xi) \\
     &    + \a(1- a \cos (\phi+ \xi))^{\a-1} a \cos (\phi+ \xi).   \notag
\end{align}

To show   the upper estimates, we can restrict the integral in
  \eqref{ex-basic} over  $[0,\pi]$  instead of $[0,2\pi]$, since
  we allow
  $\t$ to take negative values,  and
 $\overrightarrow\Delta_\t^2
h_\a(-\phi+\{\cdot\})=\overrightarrow\Delta_{-\t}^2
h_\a(\phi+\{\cdot\})$.   Using \eqref{exam2-1} and the identity $1
- a \cos \psi = 1- a + 2a \sin^2\frac{\psi}{2}$,  we have
\begin{align} \label{exam2-2}
 |\overrightarrow\tr_\t^2  h_\a(\phi) |=\f12\t^2  \left|
   h_\a''(\phi+\xi)\right| \le c\, a\, \t^2 \left(1-a +
           2 a \sin^2 \tfrac{\phi+ \xi}{2}\right)^{\a-1}.
\end{align}
We break the integral of $|\overrightarrow\Delta_\t^2 h_\a(\phi+\{\cdot\})|^p$
into two parts:
$$
 \int_0^\pi |\overrightarrow\Delta_\t^2 h_\a(\phi+\{\cdot\})|^p d\phi
   =   \int_{0}^{3|\t|}\cdots+\int_{3|\t|}^\pi\cdots=: I_1+I_2.
$$
If $ a\t^2\ge 1-a$, then $a\ge 1-\t^2\ge a_0>0$, and   using  the
definition of  $
  \overrightarrow\Delta_\t^2$, we obtain
$|\overrightarrow \tr_\t^2 h_\a (\phi+\{\cdot\})|  \le c (1-a
+a\t^2)^\a$ for $|\phi|\leq 3|\t|$, which in turn implies
$$ I_1 \leq C \int_0^{3|\t|}(1-a
+a\t^2)^{\a p}\, d\phi\sim |\t|^{2\a p +1} a^{\a p} \sim a^{p}
|\t|^{2p} (1-a+a\t^2)^{(\a-1) p+\f12}. $$ On the other hand,  if $
a\t^2\leq 1-a$ then using \eqref{exam2-2}, we have
$\left|\overrightarrow\Delta_\t^2 h_\a(\phi+\{\cdot\})\right| \le
c a \t^2 (1-a)^{\a-1 }$ for $|\phi|\leq 3|\t|$, and hence
\begin{align*}
 I_1   & \le  c \,a^p |\t|^{2p+1} (1-a)^{(\a-1) p} \leq
c a^{p} |\t|^{2p} (1-a+a\t^2)^{(\a-1) p+\f12},
\end{align*}
where the last step uses the fact that $ 1-a\sim 1$
when $0<a\leq \f12$.  Finally,  using \eqref{exam2-2}, we deduce
\begin{align*}
 I_2   & \le  c \,a^p \t^{2p}\int_{3\t}^\pi (1-a + a \phi^2)^{(\a-1)p} d\phi
 \sim a^{p-\f12} \t^{2p}\int_{3\sqrt{a}\t}^{\sqrt{a}\pi }
   (\sqrt{1-a} +  \phi)^{2(\a-1)p} d\phi,
\end{align*}
which is estimated by the desired upper  bounds. This
completes the proof of the upper estimates.

For  the proof of the lower estimates,  we shall use $\d_\a'$ or
$\d_\a''$ to denote a sufficiently small positive constant which
depends only on $\a$, and may vary at each occurrence.
Note that if $1-a\ge \d_\a'>0$ then the desired lower estimates follow
immediately since, \mbox{ by \eqref{exam2-1}},
$$
 \left|\overrightarrow\tr_\t^2 h_\a(\phi)\right|=\f12\t^2  \left|
   h_\a''(\phi+\xi)\right|\ge c_\a a \t^2
$$
whenever  $ 3|\t|\leq \phi\leq \d_\a'$. Thus, for the rest of the proof, we
may assume that $|\t|+\sqrt{1-a}\leq \d_\a''$.
We claim that for  $\a\neq \f12$,  there exists a constant $c_\a>2$ such that
\begin{align} \label{exam2-31}
 \left|\overrightarrow\tr_\t^2  h_\a(\phi) \right | = \f12\t^2  \left|
   h_\a''(\phi+\xi)\right| \ge c\, \t^2 \phi^{2\a-2}
\end{align}
whenever $c_\a ( |\t|+ \sqrt{1-a})\leq \phi \leq \d_\a'$. Indeed,
  setting  $\psi=\phi +\xi$, and   using  \eqref{exam2-1}, we obtain,
   for $c_\a ( |\t|+ \sqrt{1-a})\leq \phi \leq \d_\a'$,
\begin{align*}
 \left|h_\a''(\psi)\right | & =a |\a| \left( 1-a+a\sin^2 \tfrac \psi 2\right)^{\a-2}
 \left |\a a \sin^2\psi -a+\cos\psi \right |\\
 &=a |\a| \left( 1-a+a\sin^2 \tfrac \psi 2\right )^{\a-2} \left  (
 \a-\tfrac12)\psi^2 + O(1-a) + O (\psi^3)\right|\\
 &\ge c |\a| |\a-\tfrac12| \phi^{2\a-2}
\end{align*}
provided that $(1- \frac{2}{c_\a}) \ge c \d_\a'$. The assertion
\eqref{exam2-31} then follows.
Now raising \eqref{exam2-31} to the power $p$, and
integrating it with respect to $\phi$ over $c_\a ( |\t|+
\sqrt{1-a})\leq \phi \leq \d_\a'$ gives the desired lower
estimates in \eqref {ex-basic} for $\a\neq \f12$.  The lower
estimate  for $\a=\f12$ can be proved similarly. Indeed, setting
$\psi=\phi+\xi$, and using \eqref{exam2-1}, we obtain
$$
h_{1/2}''(\psi) = \tfrac 14 ( 1-a\cos \psi)^{-\f32} a \Bigl( 2\cos
\psi -a\cos^2\psi-a\Bigr)\ge \tfrac 14 ( 1-a\cos \psi)^{-\f12}
a\cos\psi,
$$
provided that  $1-a\ge 2\sin^2 \f \psi2$. Thus, if
$3|\t|\leq \phi \leq \sqrt{(1-a)/2} $ then $ 1-a\ge 2\phi^2 \ge \f
12 \psi^2 \ge 2\sin^2\f \psi2$, and hence
\begin{align*}
 \left|\overrightarrow\tr_\t^2  h_{1/2}(\phi)\right | = \tfrac12\t^2  \left|
   h_{1/2}''(\phi+\xi)\right| \ge c|a| \t^2 (1-a)^{-\f12}.
\end{align*}
Integrating the $p$-th power of this inequality with respect to $\phi$
over $3|\t|\leq \phi \leq \sqrt{(1-a)/2}$ and using $1-a+a |\t| \sim 1-a$,
which holds for $c \t^2 \le 1-a$, give the desired lower estimate for
$\a=\f12$.
\end{proof}

As an application of Lemma \ref{lem-0-1}, we prove the asymptotic
of $\o_2(g_\a,t)_{p,\mu}$ in Example 5.11. The proof also suggests
what to come in the following two sections.

\begin{lem}\label{lem-4-3}
Let $h_\a(s,\phi): = (1-s \cos \phi)^\a$, $\a \ne 0$, $\mu> 0$ and
let
\begin{equation}\label{interal-4-3}
 \Omega_2(h_\a, \t)_{p,\mu}:= \left( \int_0^1 s  \int_0^{2\pi}
     \left |\overrightarrow \tr_\t^2 h_\a(s,\phi) \right|^p d \phi
    (1- s^2)^{\mu-1}ds \right)^{1/p},
\end{equation}
where for $ p = \infty$ it is defined as maximum of $ |
\overrightarrow \tr_\t^2 h_\a(s,\phi) |$ over $0 \le s \le 1, 0
\le \phi \le 2\pi$. Then there is a $t_0 > 0$ such that for $0 < \t
< t_0$, $\mu \ge 0$ and $1 \le p \le \infty$,
\begin{align}\label{eq:exam2*}
     \Omega_2(h_\a, \t)_{p,\mu} \sim \begin{cases}
       |\t|^{2 \a + \frac{2 \mu + 1}{ p}}, & - \frac{2\mu+1}{2p} < \a < 1 - \frac{2\mu+1}{2p},\\
       \t^{2} \Bigl|\log |\t| \Bigr|^{1/p}, &  \a = 1 - \frac{2 \mu+1}{2p}, \quad p \ne \infty,\\
       \t^{2},  &  \a > 1 - \frac{2 \mu +1}{2p}.
    \end{cases}
\end{align}
\end{lem}

\begin{proof}
Again we only consider $1 \le p < \infty$.  If $ \a < 1 -
\frac{1}{2p}$ and $\a\neq \f12$,  then  we apply \eqref{ex-basic}
with $a = s$ to obtain that
\begin{align*}
  \Omega_2(h_\a,\t)_p^p &
  \sim |\t|^{2p} \int_0^1 s^{p+1} (1-s+s\t^2)^{(\a-1)p+\frac12}
     (1-s^2)^{\mu - 1} ds\\
     & \sim |\t|^{2p} \left[\int_0^{1-\t^2} s^{p+1} (1-s)^{(\a-1)p+\mu-\frac12}
      ds+
     \int_{1-\t^2}^1  |\t|^{(2\a-2)p+1}
     (1-s)^{\mu - 1} ds\right],
\end{align*}
which, when integrated out according to $(\a-1) p + \mu -1/2 <
-1$, $= -1$, and $> -1$, is easily seen to be equivalent to the
$p$-th power of the right hand side of \eqref{eq:exam2*}.
Similarly, by \eqref{ex-basic} applied to  $a = s$, we have, for
$\a=1-\f1{2p}\neq \f12$,
\begin{align*}
  \Omega_2(h_\a,\t)_p^p &
  \sim |\t|^{2p} \int_0^1 s^{p+1} |\log (1-s+s\t^2)|
     (1-s^2)^{\mu - 1} ds\\
     & \sim |\t|^{2p} \left[\int_0^{1-\t^2} s^{p+1} (1-s)^{\mu-1}\log
     \f 1{1-s} \, ds +\int_{1-\t^2}^1  (1-s)^{\mu-1}\Bigl|\log
     |\t|\Bigr|\, ds\right]\\
     &\sim |\t|^{2p},
\end{align*}
whereas for $\a>1-\f1{2p}$,
$$
  \Omega_2(h_\a,\t)_p^p   \sim |\t|^{2p} \int_0^1 s^{p+1}
     (1-s^2)^{\mu - 1} ds\sim |\t|^{2p}.
$$
Finally, in the case when  $ \a=\f12$,  the desired upper estimate
for  $\Omega_2(h_{1/2},\t)_p^p$ can be obtained exactly as above,
while the lower estimate for $\Omega_2(h_{1/2},\t)_p^p$ can be
obtained using the second statement of Lemma \ref{lem-0-1}:
\begin{align*}
  \Omega_2(h_{1/2},\t)_p^p&\ge \int_0^{1-36\t^2} s
\int_0^{2\pi}
     \left |\overrightarrow \tr_\t^2 h_{1/2}(s,\phi) \right|^p d \phi
    (1- s^2)^{\mu-1}ds\\
    &\ge c |\t|^{2p} \int_0^{1-36\t^2} s^{p+1}(1-s)^{-\f p2-\f12+\mu}\,
    ds,
\end{align*}
which, by an easy calculation,  gives  the desired lower estimate
for the case of $\a=\f12$.
\end{proof}

\begin{rem} \label{rem8.1}
In the previous two lemmas we considered only the second order
difference. Our proof can be adopted to give the upper estimates
for $r > 2$. For examples, let $r \ge 2$ and define
$$
 \Omega_r(h_\a, \t)_{p,\mu}:= \left( \int_0^1 s  \int_0^{2\pi}
     \left |\overrightarrow \tr_\t^r h_\a(s,\phi) \right|^p d \phi
    (1- s^2)^{\mu-1}ds \right)^{1/p}
$$
for $1 \le p < \infty$ and the usual convention for $p =\infty$. Then we can show
that
\begin{align*}
     \Omega_r(h_\a, \t)_{p,\mu} \le  c \begin{cases}
       |\t|^{2 \a + \frac{2 \mu + 1}{ p}}, & - \frac{2\mu+1}{2p} < \a < \frac{r}{2} - \frac{2\mu+1}{2p},\\
       \t^{2} \Bigl|\log |\t| \Bigr|^{1/p}, &  \a =  \frac{r}{2} - \frac{2 \mu+1}{2p}, \quad p \ne \infty,\\
       \t^{2},  &  \a >  \frac{r}{2} - \frac{2 \mu +1}{2p}.
    \end{cases}
\end{align*}
Although we believe that the lower estimate should also hold, it
is much more difficult to establish. For this reason, we only
considered $r =2$.

For the same reason and because the computation for $r=2$ is
already rather involved, we shall consider only $r =2$ in most of
our examples in the next two sections. In all cases, our method
can be adopted to establish the upper estimates for all $r\ge 1$.
\end{rem}

\section{Computational Examples on the unit sphere}
\setcounter{equation}{0}

In this section we compute the modulus of smoothness $\o_r(f,t)_p =
\o_r(f,t)_{L^p(\SS^{d-1})}$ defined in \eqref{eq:modulus} and the best
approximation $E_n(f)_p : = E_n(f)_{L^p(\SS^{d-1})}$ of \eqref{eq:bestEn}.

\subsection{Computation of moduli of smoothness}

We start with a simple example that follows directly from the modulus
of smoothness for trigonometric functions.

\begin{exam} \label{example1}
For $x \in \SS^{d-1}$ and $d \ge 3$, let $f_\a(x) = x^\a$ with $\a
= (\a_1,\ldots,\a_d)\neq 0$. If $0 \le \a_i < 1$ for $1 \le i \le
d$, then for $r \ge 2$ and $1\leq p\leq \infty$,
\begin{equation}\label{eq:exam1}
 \o_r (f,t)_{L^p(\SS^{d-1})}\sim t^{\d + 1/p}, \qquad \d = \min_{\a_i \ne 0} \{\a_1,\ldots,\a_d\}.
\end{equation}
\end{exam}

Indeed, we only need to consider $\tr_{1,2,\t}^2 f$, which, by
\eqref{Delta_ij}, can be expressed as a forward difference
$$
     \tr_{1,2,\t}^r f_\a(x)    =   x_3^{\a_3} \cdots x_d^{\a_d} s^{\a_1+\a_2}
        \overrightarrow \tr_\t^r  \left[(\cos \phi)^{\a_1} (\sin \phi)^{\a_2}\right]
$$
where $(x_1,x_2) = (s \cos \phi, s\sin \phi)$. Hence, by
\eqref{Int_m=2} we obtain
\begin{align*}
  \|\tr_{1,2,\t}^r f_\a \|_p = c  \left(\int_0^{2\pi}
       \left| \overrightarrow \tr_\t^r  \left[(\cos \phi)^{\a_1} (\sin \phi)^{\a_2}\right]   \right|^p
          d\phi \right)^{1/p}.
\end{align*}
Furthermore, using the well known relation
$$
   \overrightarrow \tr_\t^r (f g)(\phi) =  \sum_{k=0}^r \binom{n}{k}
   \overrightarrow \tr_\t^k f (\phi)
        \overrightarrow \tr_\t^{r-k} g(\phi + k \t),
$$
we can consider the differences for $\cos (\phi+\cdot)$ and $\sin
(\phi+\cdot)$ separately. Since the sine and cosine functions
cannot be both large or both small, we can divide the integral
domain accordingly and estimate the integral in the $L^p$ norm.
Furthermore, in our range of $\a_i$, we only need to consider the
second difference ($r=2$) upon using (1) of Proposition \ref{modulus}.
The equation \eqref{eq:exam1} also holds for $r =1$ and $p = \infty$.

Our second example is more interesting and appears to be
non-trivial.

\begin{exam} \label{example2}
For $d \ge 3$ and $\a\neq 0$, let $g_\a(x) = (1-x_1)^\a$,  $x =
(x_1,\ldots, x_d) \in \sph$. Then for $1 \le p \le \infty$,
\begin{equation}\label{eq:exam2}
    \o_2 (g_\a,t)_{L^p(\SS^{d-1})}  \sim \begin{cases}
       t^{2 \a + \frac{d-1}{p}}, & - \frac{d-1}{2p} < \a < 1 - \frac{d-1}{2p}
       , \\
       t^{2} |\log t |^{1/p}, &  \a = 1 - \frac{d-1}{2p}, \quad p \ne \infty,\\
       t^{2},  &  \a > 1 - \frac{d-1}{2p}.
      \end{cases}
\end{equation}
For $\a =0$, $\o_2 (g_\a,t)_p =0$.
\end{exam}

If neither $i$ nor $j$ equals to $1$, then $\tr^2_{i,j,\t} g_\a(x)
=0$. Thus, we only need to consider $\tr^2_{1,j,\t} g_\a$ and we
can assume $j =2$.  Since $x \in \sph$ and $d \ge 3$ imply that
$(x_1,x_2) \in \BB^2$, by \eqref{Delta_ij},
\begin{align*}
 \| \tr_{1,2,\t}^2 g_\a \|_p^p & = c\int_{B^2} |\tr_{1,2,\t}^2 g_\a(x_1,x_2)|^p
      (1-x_1^2-x_2^2)^{\mu-1} dx \\
 & = c\int_0^1s  \int_0^{2\pi}  \left|  \overrightarrow \tr_\t^2 (1- s \cos \phi)^\a
   \right |^p  d\phi (1-s^2)^{\mu-1} ds,
\end{align*}
where $\mu  = \frac{d-2}{2}$ and the forward difference acts on
$\phi$; for $p = \infty$ the integral is replaced by the maximum
taken over $0 \le s \le 1$ and $0\le \phi \le 2 \pi$. Hence,
\eqref{eq:exam2} follows from Lemma \ref{lem-4-3}.

Our next example is more complicated and it should be compared
with \eqref{eq:exam2}. In particular, we note that its asymptotic
order is independent of $d$, in contrast to the order in
\eqref{eq:exam2}.

\begin{exam} \label{ex:sphere-2}
Let $d\ge 3$ and let  $f(x)=(x_1^2+x_2^2)^\a$ for $x\in \SS^{d-1}$
and $\a\neq 0$. Then for $1\leq p\leq \infty$,
$$
 \omega_2(f,t)_{L^p(\SS^{d-1})} \sim\begin{cases}
t^{2\a +\f 2p},&\   \  \text{if $-\f1p<\a<1-\f 1p$;}\\
 t^2 |\log t|^{\f 1p}, & \   \  \text{if $\a=1-\f 1p $;}\\
 t^2, & \   \  \text{if $\a>1-\f 1p$.}
\end{cases}
$$
\end{exam}

\begin{proof}
 Since $\tr^2_{i,j,\t} f(x)=0$ if $(i,j)=(1,2)$ or $ 3\leq i<j\leq d$,
 it  suffices to consider $ \tr_{1,3,\t}^2 f(x)$.
 For a fixed $x\in \SS^{d-1}$, let $ g_x(t) = f(Q_{1,3,t} x).$
Clearly, $ g_x(t) = ( v(t) ^2 +x_2^2)^\a$, where
 $v(t)\equiv v_x(t)=x_1\cos t -x_3 \sin t$. A straightforward
 computation shows that
 \begin{align} \label{0-1-30}
  g_x''(t) = & 4\a (\a-1) \bigl( v(t)^2 +x_2^2\bigr)^{\a-2} (v(t) v'(t))^2 \\
   & +2\a \left(v(t)^2+x_2^2\Bigr)^{\a-1} \Bigl[ ( v'(t))^2 + v(t) v''(t)\right]. \notag
\end{align}

Let us  start  with the proof of the lower estimate.
 Setting $c_\a=8(1+|\a|)$ and
$$  E_\t=\Bigl\{ x\in \SS^{d-1}:\ \f 1{4}\ge |x_2|\ge 2\sqrt{c_\a} |x_1|\ge 4\sqrt{c_\a}
|\t|,\ \ \text{and}\ \ |x_3|\ge \f 12\Bigr\}, $$ we assert  that
for $\a\neq 0$,
\begin{equation}\label{0-2-30}
|g_x''(t)| \ge c |x_2|^{2\a-2}, \   \  \text{whenever  $|t|\leq
2|\t|$ and $x\in E_\t$}, \end{equation} where $c$ is a  positive
constant depending only on $\a$. \eqref{0-2-30} together with the
mean value theorem will imply that for   $ \t\in (0, \d_\a]$
$$ |\tr_{1,3,\t}^2 f (x)|\ge C
|\t|^2  |x_2|^{2\a-2},\   \  \text{whenever $x\in E_\t$}.$$
Integrating the $p$-th power of the last inequality over $E_\t$ will give the
desired lower estimates.
To show the assertion \eqref{0-2-30}, we observe that if $x\in
E_\t$ and $|t|\leq 2|\t|$, then  $|v(t)|\leq |x_1|+|t|\leq 2|x_1|
\leq  |x_2|/\sqrt{c_\a}$, which implies
\begin{align*}
 4 |1-\a | (v(t)^2+x_2^2)^{\a-2} ( v(t)v'(t))^2
   & \leq 4 (1+|\a|)(v(t)^2+x_2^2)^{\a-2}
          \f { v(t)^2 +x_2^2}{c_\a+1} ( v'(t))^2\\
   &\leq ( v(t)^2+x_2^2)^{\a-1}(v'(t))^2.
\end{align*}
Thus, using \eqref{0-1-30}, we deduce that if  $x\in E_\t$ and
$|t|\leq 2|\t|\leq 2\d_\a$, then
\begin{align*}
|g_x''(t)| &\ge 2|\a| \left(v(t)^2+x_2^2\right)^{\a-1} \left [ (v'(t))^2 -
       |v(t) v''(t)|\right] - |\a|( v(t)^2 +x_2^2)^{\a-1} (v'(t))^2 \\
& = |\a| \left( v(t)^2 +x_2^2\right)^{\a-1}\left[ v'(t)^2 -2|v(t)v''(t)|\right]\\
& \ge c |\a| \left( v(t)^2 +x_2^2\right)^{\a-1}\left[ ( |x_3|-|t|)^2 - 2(|x_1|+|t|)^2 \right]\\
& \ge  c |\a| \left( v(t)^2 +x_2^2\right)^{\a-1}\sim (x_1^2+x_2^2)^{\a-1}\sim |x_2|^{2\a -2}
\end{align*}
proving the desired assertion  \eqref{0-2-30}.

For the upper estimate, it is easy to see by \eqref{0-1-30} that if
 $\sqrt{ x_1^2+x_2^2}\ge 4 |t|$ then
 $$
 |g_x''(t)| \leq c  ( v(t) ^2 +x_2^2)^{\a-1}\sim  (x_1^2+x_2^2)^{\a-1},
 $$
 which, using the mean value theorem,  implies that
 \begin{equation}\label{0-3}
   |\tr_{1,3,\t}^2 f (x)|\leq c |\t|^2 (x_1^2+x_2^2)^{\a-1}
 \end{equation}
whenever $\sqrt{ x_1^2+x_2^2}\ge 8  |\t|$.
On the other hand, however,  if $\sqrt{ x_1^2+x_2^2}\leq  8 \t$, then
using the  definition of $\tr_{i,j,\t}^2$, we have
\begin{align}\label{0-4}
|\tr_{1,3,\t}^2  f (x)| & \leq c |\t|^{2\a}+ c |x_2|^{2\a}.
\end{align}
Now we break the integral into two parts:
\begin{align*}
\int_{\SS^{d-1}}|\tr_{1,3,\t}^2 f(x)|^p d\s(x) & = \int_{ \{ x\in
\SS^{d-1}: x_1^2+x_2^2 \leq 64 \t^2\}}\cdots +\int_{ \{ x\in
\SS^{d-1}: x_1^2+x_2^2 > 64 \t^2\}}\cdots\\
 & \equiv I_1+I_2.
\end{align*}
Using \eqref{0-4} and the condition $\a p+1>0$, we
have
\begin{align*}
I_1& \leq c \int_{x_1^2+x_2^2 \leq  64 \t^2}(\t^{2\a p}
+|x_2|^{2\a p}) (1-x_1^2-x_2^2)^{\f {d-4}2}\, dx_1 dx_2\leq c |\t|^{2\a p+2},
\end{align*}
whereas using \eqref{0-3} gives
\begin{align*}
  I_2&\leq c |\t|^{2p} \int_{64\t^2\leq x_1^2+x_2^2 \leq 1 }
     |x_1^2+x_2^2|^{(\a -1)p}(1-x_1^2-x_2^2)^{\f d2 -2}\, dx_1 dx_2\\
& \leq c |\t|^{2p} \int_{8|\t|}^1 r^{(2\a -2)p+1} (1-r^2)^{\f d2-2}\, dr,
\end{align*}
which, by a simple calculation, leads to the desired upper
estimates.
\end{proof}

Our last example includes a family of functions and will be useful
in the next section. Note that the asymptotic orders in \eqref{0-7-10} and
\eqref{0-7-10-1} below are different for $\|y_0\| =1$ and $\|y_0\| < 1$,
as can be expected.

\begin{exam} \label{ex:sphere}
Let     $y_0$ be a fixed point in $\BB^d$, let
 $0\neq \a>-
\frac{d-1}{2p}$,  and let $f_\a : \SS^{d-1}\to \RR$ be given  by
$f_\a(x) := \|x - y_0 \|^{2 \a}$. If $\a\neq 1 - \frac{d-1}{2p}$
then
\begin{equation}\label{0-7-10}
    \omega_2(f_\a,t)_{L^p(\SS^{d-1})}\sim \|y_0\|t^2
        (t+1-\|y_0\|)^{2 (\a-1) + \frac{ d-1}{p}}+\|y_0\|t^2,
    \end{equation}where the constants of equivalence are independent of $t$ and
$y_0$. Moreover, if  $\a = 1 - \frac{d-1}{2p} $, then
\begin{equation}\label{0-7-10-1}
 c^{-1} \|y_0\|t^2\leq  \omega_2(f_\a,t)_{L^p(\SS^{d-1})}
\leq c\|y_0\| t^2|\log (t+1-\|y_0\|)
|^{\frac{1}{p}},\end{equation} where $c$ is a  positive constant
independent of $y_0$ and $t$.
\end{exam}

\begin{proof}
We start with the proof of the upper estimates
\begin{equation*}
\|\tr_{i,j, \t}^2 f _\a\|_{L^p(\SS^{d-1})}\leq c \Phi_\a(|\t|),\ \
1\leq i<j\leq d,\   \  |\t|\leq \d_\a,
\end{equation*}
where $\Phi_\a(t)$ denotes   the desired upper bound of
$\omega_2(f_\a,t)_{L^p(\SS^{d-1})}$  defined by the right hand
side of either  \eqref{0-7-10} or \eqref{0-7-10-1}. By symmetry,
it suffices to prove this inequality for $(i,j)=(1,2)$. We write $
y_0=(y_{0,1}, y_{0,2},\cdots, y_{0,d})=ry_0'$ with $0\leq r\leq 1$
and $y_0'\in \SS^{d-1}$. Let
$$
E_1:=\{ x\in \SS^{d-1}:\  \ \|x-y_0\|\leq 4 |\t|\} \,\,  \hbox{and} \,\,
E_2 :=\{ x\in \SS^{d-1}:\  \ \|x-y_0\|>  4 |\t|\}.
$$
We break the integral of $|\tr_{1,2,\t}^2 f_\a (x)|^p$ into two parts:
$$
\int_{\SS^{d-1}} |\tr_{1,2,\t}^2 f_\a (x)|^p\, d\s(x) =
   \int_{E_1}\cdots+\int_{E_2}\cdots=: I(E_1)+I(E_2).
$$
To estimate $I(E_1)$, we assume, without loss of generality, that
$r=\|y_0\|\ge 1-4|\t|>\f12$ since otherwise $E_1=\emptyset$. Then,
for $x\in E_1$, $\|x-y_0'\| \le \|x-y_0\| +(1-r) \le 8 \t$, which
is equivalent to $1- \la x,y_0'\ra \le 32 \t^2$ or $\arccos \la
x,y_0'\ra \le c \t$. Hence, we can deduce from the definition that
\begin{align*}
I(E_1) &\leq 4 \sup_{|t|\leq 2|\t|} \int_{E_1} |f_\a(Q_{1,2,t}
x)|^p\, d\s(x) \\
&  = 4 \sup_{|t|\leq 2|\t|} \int_{Q_{1,2,-t} (E_1)}|f_\a( x)|^p\, d\s(x)
\leq 4 \int_{ \|x-y_0'\|\leq 10|\t|} |f_\a (x) |^p\, d\s(x) \\
& = c \int_{ \|x-y_0'\|\leq  10 |\t|} \left | (1-r)^2 + 2r
(1- \la x,  y_0' \ra)\right|^{\a p}\, d\s(x) \\
 & \leq c \int_0^{c|\t|} \left [ |\t|^{2\a p} + (1-\cos v)^{\a p}\right]\sin^{d-2} v\, dv \\
 & \le c |\t|^{2\a p+d-1}\sim r^p |\t|^{2p}( 1-r+|\t|)^{2(\a-1)
 p+d-1},
\end{align*}
where the last step uses the assumption $1-r\leq 4|\t|<\f12$. To
estimate $I(E_2)$, we set $g_{x, y_0}(t)=f_\a(Q_{1,2,t} x)$ for a
fixed $x\in E_2$.  We then claim that
\begin{equation}\label{9-9-10-9}
|g_{x, y_0}''(t)|\leq c r\|x-y_0\|^{2\a-2},\quad \text{whenever
$|t|\leq 2|\t|$ and $x\in E_2$.}
\end{equation}
To see this, recall that $Q_{1,2,t} x= (x_1 (t), x_2(t),
x_3,\cdots, x_d)$, with $x_1(t) = x_1\cos t-x_2\sin t$ and $x_2(t)
=x_1\sin t+x_2\cos t$. Thus, a straightforward calculation shows
that
\begin{align}
 g_{x,y_0}''(t)  = & 4\a(\a-1) \|Q_{1,2,t} x -y_0\|^{2\a -4} \left (
  (x_1(t)-y_{0,1}) x_1'(t) + (x_2(t)-y_{0,2}) x_2'(t)\right)^2\notag \\
&+ 2\a\|Q_{1,2,t} x -y_0\|^{2\a -2}\left [ ( x_1'(t))^2 + (x_1(t)
   -y_{0,1})x_1''(t) + (x_2'(t))^2  \right. \notag\\
& \left.  + (x_2(t)-y_{0,2}) x_2''(t)\right].\label{9-10-10}
\end{align}
For  $x\in E_2$ and $|t|\leq 2|\t|$,  we have  $\|Q_{1,2,t}x - x
\| = 2 \sqrt{x_1^2+x_2^2} |\sin \f{t}2| \le 2 |\t| \le
\|x-y_0\|/2$, and  hence,  $\|Q_{1,2,t} x -y_0\|\sim \|x-y_0\|$.
Therefore, using \eqref{9-10-10}, we conclude that
$$ |g_{x, y_0}''(t)|\leq c \|x-y_0\|^{2\a-2},\quad
\text{whenever $|t|\leq 2|\t|$ and $x\in E_2$,}$$ which  proves
the claim  \eqref{9-9-10-9} when $r=\|y_0\|\ge 1/2$. On the other
hand, however, since the function $(x, y, t)\mapsto g_{x, y}''(t)$
is continuously differentiable on $\{ (x, y, t):\ x\in \sph, \
\|y\|\leq 1/2,  \  |t|\leq \pi\}$, it follows that for $|t|\leq
\pi$,\  \  $\|y_0\|\leq \d_\a$ and $ x\in \sph$,
$$ |g_{x, y_0}''(t)|=|g_{x, y_0}''(t)-g_{x, 0}''(t)|\leq c
\|y_0\|=cr,$$ where the second step uses the fact that $g_{x,
0}(t)\equiv 1$ for $x\in\sph$. This proves the claim
\eqref{9-9-10-9} for the case of $r<1/2$.

Now using \eqref{9-9-10-9} and  the mean value theorem, we have,
for some $\xi_\t$ between $0$ and $2\t$,
$$
\left|\tr_{1,2,\t}^2 f (x) \right |=\tfrac 12 \t^2
|g_{x,y_0}''(\xi_\t)|\leq cr \t^2
     \|x-y_0\|^{2\a-2}.
$$
Integrating the $p$-th power of the last inequality over $E_2$ gives
\begin{align*}
I(E_2) &\leq c r^p |\t|^{2p} \int_{E_2}\left ( ( 1-r)^2 +2r (
1- \la x,  y_0' \ra)\right)^{(\a-1)p}\, d\s(x).
\end{align*}
Thus, if $1-r \leq |\t|$, then $3|\t|\leq \|x-y_0'\| \le \arccos
\la x,y_0'\ra$, so that
\begin{align*}
I(E_2)&\leq c r^p |\t|^{2p} \int_{\{ x\in \sph: \arccos \la x,
y_0' \ra \ge 3|\t|\}  }\left( ( 1-r)^2 +2r ( 1- \la x, y_0'\ra\right)^{(\a-1)p}\,
    d\s(x)\\
&\sim r^p |\t|^{2p} \int_{3|\t|}^\pi ( 1-r +
s)^{2(\a-1)p}s^{d-2}\,ds \sim (\Phi_\a(|\t|))^p;
\end{align*}
whereas if  $1-r \ge  |\t|$, then $1-r+|\t| \sim 1-r$ and considering $r \ge 1/2$ and
$r \le 1/2$, respectively, we get
\begin{align*}
 I(E_2)&\leq c r^p
|\t|^{2p} \int_{\sph }\left( ( 1-r)^2 +
( 1- \la x, y_0'\ra\right)^{(\a-1)p}\, d\s(x)+c r^p|\t|^{2p}\\
&\sim r^p |\t|^{2p}\int_0^\pi (1-r+t)^{2(\a-1)p} \sin^{d-2} t\,dt+
  r^p |\t|^{2p}\sim (\Phi_\a(|\t|))^p.
\end{align*}
Putting the above together, we complete the proof of the upper estimates.

Next, we turn to the proof of the lower estimates. Without loss of
generality, we may assume that $|y_{0,1}|\ge
\|y_0\|/\sqrt{d}=r/\sqrt{d}$. We then claim that  if $|\t|\leq
\d_\a$ and $\a\neq 1-\f{d-1}{2p}$,
\begin{equation}\label{9-11}
\|\tr^2_{1, j, \t} f_\a\|_{L^p(\sph)}\ge c  \Phi_\a(|\t|) \quad
    \text{ for $2\leq j \leq d$}
\end{equation}
from which the desired lower estimate will follow. By symmetry, it
is enough to consider  $\tr^2_{1, 2, \t} f_\a$.   For $d > 3$,
using the formula \eqref{Int_m=2},  we can write
\begin{align}\label{0-9-jan}
  &\int_{\SS^{d-1}} |\tr_{1,2,\t}^2 f_\a (x)|^p d\sigma(x) \\
    &= \int_0^{\f \pi2} \cos \b (\sin \b)^{d-3}
     \int_{\SS^{d-3}} \int_{0}^{2\pi}
    \left |\overrightarrow{\tr}_{\t}^2 g_{\b,\xi}(\phi + \{\cdot\})
    \right |^p  d\phi  d\sigma (\xi) d\b,\notag
\end{align}
where  $g_{\b,\xi} (\phi) = f_\a(\cos\b\cos\phi, \cos\b\sin \phi,
\xi\sin\b )$ for $\xi \in \SS^{d-3}$, $\b \in [0,\f\pi2]$ and
$\phi \in [0, 2\pi]$. For $d = 3$, we need to use
\eqref{IntS-Bm=1}, which is an easier case. We shall assume $d >
3$ in the rest of the proof.

We write $y_0 =ry_0'=r (\cos\g \cos \phi_0, \cos\g \sin \phi_0,
v\sin\g)$, where $v\in \SS^{d-3}$ and $0\leq\g\leq \d_d<\f\pi2$,
the latter follows from $|y_{0,1}| \ge r /\sqrt{d} \ge (1-\d) /\sqrt{d}$.
 Then
\begin{align*}
  g_{\b,\xi}(\phi) & = \left(1+ r^2 - 2r\cos\g \cos\b \cos (\phi-\phi_0) -
     2 r\sin\b \sin\g\la \xi, v\ra \right)^\a \\
        & = A^\a (1 - a \cos (\phi-\phi_0))^\a,
\end{align*}
where $A: =1+ r^2  -  2 r(\sin \b\sin\g)\la \xi, v\ra$ and $a =2
r\cos\g\cos\b/A$. Since $0\leq \sin \b \sin \g \leq \sin \d_d<1$,
we have
$$
A= (1-r)^2 +2r(1-(\sin\b\sin\g)\la \xi, v\ra)\sim 1,
$$
and
\begin{align*}
1-a &=A^{-1} \left( (1-r)^2 +4r \sin^2\f {\b-\g}2 +2r \sin\b\sin\g
   (1-\la \xi,v\ra)\right) \\
&  \sim (1-r)^2 +r|\b-\g|^2 +r\sin\b\sin\g (1-\la \xi,v\ra)\\
&\sim (1-r)^2 +|\b-\g|^2 +\b\g (1-\la \xi,v\ra).
\end{align*}
In particular, there exists a constant $c_1\ge 2$ so that
$1-r+|\b-\g|\ge c_1 |\t|$ implies $1-a\ge 36 \t^2$.

If  $0\leq r < 1-\d <1$ for some small positive absolute constant
$\d$, then applying \mbox{Lemma \ref{lem-0-1}}, we deduce that for
$|\t|\leq \d/c_1$,
\begin{align*}
 &\int_{0}^{2\pi}
 \left |\overrightarrow{\tr}_{\t}^2 g_{\b,\xi}(\phi + \{\cdot\})\right |^p  d\phi
      \ge c  |\t|^{2p} A^{\a p} |a|^p   (1- a )^{(\a-1)p+1/2}  \\
      & \qquad \ge c  r^p |\t|^{2p}  (\cos\b)^p ,
\end{align*}
which, using  \eqref{0-9-jan},  implies the desired lower estimate
\eqref{9-11} in the case of $0\leq r \leq 1-\d$.

For the rest of the proof,   assume that   $1-\d<r\leq 1$. We
shall further assume that  $\a< 1-\f{d-1}{2p}$, as the case
$\a>1-\f{d-1}{2p}$ can be treated similarly, and in fact, is much
easier.  Since  $1-\f {d-1}{2p}<1-\f1{2p}$, we can apply
\mbox{Lemma \ref{lem-0-1}} to  deduce that for $\b\ge\g+ c_1
|\t|+1-r$,
\begin{align*}
 &\int_{0}^{2\pi}
 \left |\overrightarrow{\tr}_{\t}^2 g_{s,\xi}(\phi + \{\cdot\})\right |^p  d\phi
      \ge c  |\t|^{2p} A^{\a p} |a|^p   (1- a )^{(\a-1)p+1/2}  \\
      & \qquad \sim |\t|^{2p}  (\cos\b)^p \left[
      (1-r)^2 +|\b-\g|^2 +\g \b(1-\la \xi,v\ra)\right]^{(\a-1)p+1/2},
\end{align*}
where the last step uses the facts that $A\sim 1$, $\f12\leq r\leq
1$ and $0\leq \g\leq \d_d<\f\pi2$. Moreover, we can further assume
that $\b<\b_0<\f\pi2$ since $\g\leq \d_d<\f \pi2$.
 Thus, applying \eqref{0-9-jan}
and the formula
$$
 \int_{\SS^{d-3}} \Psi(\la \xi, u \ra) d\sigma(\xi) = \sigma_{d-4}
  \int_{-1}^1 \Psi(\|u\| z)
         (1-z^2)^{\frac{d-5}{2}}dz\, d\b,
$$
we conclude that
\begin{align}\label{0-10-jan}
  I:=& \int_{\SS^{d-1}} |\tr_{1,2,\t}^2 f_\a (x)|^p d\sigma(x)
  \ge c   |\t|^{2p} \int_{\g+c_1|\t|+1-r}^{\b_0}  \b ^{d-3} .\notag \\
  &\times \left[\int_{-1}^1 (1-z^2)^{\frac{d-5}{2}}  \left( (1-r)^2+ |\b - \g|^2 +
         \g\b  (1-z)  \right )^{(\a-1)p+1/2}dz\right]   d\b
\end{align}
To estimate $I$, we consider two cases.

{\it Case 1.} $0\leq\g \leq 1-r +c_1|\t|$. Then, for $\b$ in the integral,
$\g\b  (1-z)\leq \pi( 1-r + c_1 |\t|) \le \b - \g$, and hence
\begin{align*}
I & \ge c |\t|^{2p} \int_{\g+c_1|\t|+1-r}^{\b_0}  \b ^{d-3} (\b -\g)^{2(\a-1)p+1}\,  d\b\\
  & \ge c |\t|^{2p} \int_{2c_1|\t|+2(1-r)}^{\b_0}   \b^{2(\a-1)p+d-2}\, d\b
  \ge c \,\Phi_\a(|\t|)^p.
\end{align*}

{\it Case 2. }  $\g >1-r+c_1|\t|$. If $\b\leq \min\{3\g, \b_0\}=:\b_1$ then
$\f {(\b-\g)^2}{4\g\b} \leq \f {\b^2}{ 4 \g \b} \leq \f 34<1.$
Thus, we can restrict the domain of the integral to
$$
\g+c_1|\t|+1-r\leq \b\leq \b_1, \quad
    \f {(\b-\g)^2}{4\g\b}\leq 1-z \leq 1,
$$
and then obtain from \eqref{0-10-jan} that
\begin{align*}
I &\ge c |\t|^{2p} \int_{\g +c_1|\t|+1-r} ^{\b_1} \b^{d-3} (\g \b)^{(\a-1)p+\f12}
   \int_{\f {(\b-\g)^2}{4\g \b} } ^1 z^{\f {d-4}2   +(\a-1)p}\, dz\, d\b\\
& \sim |\t|^{2p} \int_{\g+c_1|\t|+1-r}^{\b_1} ( \b/\g)^{ - \frac{d-3}{2}}
(\b-\g)^{d-2+2(\a-1)p}\, d\b\\
&  \ge |\t|^{2p} \int_{\g+c_1|\t|+1-r}^{\b_1} (\b-\g)^{d-2+2(\a-1)p}\, d\b
  \sim |\t|^{2p} ( |\t|+1-r)^{2(\a-1)p+d-1},
\end{align*}
provided that $d-1+2(\a-1)p<0$.
\end{proof}

\subsection{Examples of best approximation on the sphere}
Our computational examples, together with  Theorem \ref{thm:best}
and its corollary,  immediately leads to the following examples on
the asymptotic order of $E_n(f)_{L^p(\SS^{d-1})}$.

\begin{exam} \label{example3}
For $d \ge 3$ let $f_\a(x) = x^\a$ with $\a=(\a_1,\ldots, \a_d)
\ne 0$. If $0 \leq  \a_i < 1$ for $1 \le i \le d$, then for
$n\in\NN$,
$$
     E_n(f_\a)_{L^p(\SS^{d-1})} \sim n^{- \d -1/p}, \qquad \d
                    = \min_{\a_i \ne 0} \{\a_1,\ldots, \a_d\}.
$$
\end{exam}

\begin{exam} \label{example4}
For $d \ge 3$ let $g_\a(x) = (1-x_1)^\a$, $x = (x_1,\ldots, x_d)
\in \sph$. Then for $ - \frac{d-1}{2p} < \a < 1 - \frac{d-1}{2p}$
and $\a\neq 0$ ,
\begin{equation*}
   E_n(f)_{L^p(\SS^{d-1})} \sim n^{-2 \a - \frac{d-1}{p}}, \qquad 1 \le p \le \infty.
\end{equation*}
\end{exam}

It is interesting to compare the two examples. As functions
defined on $\RR^d$, the functions $x_1^\a$ and $(1-x_1)^\a$ have
the same smoothness and a reasonable modulus of smoothness would
confirm that. As functions on the sphere $\sph$, however, they
have different orders of smoothness as seen in Examples
\ref{example1} and \ref{example2}, and their errors of best
approximation are also different as seen in Examples
\ref{example3} and \ref{example4}.

For $\a \ge 1 - \frac{d-1}{2p}$, the asymptotic order of
$\o_2(g_\a, t)_p$ in \eqref{eq:exam2} does not lead to the
asymptotic order of $E_n(f)_p$, since our inverse theorem in
\eqref{inverse} is of weak type. This remark also applies to other
examples below.

From our Examples \ref{ex:sphere-2} and \ref{ex:sphere}, we also
obtain the following results:

\begin{exam} \label{example5}
For $d \ge 3$ let $f_\a(x) = (x_1^2+ x_2^2)^\a$, $x = (x_1,\ldots,
x_d) \in \sph$. Then for $ - \frac{1}{p} < \a < 1 - \frac{1}{p}$
and $\a\neq 0$,
\begin{equation*}
   E_n(f)_{L^p(\SS^{d-1})} \sim n^{-2 \a - \frac{2}{p}}, \qquad 1 \le p \le \infty.
\end{equation*}
\end{exam}

\begin{exam} \label{example6}
Let $y_0$ be a fixed point in $\BB^d$, let $\a\neq 0$,  and let
$f_\a : \SS^{d-1}\to \RR$ be defined by $f_\a(x) := \|x - y_0
\|^{2 \a}$. If $ - \frac{d-1}{2p}<\a< 1 - \frac{d-1}{2p}$, then
\begin{equation*}
   E_n(f)_{L^p(\SS^{d-1})} \sim n^{-2}\|y_0\|
        (n^{-1}+1-\|y_0\|)^{2 (\a-1) + \frac{ d-1}{p}}.
\end{equation*}
\end{exam}

In particular, if $\|y_0\| =1$, then $f_\a$ has a singularity and
the asymptotic order is $n^{-2\a -\frac{ d-1}{p}}$ instead of
$n^{-2\a}$.

\section{Computational Examples on the unit ball}
\setcounter{equation}{0}

In this section we compute the modulus of smoothness $\o_r(f,t)_{L^p(\BB^d)}:
 =
\o_r(f,t)_{p, 1/2}$ defined in \eqref{omegaB-def} and the best
approximation $E_n(f)_{L^p(\BB^d)}:=E_n(f)_{p,1/2}$ of
\eqref{best-Ball}, both with constant weight function.

\subsection{Computation of moduli of smoothness}
Since $\o_r(f,t)_{L^p(\BB^d)}$ is closely related to
$\o_r(f,t)_{L^p(\SS^d)}$ according to Lemma \ref{lem:moduliBS},
our first three examples are derived directly from those in the
previous section.

\begin{exam} \label{ex-ball-1}
For $\a\neq 0$, define $f_\a: \BB^d\to \RR$ by  $f_\a(x) = (1-
\|x\|^2 + \|x - y_0\|^2)^\a$, where $y_0$ is a fixed point on
$\BB^d$. If $\a\neq 1-\f {d+1}{2p}$, then
$$
\omega_2(f_\a,t)_{L^p(\BB^d)} \sim t^2\|y_0\| ( t+1-\|y_0\|)^{
2(\a-1)
   +\f {d+1}p} +t^2\|y_0\|,
$$
where the constants of equivalence are independent of $y_0$ and
$t$. Moreover, if $\a=1-\f{d+1}{2p}$, then
$$
  c_\a^{-1} t^2 \|y_0\|\leq
   \omega_2(f_\a,t)_{L^p(\BB^d)} \leq c_\a t^2\|y_0\|
    |\log( t+1-\|y_0\|)|^{\frac{1}{p}},
$$
where $c_\a$ is independent of $t$ and $y_0$.
\end{exam}

\begin{proof}
Let $F_\a: \SS^{d+1}\to \RR$ be defined by
$$
    F_\a(x,x_{d+1}, x_{d+2}) = f_\a(x)= (\|x-y_0\|^2 + x_{d+1}^2+x_{d+2}^2)^\a =
    \|X-Y_0\|^{2\a},
$$
where $X=(x, x_{d+1}, x_{d+2})\in \SS^{d+1}$, $x\in \BB^d$, and $Y_0=(y_0, 0,0)
\in \BB^{d+2}$. Since the moduli of smoothness of $ F_\a$ on $\SS^{d+1}$ were
computed in Example \ref{ex:sphere}, the stated result follows
from Lemma \ref{lem:moduliBS}.
\end{proof}

Similarly,we can deduce directly  from Examples  \ref{example1} and \ref{ex:sphere-2}
the following results:

\begin{exam}  \label{ex-ball-2}
For $\a\neq 0$, let $f_\a(x)=(1-\|x\|^2)^\a$ for $x\in \BB^d$. Then
$$
 \omega_2(f_\a,t)_{L^p(\BB^d)} \sim\begin{cases}
t^{2\a +\f 2p},&\   \  \text{if $-\f1p<\a<1-\f 1p$;}\\
 t^2 |\log t|^{\f 1p}, & \   \  \text{if $\a=1-\f 1p $;}\\
 t^2, & \   \  \text{if $\a>1-\f 1p$.}
 \end{cases}
 $$
\end{exam}

\begin{exam}  \label{ex-ball-3}
Let $f_\a(x)=x^\a$ for $x\in \BB^d$ and  $\a=(\a_1,\cdots,
\a_d)\neq  0$. If $0\leq \a_i< 1$ for all $1\leq i\leq d$ then for
$1\leq p\leq \infty$,
$$
 \omega_2(f_\a,t)_{L^p(\BB^d)} \sim t^{\d +\f 1p},\   \
 \d:=\min_{\a_i\neq 0}\{ \a_1,\cdots, \a_d\}.
 $$
\end{exam}

Our next example is more complicated and requires a proof.

\begin{exam} \label{ex-ball-4}
Let $\a\neq 0$, $d\ge 2$ and let $f_\a: \BB^d\to \RR$ be given by
$f_\a(x) = \|x - e_0\|^{2\a}$, where   $e_0 = (1,0,\ldots,0)\in
\BB^d$. Then
\begin{equation}\label{3-15-10}
   \omega_2(f_\a,t)_{L^p(\BB^d)} \sim \begin{cases}
        t^{2 \a + \frac{ d}{p}}, & - \frac{d}{2 p} < \a < 1 - \frac{d}{2 p},\\
        t^{2} |\log t |^{\frac{1}{p}}, &  \a = 1 - \frac{d}{2p}, \quad p \ne \infty,\\
        t^{2},  &  \a > 1 - \frac{d}{2p}.
    \end{cases}
\end{equation}
\end{exam}

Before we give the proof of \eqref{3-15-10}, several remarks are
in order. First, it is interesting to compare these examples. We
consider  the function smoother when the asymptotic order of its
modulus of smoothness is higher.  Example \ref{ex-ball-2} has a
singularity at $\|x\| =1$, the boundary of $\BB^d$, and is a
radial function, for which the asymptotic order is independent of
the dimension $d$. Example \ref{ex-ball-4} has a singularity at $x
= e_0$, also on the boundary of the ball, but it is smoother than
the one in Example \ref{ex-ball-2} for $d > 2$ and $\a < 1-
\frac{d}{2p}$.  Furthermore, Example \ref{ex-ball-1} with $y =e_0$
also has a singularity at $x =e_0$ and its formulation is like the
addition of the other two cases; it is, nevertheless, the
smoothest one among the three functions. This does not seem to be
intuitively evident. Second, the comparison of these cases shows
the effect of the part of the modulus of smoothness in the Euler
angles. In fact, as the proof below will show, the part defined via 
difference in Euler angles in the definition of $\o_2(f,t)_{L^p(\BB^d)}$ 
in \eqref{omegaB-def} is dominating for Example \ref{ex-ball-4}. 
We also note that the asymptotic order of  Example
\ref{ex-ball-2} is independent of the dimension. Finally, we
should mention that the reason we restrict to $e_0$ in the last
example is given after the proof in Remark \ref{rem10.1}.

\medskip\noindent
{\it Proof of Example \ref{ex-ball-4}.}
The proof of \eqref{3-15-10} proceeds in  three steps. The first step deals with
the difference in the Euler angles, which can be done, in fact, more generally.

{\it Step 1.} Let $x_0 \in \SS^{d-1}$ and $f_\a(x) =
\|x-x_0\|^{2\a}$, which includes the case of $x_0 = e_0$ as a
special case. We prove that for $1\leq i<j\leq d$,
\begin{equation}\label{0-7-29}
   \|\tr^2_{i,j,\t} f_\a\|_{L^p(\BB^d)} \sim \begin{cases}
        |\t|^{2 \a + \frac{ d}{p}}, & - \frac{d}{2 p} < \a < 1 - \frac{d}{2 p},\\
        |\t|^{2} \bigl|\log |\t| \bigr|^{\frac{1}{p}}, &  \a = 1 - \frac{d}{2p}, \quad p \ne \infty,\\
        \t^{2},  &  \a > 1 - \frac{d}{2p}.
        \end{cases}
\end{equation}
For $x \in \BB^d$, write $x = \|x\| x'$, $x' \in \SS^{d-1}$. We then have
$$
   \|x - x_0\|^2 = 1+\|x\|^2 - 2 \|x\| \la x',x_0\ra =  \left \|x' - \|x\|x_0 \right\|^2.
$$
Let $g_{\a, r} (x') = \|x' - y_0\|^{2\a}$, $x' \in \SS^{d-1}$ and $y_0 = r x_0$. Then the
above equation shows
$$
   f_\a(x) =  g_{\a, \|x\|} (x'), \qquad x' \in \SS^{d-1}.
$$
Since $\tr_{i,j,\t}^r$ commutes with $\|x\|$,  by using polar coordinates and the proof
of Example 2.1, we obtain for $\a \ne 0, 1- \frac{d-1}{2p}$ that
\begin{align*}
   \int_{\BB^d} |\tr_{1,2,\t}^2 g_\a (x)|^p dx &= \int_0^1 r^{d-1} \int_{\SS^{d-1}}
      |\tr_{1,2,\t}^2 f_{\a,r} (x')|^p d\sigma(x') dr \\
     & \le  c  \int_0^1 r^{d-1+p} \left( |\t|^{2p}
     (1+ \t - r)^{2(\a-1)p +d-1} + |\t|^{2p} \right)  dr \\
     & \le c|\t|^{2\a p + d},
\end{align*}
if $2(\a-1)p + d < -1$ or $\a < 1- \frac{d}{2p}$; moreover, the same computation also
gives lower bound if we select a pair of $i,j$ as in Example 2.1, the choice of which
depends only on $x_0$ and is independent of $\|x\|$. The other cases of $\a$ can be
handled similarly.

{\it Step 2.} Next we consider the term $\tr_{1,d+1,\t}^2 \wt f_\a$. We show that  for
$\a>- \frac{d}{2 p}$
\begin{align}
 \left (\int_{\BB^{d+1}}
    |\tr^2_{1, d+1, \t} \wt f_\a(x)|^p \f {dx}{\sqrt{1-\|x\|^2}}\right)^{\f1p}\leq c
        |\t|^{2 \a + \frac{ d+1}{p}} + c \t^2,\label{3-17-10}
\end{align}
where $\wt f _\a(x) =f_\a (x')$ for $x=(x', x_{d+1})\in  \BB^{d+1}$. Let
$$
     E_{2,1}:= \{ x\in \BB^{d+1}: 1- c^2\t^2\leq x_1 \leq 1\}, \quad
     E_{2,2}:= \{ x\in \BB^{d+1}: -1\leq x_1 \leq 1-c^2\t^2\},
$$
where $c>1$ is a sufficiently large absolute constant; we break the integral
into two parts,
$$
\int_{\BB^{d+1}} |\tr^2_{1, d+1, \t} \wt f_\a(x)|^p \f {dx}{\sqrt{1-\|x\|^2}}
 =  \int_{E_{2,1}}\cdots+\int_{E_{2,2}}\cdots\equiv
I_2(E_{2,1})+I_2(E_{2,2}),
$$

To estimate $I_2(E_{2,1})$,  observe that for $x\in E_{2,1}$,
$|x_j|\leq \sqrt{1-x_1^2}\leq c |\t|$ for all $2\leq j \leq d+1$,
which implies, upon using $(1-x_1 \cos \psi - x_{d+1}\sin\psi) \le
c |\t|$  for $\psi \le 2 |\t|$, that $|f_\a(x_1 \cos \psi +
x_{d+1} \sin \psi, x_2,\ldots,x_d)| \le c
(|\t|^{2\a}+|x_2|^{2\a})$,  so that, by the definition,  $|\tr_{1,
d+1, \t}^2 \wt f_\a (x)|\leq
 c (|\t|^{2\a}+|x_2|^{2\a})$. Thus,
\begin{align*}
I_2(E_{2,1}) & \leq c \int_{1-c^2\t^2}^1 (1-x_1^2)^{\f
{d-1}2}\left[\int_{-1}^1 \left( |\t|^{2\a p} +  \left| \sqrt{1-x_1^2}
 s \right|^{2\a p} \right) (1-s^2)^{\f {d-2}2}\, ds \right] dx_1\\
& \leq c |\t|^{2\a p+d+1}.
\end{align*}
For the estimate of
$I_2(E_{2,2})$,  we write $ g_x(u)= \wt f_\a(Q_{1, d+1, u} x)$ for
a fixed $x\in \BB^{d+1}$.  That is,  $g_x (u)=\left( t_x(u)^2 +
\sum_{j=2}^d x_j^2\right)^\a $ with  $ t_x(u)=1-x_1\cos u +x_{d+1}
\sin u$. A straightforward calculation shows
\begin{align}\label{0-9-30}
 g_x''(u) =& \,4 \a (\a-1) \Bigl( t_x(u)^2 + \sum_{j=2}^d x_j^2\Bigr)^{\a-2}
 ( t_x(u) t_x'(u))^2\\
&+ 2\a\Bigl( t_x(u)^2 + \sum_{j=2}^d x_j^2\Bigr)^{\a-1} \left(
(t_x'(u))^2+ t_x(u) t''_x(u)\right)\notag.
\end{align}
Observe that if $x\in E_{2,2}$ and $|\xi|\leq 2|\t|$, then
$$ |t_x'(\xi)| = |x_1\sin \xi +x_{d+1}\cos \xi|\leq 2|\t|+|x_{d+1}|\leq
c' \sqrt{1-x_1},$$ and
\begin{align*}
 |t_x(\xi)-(1-x_1)|&\leq 2\t^2 + 2 \sqrt {1-x_1^2}
|\t|\leq ( \tfrac 2{c^2} + \tfrac {2\sqrt{2}}{c}) (1-x_1)
\end{align*}
which, in particular,  implies $ |t_x(\xi)|\sim 1-x_1$ provided
that  $c$ is large enough. Thus, by \eqref{0-9-30} it  follows
that for $x\in E_{2,2}$ and $ |\xi|\leq 2|\t|$,
\begin{equation*}
|g_x''(\xi)|\leq c \biggl( (x_1-1)^2 +\sum_{j=2}^d x_j^2 \biggr)^{\a-1} (1-x_1),
\end{equation*}
which,  using the mean value
theorem, implies that if $x\in E_{2,2}$, then for some $\xi$
between $0$ and $2\t$,
\begin{equation} \label{0-9-30-1}
   |\tr^2_{1, d+1, \t} \wt f_\a (x) |=\f 12\t^2 |g_x''(\xi)|\leq
     c \t^2 (1-x_1)\biggl( (1-x_1)^2 +\sum_{j=2}^d x_j^2 \biggr)^{\a-1}.
\end{equation}
If $\a \ge 1$, then we can drop the $(...)^{\a-1}$ term and the
estimate of $I_2(E_{2,2})$ follows trivially. For $\a < 1$,
integrating the $p$-th power of the inequality \eqref{0-9-30-1}
over $E_{2,2}$ yields
\begin{align*}
I_2(E_{2,2}) & \leq c |\t|^{2p} \int_{-1}^{1-c\t^2}(1-x_1^2)^{\a p + \frac{d-1}{2}}
 \int_{\BB^d}  \biggl( (1-x_1)^2 +\sum_{j=1}^{d-1} u_j^2 \biggr)^{(\a-1)p} \frac{du}{\sqrt{1-\|u\|^2}}
  d x_1 \\
& = c |\t|^{2p} \int_{-1}^{1-c\t^2}(1-x_1^2)^{\a p + \frac{d-1}{2}}
 \int_{\BB^{d-1}}  \left( (1-x_1)^2 + \|v\|^2 \right)^{(\a-1)p} dv dx_1  d x_1
\end{align*}
by \eqref{B-B}. Hence, switching to spherical-polar coordinates, we obtain
\begin{align*}
I_2(E_{2,2}) &  \leq c |\t|^{2p} \int_{-1}^{1-c\t^2}(1-x_1^2)^{\a p + \frac{d-1}{2}}
  \int_0^1 \left (\sqrt{1-x_1} + r \right)^{2(\a-1)p + d-2} dr dx_1\\
& \leq c  |\t|^{2p} \int_{-1}^{1-c\t^2} (1-x_1)^{(2\a-1) p + d-1} d x_1
 \le  c |\t|^{2\a p+d+1}+ c  |\t|^{ 2p}
\end{align*}
where we have used, in the last step, that $2  \a p + d \le 4 \a p + 2d$.

{\it Step 3.} Finally we consider $\tr_{i,d+1,\t}^2 \wt f_\a$ for $2 \le j \le d$.
We prove that for $\a>- \frac{d}{2 p}$,
\begin{align}\label{3-20-10}
 \left(\int_{\BB^{d+1}} |\tr^2_{j, d+1, \t} \wt f_\a(x)|^p \f
{dx}{\sqrt{1-\|x\|^2}}\right)^{\f1p}\leq c
        |\t|^{2 \a + \frac{ d+1}{p}}+c\t^2.
 \end{align}
As in the case $E_{2,2}$, the case of $\a \ge 1$ is easy and we assume
$\a < 1$. Clearly, it suffices to consider $\tr^2_{2, d+1, \t} \wt f(x)$. Let
\begin{align*}
     E_{3,1}&:=\{ x\in \BB^{d+1}: 1-x_1\leq c\t^2\},\\
     E_{3,2}&:=\{ x\in \BB^{d+1}: 1-x_1\ge  c\t^2,\
      |x_2|\ge 4\sqrt{1-x_1^2}|\t|\}\\
    E_{3,3} &:=\{ x\in \BB^{d+1}:1-x_1\ge  c\t^2,\ \  |x_2|< 4\sqrt{1-x_1^2}
      |\t|\},\end{align*}
where $c$ is a sufficiently large absolute constant.
We break the integral into three parts,
\begin{align*}
\int_{\BB^{d+1}} |\tr^2_{2, d+1, \t} \wt f_\a(x)|^p \f {dx}{\sqrt{1-\|x\|^2}}
 & = \int_{E_{3,1}}\cdots+\int_{E_{3,2}}\cdots+\int_{E_{3,3}}\cdots \\
  &  \equiv
I_3(E_{3,1})+I_3(E_{3,2})+I_3(E_{3,3}).
\end{align*}
Clearly, $I_3(E_{3,1})$ can be estimated exactly as $I_2(E_{2,1})$
in Step 2.

To estimate $I_3(E_{3,2})$ and $I_3(E_{3,3})$, we set,  as in Step
2, $g_x(u) = \wt f_\a(Q_{2, d+1, u} x)$ for any fixed $x\in \BB^{d+1}$.
Since  $|x_2|, |x_{d+1}|\leq \sqrt{1-x_1^2}$, it is easy to verify  that
\begin{equation}\label{0-11}
  |g_x''(u)|\leq c\biggl( (x_1-1)^2 +t_x(u)^2 +\sum_{j=3}^d
    x_j^2\biggr)^{\a-1} (1-x_1),
\end{equation}
where  $t_x(u)=x_2\cos u-x_{d+1} \sin u$. Observe that if $x\in E_{3,2}$
and $|u|\leq 2|\t|$ then
\begin{align*}
    |x_2-t_x(u)|&=|x_2(1-\cos u)+ x_{d+1}\sin u| \\
       & \leq \sqrt{x_2^2+x_{d+1}^2} \sqrt{(1-\cos u)^2 +\sin^2 u}
     \leq  \sqrt {1-x_1^2}|u|\leq \f12 |x_2|.
 \end{align*}
This  implies $|t_x(u)|\sim |x_2|$ for $x\in E_{3,2}$.
      Thus, using \eqref{0-11}, we have, for  $x\in E_{3,2}$,
\begin{equation}\label{3-22-jan}
   |\tr_{2, d+1,\t}^2 \wt f_\a(x)|=\f 12\t^2  |g_x''(\xi)|\leq
    c\t^2 \biggl( (x_1-1)^2+\sum_{j=2}^d x_j^2\biggr)^{\a-1}
(1-x_1),
\end{equation}
where $\xi$ is a number between $0$ and $2\t$. In particular, this allows
us to estimate $I_3(E_{3,2})$ exactly as  in Step 2.

It remains  to estimate  $I_3(E_{3,3})$. Using \eqref{0-11} and the mean
value theorem, we have, for $x\in  E_{3,3}$,
\begin{equation} \label{0-12}
|\tr_{2, d+1,\t}^2 \wt f_\a(x)|=\f 12\t^2  |g_x''(\xi)|\leq
      c\t^2 \biggl( (x_1-1)^2+\sum_{j=3}^d x_j^2\biggr)^{\a-1} (1-x_1),
\end{equation}
where $\xi$ is a number between $0$ and $2\t$. For $m \ge 2$ let
$E(m) := \{ (x_1,\ldots,x_m)\in \RR^m:  1-x_1\ge       c\t^2, |x_2|\leq 4 \sqrt{1-x_1^2}|\t|\}$.
Integrating the $p$-th power of \eqref{0-12} over $E_{3,3}$ and using \eqref{B-B},
we obtain
\begin{align*}
 I_3(E_{3,3})& \leq c |\t|^{2p} \int_{E(d)} (1-x_1)^p \Bigl( (1-x_1)^2
      +\sum_{j=3}^d x_j^2\Bigr)^{(\a-1)p}\, dx\\
&\leq c|\t|^{2p}\int_{E(2)} (1-x_1)^p \left[\int_{\BB^{d-2}} \left((1-x_1)^2
   +(1-x_1^2-x_2^2)\|v\|^2\right)^{(\a-1)p}\, dv\right] \\
& \qquad\qquad \qquad  \times (1-x_1^2-x_2^2)^{\f {d-2}2}\, dx_1dx_2+ c|\t|^{2p}.
\end{align*}
Note that if $0\leq x_1 \leq 1-c\t^2$ and $|x_2|\leq 4\sqrt{1-x_1^2} |\t|$
with $c\ge 32$,  then $|x_2|\leq 4|\t|$ and $1-x_1^2 \ge 1-x_1 \ge c\t^2\ge 2x_2^2$,
which implies $1-x_1^2-x_2^2\sim 1-x_1^2$. Thus,
\begin{align*}
 I_3(E_{3,3})&\leq c |\t|^{2p+1} \int_{0}^{1-c\t^2} (1-x_1)^{\a p+\f {d-1}2}
    \left[\int_{0}^1 \Bigl( \sqrt{1-x_1} +t\Bigr)^{2(\a-1)p +d-3}\,dt\right]\, dx_1\\
      & \qquad +c|\t|^{2p}\\
 & \leq c |\t|^{4\a p+2d}+c |\t|^p\leq c  |\t|^{2\a p+d+1}+ c |\t|^{2p}.
\end{align*}
Putting the above together, we have established \eqref{3-20-10}. The proof is complete.
\qed

\begin{rem} \label{rem10.1}
Our proof in Step 1 works for the more general case of $f_\a(x) =
\|x - x_0\|^{2\a}$, $x_0 \in \SS^{d-1}$. We notice that the steps
2 and 3 have smaller estimate, so that the dominating term is in
Step 1. We expect that \eqref{3-15-10} holds for $f_\a(x) =
\|x-x_0\|^{2\a}$. For the case of $\a < 1 - \frac{d}{2p}$, this is
indeed the case, as can be derived from our direct and the inverse
theorem, and the rotation invariance of $E_n(f)_{L^p(\BB^d)}$,
see \eqref{eq:ex-ball-6-1} at the end of the next subsection.
\end{rem}

\subsection{Examples of best approximation on the ball}
Our computational examples and Theorem \ref{thm:JacksonB}
immediately lead to the  following examples on the asymptotic
order of $E_n(f)_{L^p(\BB^d)}$. We give two examples, one
corresponds to Example \ref{ex-ball-1} and the other corresponds
to  Example \ref{ex-ball-4}.

\begin{exam} \label{ex-ball-5}
For $\a \ne 0$, let $f_\a(x) = (1-\|x\|^2)^\a$. Then for $-\frac{1}{p} < \a < 1 - \frac{1}{p}$,
$$
         E_n(f_\a)_{L^p(\BB^d)} \sim n^{- 2 \a - \frac 2 p}.
$$
\end{exam}

\begin{exam} \label{ex-ball-6}
For $\a \ne 0$, $d \ge 2$, let $f_\a(x) = \|x- e_0\|^{2\a}$, where
$e_0 = (1,0,\ldots,0)$. For $-\frac{d}{2p} < \a < 1 - \frac{d}{2
p}$,
\begin{equation}\label{eq:ex-ball-6}
     E_n(f_\a)_{L^p(\BB^d)} \sim n^{- 2 \a - \frac d p}.
\end{equation}
\end{exam}

Although our moduli of smoothness on the ball are not rotationally
invariant, the best approximation $E_n(f)_{L^p(\BB^d)}$ is; that
is, $E_n(f)_{L^p(\BB^d)} = E_n(f(\rho {\cdot}))_{L^p(\BB^d)}$ for
$\rho \in O(d)$. This implies, since every point $x_0$ on
$\SS^{d-1}$ can be rotated to $e_0$, that \eqref{eq:ex-ball-6}
holds for $f_{\a,x_0} (x) : = \|x - x_0\|^{2\a}$. In particular,
Theorem \ref{thm:JacksonB} shows then
\begin{equation} \label{eq:ex-ball-6-1}
   \o_2 (f_{\a,x_0},t)_{L^p(\BB^d)} \sim t^{2 \a  + \frac d p}, \qquad
          -\tfrac{d}{2p} < \a < 1 - \tfrac{d}{2 p},
\end{equation}
as we indicated in Remark \ref{rem10.1}.

\medskip\noindent
{\it Acknowledgement.} The authors thank a referee for his careful 
reading and helpful suggestions, especially for his suggestion of 
including explicit examples on the ball.

\enddocument